\documentclass[11pt]{article} 

\usepackage{jmlr2e}
\usepackage[utf8]{inputenc}
\usepackage{times}

\usepackage[a4paper]{hyperref}
\usepackage{xcolor}
\xdefinecolor{linkcol}{rgb}{0,0,0.4} 
\xdefinecolor{citecol}{rgb}{0.21,0.5,0.01} 
\xdefinecolor{redblack}{rgb}{0.6,0,0}
\xdefinecolor{myorange}{rgb}{0.15,0.15,0.15}
\xdefinecolor{myorange2}{rgb}{0.99,0.42,0.01}
\hypersetup
{
	bookmarksopen=true,
	pdftitle="Local DKW confidence bands",
	pdfauthor="Odalric-Ambrym Maillard", 
	pdfsubject="", 
	pdfmenubar=true, 
	pdfhighlight=/O, 
	colorlinks=true, 
	pdfpagemode=None, 
	pdfpagelayout=SinglePage, 
	pdffitwindow=true, 
	linkcolor=linkcol, 
	citecolor=citecol, 
	urlcolor=linkcol 
}

\usepackage{changepage}

\firstpageno{1}
\usepackage[top=1in, bottom=1.25in, left=1in, right=1in]{geometry}
\usepackage{fancyhdr}
\usepackage{appendix} 

\usepackage{graphicx}

\usepackage{algorithm,algorithmic} 

\usepackage{amsmath,amsfonts,amssymb}
\usepackage{math_letters} 

\newtheorem{assumption}{Assumption}

\newenvironment{mydefinition}[1]{
	
	\begin{minipage}{0.91\textwidth}
		\vspace{2mm}	
		\begin{definition}[#1]
		}{
		\end{definition}
		\vspace{0mm}
	\end{minipage}
	
}

\newenvironment{myassumption}[1]{
	
	\begin{minipage}{0.91\textwidth}
		\vspace{2mm}	
		\begin{assumption}[#1]
		}{
		\end{assumption}
		\vspace{0mm}
	\end{minipage}
	
}

\newenvironment{mylemma}[1]{
	
	\begin{minipage}{0.91\textwidth}
		\vspace{2mm}	
		\begin{lemma}[#1]
		}{
		\end{lemma}
		\vspace{0mm}
	\end{minipage}
	
}

\newenvironment{mytheorem}[1]{
	
	\begin{minipage}{0.91\textwidth}
		\vspace{2mm}	
		\begin{theorem}[#1]
		}{
		\end{theorem}
		\vspace{0mm}
	\end{minipage}
	
}

\newenvironment{mycorollary}[1]{
	
	\begin{minipage}{0.91\textwidth}
		\vspace{2mm}	
		\begin{corollary}[#1]
		}{
		\end{corollary}
		\vspace{0mm}
	\end{minipage}
	
}

\newenvironment{myproof}[1]{
	
	\begin{adjustwidth}{0.1cm}{0.1cm}	
		
		\noindent
		\hrulefill
		
		\noindent{\bf Proof #1: }				
		
		\vspace{-2mm}
		\noindent
		\hrulefill
		
	}
	{$\hfill\square$
		
		\vspace{-2mm}
		\noindent
		\hrulefill
	\end{adjustwidth}
	\bigskip
}

\newcommand{\CVAR}{\texttt{CVaR}}

\newcommand{\RTCE}{\underline{\texttt{RTCE}}}
\newcommand{\LTCE}{\overline{\texttt{LTCE}}}
\newcommand{\lCVAR}{\CVAR^\ell} 
\newcommand{\rCVAR}{\CVAR^r} 

\newcommand{\uu}{\underline{u}}
\newcommand{\ou}{\overline{u}}
\newcommand{\ux}{\underline{x}}
\newcommand{\ox}{\overline{x}}
\newcommand{\ldelta}{\delta} 
\newcommand{\rdelta}{\tilde \delta} 
\newcommand{\lepsilon}{\epsilon} 
\newcommand{\repsilon}{\tilde \epsilon} 

\renewcommand{\ln}{\log}

\usepackage{todonotes}

\newcommand{\includeg}[1]{
	\includegraphics[width=0.45\linewidth, trim = {0.9cm 0.8cm 1cm 0.8cm}, clip]{figs/#1}
}
\graphicspath{{figs/}}

\title{Local Dvoretzky-Kiefer-Wolfowitz confidence bands}
\author{
	\name{Odalric-Ambrym Maillard}\\
	\addr{Université de Lille, Inria, CNRS, Centrale Lille, UMR 9189 – CRIStAL, F-59000 Lille, France}\\
	\email{odalric.maillard@inria.fr}
}
\editor{}

\begin{document}
%
	\begin{center}
		{\Large {{\bf Local Dvoretzky-Kiefer-Wolfowitz confidence bands}}}

\bigskip
	{\bf Odalric-Ambrym Maillard}\\
	{\small 
		\textit{Université de Lille, Inria, CNRS, Centrale Lille}\\\textit{UMR 9189 – CRIStAL, F-59000 Lille, France}}\\
{\small\textsc{odalric.maillard@inria.fr}}
	\end{center}
\begin{abstract}
In this paper, we revisit the concentration inequalities for the supremum of the cumulative distribution function (CDF) of a real-valued continuous distribution as established by  Dvoretzky, Kiefer, Wolfowitz and revisited later by Massart in two seminal papers.
We focus on the concentration of the \textit{local} supremum over a sub-interval, rather than on the full domain. 
That is, denoting $U$ the CDF of the uniform distribution over $[0,1]$  and $U_n$ its empirical version  built from $n$ samples,  we study
$\Pr\Big(\sup_{u\in [\uu,\ou]} U_n(u)-U(u) >  \epsilon\Big)$ for different values of $\uu,\ou\in[0,1]$.
Such local controls naturally appear for instance when studying estimation error of spectral risk-measures (such as the conditional value at risk), where $[\uu,\ou]$ is typically $[0,\alpha]$ or $[1-\alpha,1]$ for a risk level $\alpha$, after reshaping the CDF $F$ of the considered distribution into $U$ by the general inverse transform $F^{-1}$.
Extending a proof technique from Smirnov, we provide exact expressions of the local quantities
$\Pr\Big(\sup_{u\in [\uu,\ou]} U_n(u)-U(u) >  \epsilon\Big)$ and $\Pr\Big(\sup_{u\in [\uu,\ou]} U(u)-U_n(u) >  \epsilon\Big)$
for each $n,\epsilon,\uu,\ou$. Interestingly these quantities, seen as a function of $\epsilon$, can  be easily inverted numerically into functions of the probability level $\delta$. Although not explicit, they can be computed and tabulated. We plot such expressions and compare them to the classical bound  $\sqrt{\frac{\ln(1/\delta)}{2n}}$ provided by Massart inequality. We then provide an application of such result to the control of generic functional of the CDF, motivated by the case of the conditional value at risk.
Last, we extend the local concentration results holding individually for each $n$ to time-uniform concentration inequalities holding simultaneously for all $n$, revisiting a reflection inequality by James, which is of independent interest for the study of sequential decision making strategies.
\end{abstract}


\section{Introduction}
Let $X$ be a real-valued random variable. The cumulative distribution function (CDF) $F$, that associates to each $x\!\in\! \Real$
the quantity $F(x)=\Pr(X\leq x)$ has been at the heart of statistics since its early ages, as $F$ characterizes the law of $X$.
The quantile function $Q(\delta)= \inf\{ x \!\in\!\Real : F(x) \geq \delta\}$ enables to simulate the random variable, since $X\stackrel{\cL}{=}Q(Y)$, where $Y$ is uniform on $[0,1]$. More generally one can reduce the study of the supremum of $F_n(x)-F(x)$ over a set $\cX$, where $F_n$ is the empirical version of $F$ built from $n$ i.i.d. samples of a distribution with CDF $F$ to the supremum of $U_n(u)-U(u)$ over its image $F(\cX)$, where $U$ denotes the CDF of the uniform distribution on $[0,1]$ and $U_n$ denotes its empirical version built from $n$ i.i.d. samples. 
The CDF is at the heart of the Glivenko-Cantelli theorem -- sometimes called the fundamental theorem of statistics --  that states that $\Pr(\lim_{n\to\infty} \sup_{x\in\Real}|F_n(x)-F(x)| =0)=1$. This result led to the definition of $P$-Glivenko-Cantelli classes of functions $\cF$,  for which $\displaystyle{\lim_{n\to\infty} \sup_{f\in\cF} |P_n(f)-P(f)|\to 0}$ almost surely, where $P(f)=\Esp[f(X)]$ denotes the measure associated to the random variable $X$, and $P_n$ denotes the empirical measure built from $n$ i.i.d. samples. This definition had a prominent role in the development of function process theory as Glivenko-Cantelli theorem shows that the class of functions $\cF = \{x\mapsto \indic{x\leq t}: t\!\in\!\Real\}$ is an example being  $P$-Glivenko-Cantelli for all probability measure $P$ on $\Real$, which opened the quest for other such classes.  
While Glivenko-Cantelli classes are nice, great efforts have been put on obtaining not only asymptotic results but further understand the speed of convergence of the supremum towards $0$. One important part of the literature focuses on Donsker classes where the supremum of $\sqrt{n}(F_n(x)-F(x))$ is studied as a random process (See \cite{kolmogorov1933sulla}, wrongly extended by \cite{donsker1952} but later corrected; we refer the interested reader to \cite{billingsley1968convergence},  \cite{pollard2012convergence}, \cite{dudley1999uniform} or \cite{shorack2009empirical} for further details and overview of the field).
An alternative approach to this prolific field of research is the one proposed by Dvoretzky–Kiefer–Wolfowitz in their seminal paper \cite{dvoretzky1956asymptotic} that looks at deviation inequalities of the supremum process for each $n$. The initial result from \cite{dvoretzky1956asymptotic} is based on an exact derivation from \cite{smirnov1944approximate}, and shows that
\beqan
\Pr\Big(\sup_{u\in [0,1]} U_n(u)-U(u) >  \epsilon\Big) \leq Ce^{-2n\epsilon^2}\,,
\eeqan 
hence providing an exponentially decreasing upper bound on the deviation probability, yet for some unspecified constant $C$.
In a seminal paper, \cite{massart1990tight}  later showed that the result holds for the constant $C=1$ provided that $e^{-2n\epsilon^2}\leq 1/2$, and that this constant cannot be improved. Such a result is especially interesting as it enables the practitioner to derive tight confidence bands on CDF. Indeed, it can be used to show that for any $\delta\in[0,0.5]$, with probability higher than $1-\delta$, uniformly for all $x\in\Real$, then 
\beqan
F_n(x) -\epsilon \leq F(x) \leq F_n(x) +\epsilon \quad\text{ where } \epsilon= \sqrt{\frac{\ln(2/\delta)}{2n}}\,. 
\eeqan

\paragraph{Risk-aversion estimation}
One important example of application of CDF deviation inequalities is when considering spectral risk-measures, such as the Conditional Value At Risk (CVAR) that is popular in economy (See \cite{mandelbrot1997variation,rockafellar2000optimization}). The  definition of the \CVAR\ changes from author to author, depending on conventions, such as whether it applies to a non-negative or non-positive random variable, and whether the risk corresponds to the upper or lower tail. We choose below a non-negative random variable with focus on its lower tail for this short presentation. While \CVAR\ at risk level $\alpha\in[0,1]$ is classically defined as an optimization (see Section~\ref{sec:risk} for more details), 
when the CDF $F$ of the considered random variable is a continuous bijection,
it takes the following convenient form
\beqan
\rCVAR_{\alpha}(\nu) &\!=\!& \Esp_\nu\bigg[X \bigg| X\leq \ox_{\alpha}(\nu) \bigg] = \frac{1}{\alpha}\Esp\Big[X\indic{X\!\leq\! \ox_{\alpha}(\nu)}\Big] 
=\frac{1}{\alpha}\int_{\Real^+} \!\!\max(\alpha\!-\!F(x),0)dx\,.
\eeqan
where we introduced the (upper) Value at Risk 
 $	\ox_{\alpha}(\nu)= \inf\{ x \!\in\!\Real : F(x)>\alpha\}$.
In particular the \CVAR\ writes as a function of the CDF in the form
$\rCVAR_{\alpha}(\nu) = h^r\Big(\int_{\cX}g^r\big(F(x)\big)dx\Big)$, where $h^r$ and $g^r$ are known monotonic functions. 
Further $g^r(\beta)=(\alpha-\beta)_+$ has support $[0,\alpha)$ that is a strict subset of $[0,1]$ for $\alpha<1$, and $\alpha$ is typically small (e.g. $\alpha=0.01, \alpha=0.05$).
This property is actually not limited to the CVaR risk-measure (see e.g. spectral risk-measures \cite{acerbi2002spectral}) and suggests to focus on controlling the local deviations of the CDF in order to later control the risk-measure.
Note also that since in this case
$\sup_{x\in [0,\ox_{\alpha}(\nu)]} F_n(x)-F(x) \stackrel{\cL}{=} \sup_{u\in [0,\alpha]} U_n(u)-U(u)$ and $\alpha$ is given by the practitioner,
controlling the local supremum for the uniform distribution  ensures a control for any unknown distribution. Hence, such results are universal in a sense, which is especially interesting when the class of distributions generating the samples is unknown to the practitioner. 

\paragraph{Outline and contribution}
The purpose of this work is not to greatly extend the vast literature on function process theory, but rather to focus specifically on the \textit{local} control of the CDF concentration. Namely, we ask how to revisit the initial results from Smirnov in the case when the supremum is not considered on $\Real$ but on a sub-interval $\cX$.
Obviously the results by DKW and Massart already apply to yield confidence bands in such cases. However, making use of the existing bound  to uniformly control  the deviations on a  set $\cX$ that is ``small''  may result in unnecessary large confidence bands that may be concerning for the practitioner.
In this article,  extending the proof techniques from \cite{smirnov1944approximate}, we derive in Section~\ref{sec:local}  exact expressions for the quantities
$\Pr\Big(\sup_{u\in [\uu,\ou]} U_n(u)-U(u) >  \epsilon\Big)$  (see Theorem~\ref{thm:concLeft}) and $\Pr\Big(\sup_{u\in [\uu,\ou]} U(u)-U_n(u) >  \epsilon\Big)$ (see Theorem~\ref{thm:concRight}) for any $n\!\in\!\Nat,\epsilon>0$ and $\uu,\ou \!\in\![0,1]$. This derivation, up to our knowledge, is (perhaps surprisingly) novel.
Interestingly, these probabilities, seen as functions of $\epsilon$ can be inverted numerically, and directly yield confidence bands on the CDF.  In Section~\ref{sec:xps}, we plot these functions and their inverse,
which enables to  highlight their non-trivial behavior and compare them to the simplified bound obtained by application of DKW and Massart bound.
This also shows that deriving approximations of the exact quantities may not be necessarily required for practical usage.
Although  these plots reveal the strikingly good match between the simplified and exact bound when considering the full interval $[0,1]$, which is expected since the constant $C=1$ obtained by Massart cannot be improved,
they also reveal the conservative nature of the simplified bounds when considering a supremum over a sub-interval of the form $[0,\alpha]$ or $[1-\alpha,1]$ for small values of $\alpha$. We believe making use of the exact bounds may thus greatly impact the practitioner interested in tight bounds.
To provide greater perspective, in Section~\ref{sec:risk}, we quickly detail the case of the risk-measure known as Conditional Value At Risk. 
In Section~ \ref{sec:TU}, we finally describe a generic way to turn the concentration bounds of Section~\ref{sec:local}, valid with high probability for each $n\!\in\!\Nat$,
into concentration bounds that are  time-uniform,  that is, valid with high probability simultaneously for all $n\!\in\!\Nat$. This extension is not trivial as it seems difficult to make use of Doob's maximal inequality in this context, hence we describe an alternative way to derive such results that is of independent interest.

\tableofcontents

\newpage
\section{Local CDF concentration}\label{sec:local}

For a distribution $\nu$ with CDF $F$, we denote by $F_n$ the empirical CDF built from a sample of size $n$ from $\nu$.
We consider distributions with \textit{continuous} CDF in the sequel.
Since the uniform distribution on $[0,1]$ plays a special role, we denote its CDF by $U$ and empirical CDF with a sample of size $n$ by $U_n$.

Let us first recall the result obtained by Massart's version of DKW inequality from \cite{massart1990tight}. 
\beqan
\forall \delta_0 \!\in\![0,0.5)\quad
\Pr\bigg( \sup_{x\in\cX} F(x)-F_n(x) > \sqrt{\frac{\ln(1/\delta_0)}{2n}}\bigg) \leq \delta_0\,.\\
\forall \delta_0 \!\in\![0,0.5)\quad
\Pr\bigg( \sup_{x\in\cX} F_n(x)-F(x) > \sqrt{\frac{\ln(1/\delta_0)}{2n}}\bigg) \leq \delta_0\,.
\eeqan

We now introduce the following closely related quantities on which we shall focus
\beqan
\ldelta_{[\uu,\ou]}(n,\epsilon)&=&\Pr\Big(\sup_{u\in [\uu,\ou]} U_n(u)-U(u) >  \epsilon\Big)\quad \lepsilon_{[\uu,\ou]}(n,\delta) = \inf\{ \epsilon : \ldelta_{[\uu,\ou]}(n,\epsilon) \leq \delta\}\\
\rdelta_{[\uu,\ou]}(n,\epsilon)&=&\Pr\Big(\sup_{u\in[\uu,\ou]} U(u)-U_n(u)  >  \epsilon\Big)\quad \repsilon_{[\uu,\ou]}(n,\delta) = \inf\{ \epsilon : \rdelta_{[\uu,\ou]}(n,\epsilon) \leq \delta\}\,.
\eeqan
Since in general, $U_n(u)\!\neq\!-U_n(1-u)$, then $\ldelta_{[\uu,\ou]}(n,\epsilon)$ is in general different from $\rdelta_{[\uu,\ou]}(n,\epsilon)$. Also, since
$U_n-U$ is right continuous with a left limit, but $U-U_n$ is left continuous with a right limit, the supremum of these functions have different behavior. 
In particular, the following result, whose proof is immediate given this observation, shows that considering the second supremum should be considered with care. 
\begin{lemma}[Asymmetry]	
	On the one hand, the value of the optimization problem $\sup_{u\in [\alpha,\beta]} U_n(u)-U(u)$ is achieved at one of the points $v=u_{(k)}$ or at $v=\alpha$, where $(u_{(k)})_{k\in[n]}$ are the order samples received from the uniform distribution, such that 
	$0 \leq  u_{(1)} \leq \dots \leq u_{(n)}\leq 1$.
	\beqan
	\sup_{u\in [\alpha,\beta]} U_n(u)-U(u) = \max\bigg\{U_n(v)-U(v): v \in \{\alpha\}\cup \{ u_{(1)},\dots,u_{(n)} \}\cap[\alpha,\beta]\bigg\}\,.
	\eeqan
		
	On the other hand, the  value of the optimization problem $\sup_{u\in [\alpha,\beta]} U(u)-U_n(u)$ is not realized at any point, but can be approached from below when approaching some of the $v=u_{(k)}$ or point $v=\beta$ from below.
	\beqan
	\sup_{u\in [\alpha,\beta]} U(u)-U_n(u) = \max\bigg\{\lim_{u\to v, u<v} U(u)-U_n(u): v \in  \{\beta\}\cup\{u_{(1)},\dots,u_{(n)}\}\cap(\alpha,\beta]\bigg\}
	\eeqan
\end{lemma}
To complement this lemma, let us recall that $\sup_{u\in [\alpha,\beta]} U_n(u)-U(u)$ and $\sup_{u\in [1-\beta,1-\alpha]} U(u)-U_n(u)$ should have the same law. Nevertheless, we now provide the following result that gives a first explicit expression of $\ldelta_{[\uu,\ou]}$ and
$\rdelta_{[\uu,\ou]}$.  For the second quantity, we make use of a construction inspired from Skorokhod convergence since the considered function is not \textit{corlol} (continuous on the right with a limit on the left).

\begin{figure}
	\centering
	\includegraphics[width=0.45\linewidth]{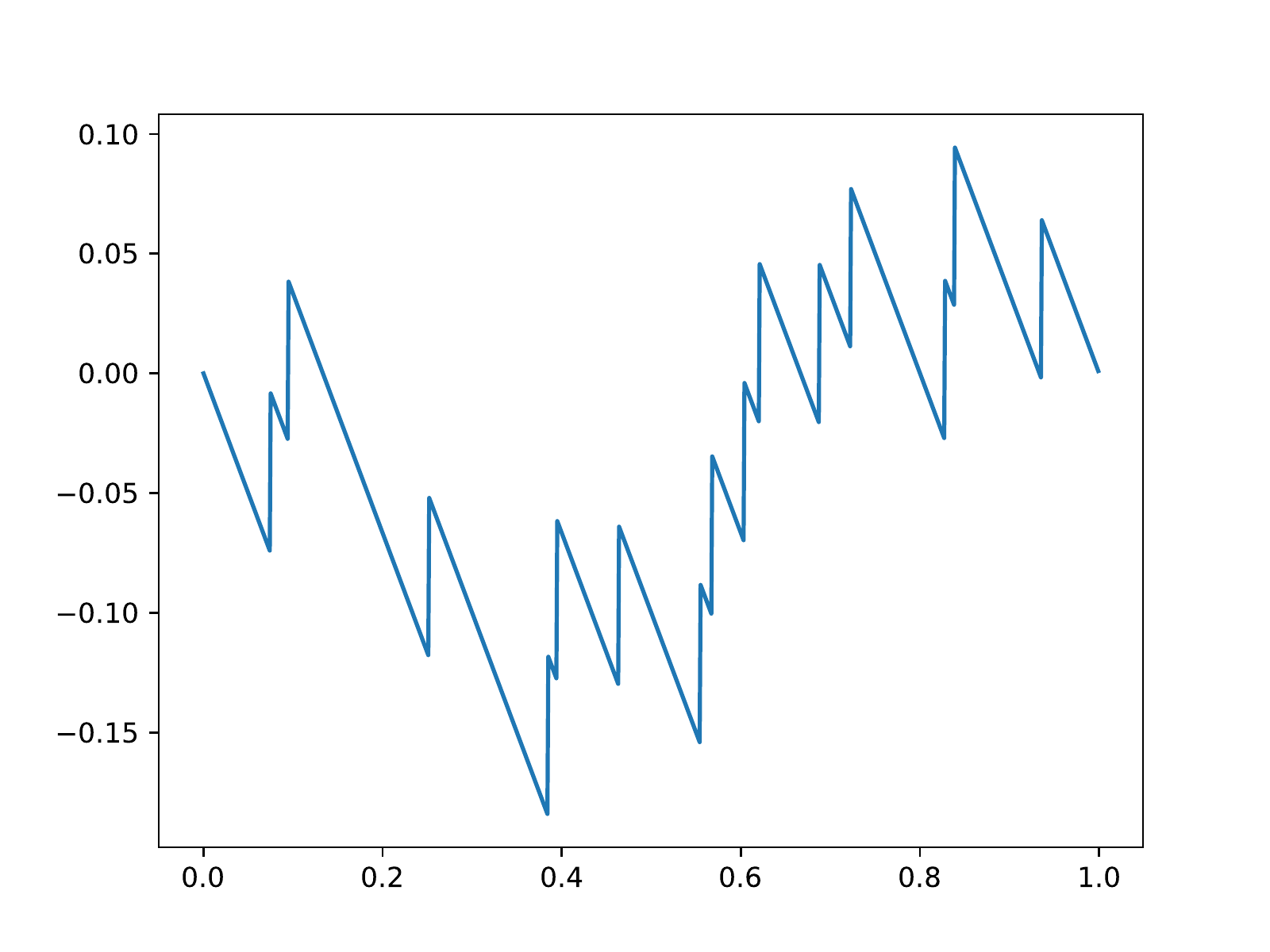}\hfill
	\includegraphics[width=0.45\linewidth]{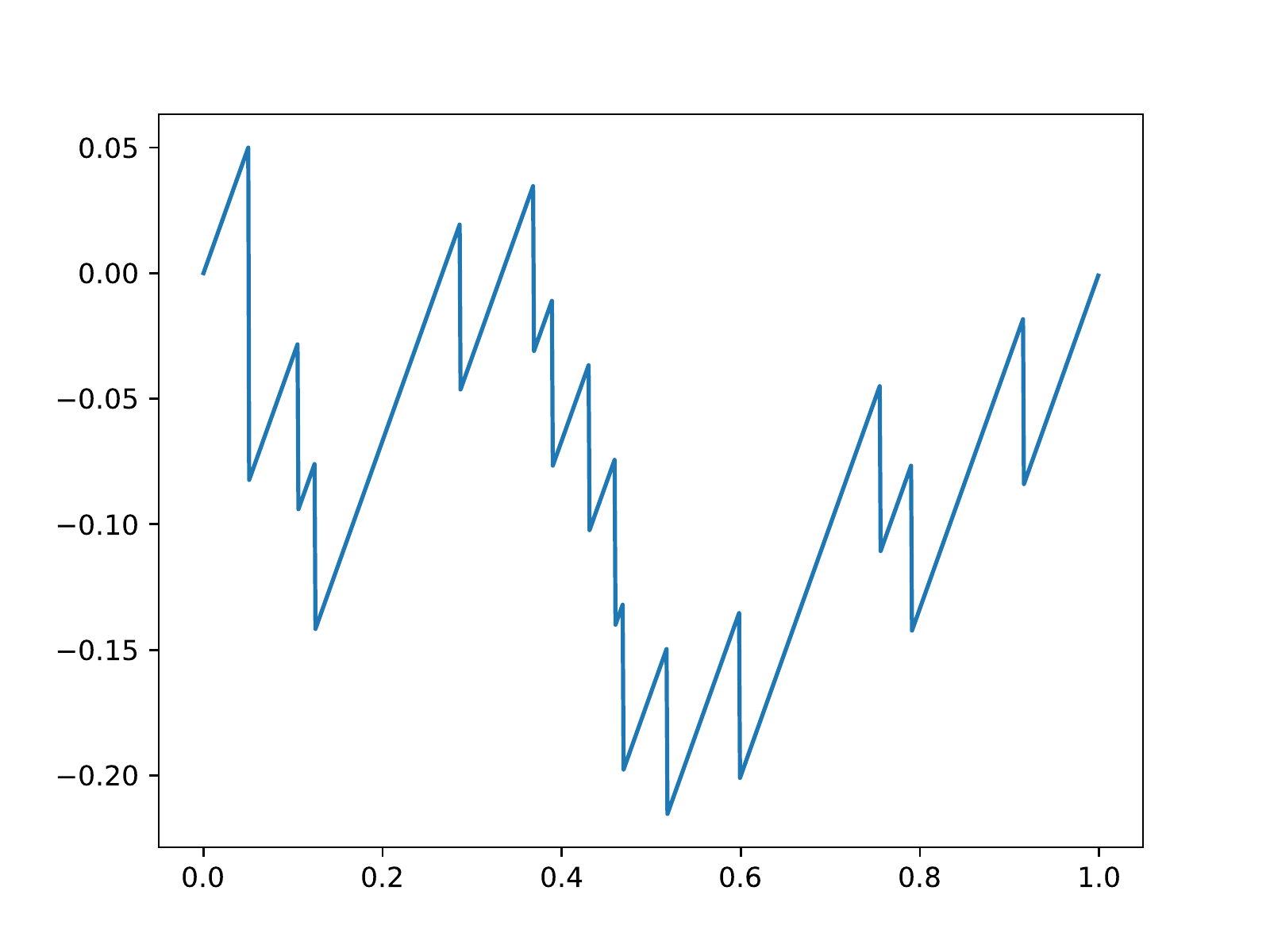}\hfill
	\caption{Illustration of the random function $u\mapsto U_n(u)-u$ (left) and 
	$u\mapsto u-U_n(u)$ (right) using $n=15$. }
	\label{fig:illustrateleft}
\end{figure}

\begin{lemma}[Exact CDF concentration]\label{lem:concfirst}Let $U$ denote the CDF of a random variable uniformly distributed on $[0,1]$ (that is, $U(x)=x$).	Let $U_n$ denote the empirical CDF built from $n$ i.i.d. samples from this distribution.
	Let us introduce the following notation	$I_k(x;a_1,\dots,a_k) = \int_{a_1}^{x}\int_{a_2}^{t_1}\dots\int_{a_k}^{t_{k-1}}dt_1\dots dt_k$
	for all $k\in\Nat_\star$ and $x \geq a_1\geq \dots\geq a_k\in\Real$, and $I_0(x;\emptyset)=1$.
	Let $\alpha,\beta$ be such that $[\alpha,\beta]\subset[0,1]$, $\epsilon>0$ and $n\in\Nat$. 
	Then it holds
	\beqan
	\Pr\Big(\sup_{u\in [\alpha,\beta]} U_n(u)-U(u)  > \epsilon\Big)&=&	\sum_{\ell=0}^{\overline{n}_{\alpha,\epsilon} -1} \binom{n}{\ell} \beta_{\ell+1,\epsilon}^{n-\ell}\ell!I_{\ell}(1;\beta_{1,\epsilon},\dots \beta_{\ell,\epsilon})\,,\\
\Pr\Big(\sup_{u\in [\alpha,\beta]} U(u)-U_n(u)  >  \epsilon\Big)
&=&\sum_{\ell=0}^{\overline{n}_{\beta,\epsilon}-1} \binom{n}{\ell}\tilde \alpha_{\ell+1,\epsilon}^{n-\ell}\ell!I_{\ell}(1;\tilde \alpha_{1,\epsilon},\dots \tilde \alpha_{\ell,\epsilon}) \,.
\eeqan
where  we introduced $\beta_{k,\epsilon} =\min(\beta,(n-k+1)/n-\epsilon)$,   $\tilde \beta=1-\beta$, $\tilde \alpha_{k,\epsilon} =\min(1-\alpha,1-\epsilon-(k-1)/n)$ as well as 
$\overline{n}_{\alpha,\epsilon} = \lceil n(1-\alpha-\epsilon) \rceil$ and $\overline{n}_{\beta,\epsilon} = \lceil n(\beta-\epsilon)\rceil$.

\end{lemma}

\paragraph{Sketch of proof of Lemma~\ref{lem:concfirst} }
We provide below a sketch of proof of the key steps leading to the control of $\sup_{u\in [\alpha,\beta]} U_n(u)-U(u)$ in Lemma~\ref{lem:concfirst}. 
The full proof is detailed in and deferred to Appendix~\ref{app:concentration}.
This enables to highlight the main ingredients of the proof, and in particular the smart use of a generic Taylor expansion by Smirnov, that we reuse in order to prove this novel result.

The first step of the proof consists in showing that
\beqan
\lefteqn{\Pr\Big(\sup_{u\in [\alpha,\beta]} U_n(u)-U(u) \leq \epsilon\Big)
}\\	 &=&n! \int\dots\int\ind\bigg\{ 0\leq  u_1 \leq \dots u_n \leq 1; \forall k\in [n],
\begin{cases}
	\text{if } u_k \in [\alpha,\beta]& \text{ then }u_k \geq \frac{k}{n}-\epsilon\\
	\text{if } k>n(\epsilon+\alpha) & \text{ then } \alpha \notin [u_k,u_{k+1})\\
\end{cases}\bigg\} du_1\dots du_n.
\eeqan

The second step enables to simplify the expression, leading, after careful rewriting, to
\beqan
\Pr\Big(\sup_{u\in [\alpha,\beta]} U_n(u)-U(u) \leq \epsilon\Big) 
&=&
n! \int_{\beta_1}^1\!\int_{\beta_2}^{t_1}
\!\dots  \!\int_{\beta_{\overline{n}_{\alpha,\epsilon}}}^{t_{\overline{n}_{\alpha,\epsilon}-1}}\!\!J_{n-\overline{n}_{\alpha,\epsilon}}(t_{\overline{n}_{\alpha,\epsilon}})dt_{\overline{n}_{\alpha,\epsilon}}\!\dots\! dt_1\!\,,
\eeqan
where $J_k(x)=\frac{x^k}{k!}$, $\beta_k = \min(\beta, (n-k+1)/n-\epsilon)$  and $\overline{n}_{\alpha,\epsilon}=\lceil n(1-\epsilon-\alpha)\rceil$.

At this point, the key trick, originating from \cite{smirnov1944approximate}, is to make use of the following  variant of the Taylor expansion 
at multiple points $a_1,a_2,\dots$, of the value of a smooth function $f$ at point $x$, given by
\beqan
f(x) &=& f(a_1) +  \sum_{\ell=1}^{n-k-1}f^{(\ell)}(a_{\ell+1})\int_{a_1}^{x}\int_{a_2}^{t_1}\dots\int_{a_k}^{t_{\ell-1}}dt_1\dots dt_\ell\\
&& + \int_{a_1}^{x}\int_{a_2}^{t_1}\dots\int_{a_{n-k}}^{t_{n-k-1}}f^{(n-k)}(t_{n-k})dt_1\dots dt_{n-k}\,,
\eeqan
applied to the function $f(x)=x^n$. This enables to rewrite the multiple integral in explicit terms. After some careful computations, this yields the desired result.

\paragraph{Explicit local CDF concentration}
In this section, we now detail the computations of the multiple integrals $I_\ell$ appearing in Lemma~\ref{lem:concfirst}.
Theorem~\ref{thm:concLeft} and Theorem~\ref{thm:concRight} below constitute the main results of this paper, as they provide the \textit{exact} value of $\ldelta_{[\uu,\ou]}(n,\epsilon)$ and $\rdelta_{[\uu,\ou]}(n,\epsilon)$. The full proof is deferred to Appendix~\ref{app:concentration}.

\begin{theorem}[Exact Left-CDF concentration]\label{thm:concLeft}Let $U$ denote the CDF of a random variable uniformly distributed on $[0,1]$ (that is, $U(x)=x$).	Let $U_n$ denote the empirical CDF built from $n$ i.i.d. samples from this distribution.
	Let $\alpha,\beta$ be such that $[\alpha,\beta]\subset[0,1]$, $\epsilon>0$ and $n\in\Nat$. Let
  $\overline{n}_{\alpha,\epsilon} = \lceil n(1-\alpha-\epsilon) \rceil$, then $n_\beta=n (1- \beta-\epsilon)$ and $m = \min\{\lfloor n_\beta\rfloor+1,\overline{n}_{\alpha,\epsilon}-1\}$. It holds when $n_\beta>0$,
	\beqan
	\lefteqn{
		\Pr\Big(\sup_{u\in [\alpha,\beta]} U_n(u)-U(u)  > \epsilon\Big)=	
		\sum_{\ell=0}^{m}\binom{n}{\ell}\Big(\!\min\Big\{1\!-\!\frac{\ell}{n}\!-\!\epsilon,\beta\Big\}\!\Big)^{n\!-\!\ell}(1-\beta)^\ell} \\
	&&+\!\!\sum_{\ell=m+1}^{\overline{n}_{\alpha,\epsilon}-1}\!\! 
	\binom{n}{\ell}\!\bigg(1\!-\! \frac{\ell}{n}\!-\!\epsilon\bigg)^{n\!-\!\ell}\bigg[\epsilon
	\bigg(\!\frac{\ell}{n}\!+\!\epsilon\!\bigg)^{\ell\!-\!1\!}\!\!+\! \sum_{j=0}^{m-1}\!\bigg[\frac{n_\beta\!-\!j}{n}\bigg]\binom{\ell}{j}\bigg(\frac{\ell\!-\!n_\beta}{n}\bigg)^{\ell\!-\!j\!-\!1\!}\!(1\!-\!\beta)^{j}\bigg]\,.
	\eeqan
	
	\vspace{-5mm}\noindent
	(Using the convention that $\displaystyle{\sum_{\ell=i}^j}$ is $0$ if $i>j$ for both sums.)
	If, on the other hand $n_\beta<0$, then

	\vspace{-6mm}\noindent
	\beqan
	\Pr\Big(\sup_{u\in [\alpha,\beta]} U_n(u)-U(u)  > \epsilon\Big)
	&=& \sum_{\ell=0}^{\overline{n}_{\alpha,\epsilon}-1} \!\binom{n}{\ell}\! \Big(1\!-\!\frac{\ell}{n}\!-\!\epsilon\Big)^{n\!-\!\ell}\!\epsilon\Big(\frac{\ell}{n}\!+\!\epsilon\Big)^{\ell-1}\,.
	\eeqan
\end{theorem}

\begin{theorem}[Exact Right-CDF concentration]\label{thm:concRight}Let $U$ denote the CDF of a random variable uniformly distributed on $[0,1]$ (that is, $U(x)=x$).	Let $U_n$ denote the empirical CDF built from $n$ i.i.d. samples from this distribution.
	Let $\alpha,\beta$ be such that $[\alpha,\beta]\subset[0,1]$, $\epsilon>0$ and $n\in\Nat$.
	Let $\overline{n}_{\beta,\epsilon}=\lceil n(\beta-\epsilon) \rceil$,
$n_\alpha=n(\alpha-\epsilon)$, 
and finally
$\tilde m = \min\{\lfloor n_\alpha\rfloor+1,\overline{n}_{\beta,\epsilon}-1\}$.
Then, when $n_\alpha>0$	
\beqan
\lefteqn{
	\Pr\Big(\sup_{u\in [\alpha,\beta]} U(u)-U_n(u)  > \epsilon\Big)=	
	\sum_{\ell=0}^{\tilde m}\binom{n}{\ell}\Big(\!\min\Big\{1\!-\!\frac{\ell}{n}\!-\!\epsilon,1\!-\!\alpha\Big\}\!\Big)^{n\!-\!\ell}\alpha^\ell} \\
&&+\!\!\sum_{\ell=\tilde m+1}^{\overline{n}_{\beta,\epsilon}-1}\!\! 
\binom{n}{\ell}\!\bigg(1\!-\! \frac{\ell}{n}\!-\!\epsilon\bigg)^{n\!-\!\ell}\bigg[\epsilon
\bigg(\!\frac{\ell}{n}\!+\!\epsilon\!\bigg)^{\ell\!-\!1\!}\!\!+\! \sum_{j=0}^{\tilde m-1}\!\bigg[\frac{n_\alpha\!-\!j\!}{n}\bigg]\binom{\ell}{j}\bigg(\frac{\ell\!-\!n_\alpha}{n}\bigg)^{\ell\!-\!j\!-\!1\!}\!\alpha^{j}\bigg]\,.
\eeqan
On the other hand, when $n_\alpha<0$, it holds	
\beqan
	\Pr\Big(\sup_{u\in [\alpha,\beta]} U(u)-U_n(u)  > \epsilon\Big)=	\!\!\sum_{\ell=0}^{\overline{n}_{\beta,\epsilon}-1}\!\! 
\binom{n}{\ell}\!\bigg(1\!-\! \frac{\ell}{n}\!-\!\epsilon\bigg)^{n\!-\!\ell}\!\bigg[\epsilon
\bigg(\!\frac{\ell}{n}\!+\!\epsilon\!\bigg)^{\ell\!-\!1\!}\bigg]\,.
\eeqan
\end{theorem}

\begin{remark}[Left-Right tails]
	Interestingly, note that it holds
	 $\rdelta_{[\uu,\ou]}(n,\epsilon)
	=	\ldelta_{[1\!-\!\ou,1\!-\!\uu]}(n,\epsilon)$. This should not be surprising since $\sup_{u\in [\alpha,\beta]} U_n(u)-U(u)$ and $\sup_{u\in [1-\beta,1-\alpha]} U(u)-U_n(u)$ indeed have the same law.
	Also we have the trivial bound $\sup_{u \in [\alpha,\beta]} U_n(u) \!-\!U(u) \in [-\beta,1\!-\!\alpha]$, as well as 
	 $\sup_{u \in [\alpha,\beta]} U(u) \!-\!U_n(u) \in [\alpha\!-\!1,\beta]$.
	 Hence $\ldelta_{[\alpha,\beta]}(n,\epsilon)=0$ for $\epsilon>1-\alpha$
	 and $\rdelta_{[\alpha,\beta]}(n,\epsilon)=0$ for $\epsilon>\beta$.
	\end{remark}

\begin{corollary}[Exact concentration in specific cases]\label{cor:specific}If $\epsilon\leq1-\alpha$, then
	\beqan\Pr\bigg( \sup_{u \in [\alpha,1]} U_n(u) -U(u)  >
	\epsilon
	\bigg) & =&\sum_{\ell=0}^{\lceil n(1\!-\!\alpha\!-\!\epsilon)\rceil\!-\!1}\epsilon\Big(\frac{\ell}{n}\!+\!\epsilon\Big)^{\ell\!-\!1} \binom{n}{\ell} \Big(1\!-\!\epsilon\!-\!\frac{\ell}{n}\Big)^{n\!-\!\ell}\,.
	\eeqan
Now, if $\epsilon\leq1$, then for $m=\min\{\lfloor n_\beta\rfloor+1,\lceil n(\beta-\epsilon) \rceil-1\}$,

	\vspace{-3mm}\noindent
	\beqan
	\lefteqn{
		\Pr\Big(\sup_{u\in [0,\beta]} U_n(u)-U(u)  > \epsilon\Big)=
		\sum_{\ell=0}^{m}\binom{n}{\ell}\Big(\!\min\Big\{1\!-\!\frac{\ell}{n}\!-\!\epsilon,\beta\Big\}\!\Big)^{n\!-\!\ell}(1-\beta)^\ell} \\
	&&+\!\!\sum_{\ell=m+1}^{\lceil n(1-\epsilon) \rceil-1}\!\! 
	\binom{n}{\ell}\!\bigg(1\!-\! \frac{\ell}{n}\!-\!\epsilon\bigg)^{n\!-\!\ell}\bigg[\epsilon
	\bigg(\!\frac{\ell}{n}\!+\!\epsilon\!\bigg)^{\ell\!-\!1\!}\!\!+\! \sum_{j=0}^{m-1}\!\bigg[\frac{n_\beta\!-\!j}{n}\bigg]\binom{\ell}{j}\bigg(\frac{\ell\!-\!n_\beta}{n}\bigg)^{\ell\!-\!j\!-\!1\!}\!(1\!-\!\beta)^{j}\bigg]\,.
	\eeqan
	
	If $\beta\geq\epsilon$, then
	
	\vspace{-6mm}\noindent	
	\beqan
		\Pr\Big(\sup_{u\in [0,\beta]} U(u)-U_n(u)  > \epsilon\Big)=\!\sum_{\ell=0}^{\lceil n(\beta-\epsilon) \rceil-1}\!\! 
	\binom{n}{\ell}\!\bigg(1\!-\! \frac{\ell}{n}\!-\!\epsilon\bigg)^{n\!-\!\ell}\!\bigg[\epsilon
	\bigg(\!\frac{\ell}{n}\!+\!\epsilon\!\bigg)^{\ell\!-\!1\!}\bigg]\!\,.
	\eeqan
	
	Now, if $\epsilon<1$,  for $\tilde m=\min\{\lfloor n_\alpha\rfloor+1,\lceil n(1-\epsilon) \rceil-1\}$,
	
	\vspace{-3mm}\noindent
		\beqan
	\lefteqn{
		\Pr\Big(\sup_{u\in [\alpha,1]} U(u)-U_n(u)  > \epsilon\Big)=	
		\sum_{\ell=0}^{\tilde m}\binom{n}{\ell}\Big(\!\min\Big\{1\!-\!\frac{\ell}{n}\!-\!\epsilon,1\!-\!\alpha\Big\}\!\Big)^{n\!-\!\ell}\alpha^\ell} \\
	&&+\!\!\sum_{\ell=\tilde m+1}^{\lceil n(1-\epsilon) \rceil-1}\!\! 
	\binom{n}{\ell}\!\bigg[\bigg(1\!-\! \frac{\ell}{n}\!-\!\epsilon\bigg)^{\!n\!-\!\ell}\bigg]\bigg[\epsilon
	\bigg(\!\frac{\ell}{n}\!+\!\epsilon\!\bigg)^{\ell\!-\!1\!}\!\!+\! \sum_{j=0}^{\tilde m-1}\!\bigg[\frac{n_\alpha\!-\!j}{n}\bigg]\binom{\ell}{j}\bigg(\frac{\ell-\!n_\alpha}{n}\bigg)^{\ell\!-\!j\!-\!1\!}\!\alpha^{j}\bigg]\,.
	\eeqan
\end{corollary}

\begin{myproof}{of Corollary~\ref{cor:specific}}
	In case $\beta=1$,then
	$n_\beta = -n\epsilon<0$ and 	
	$\overline{n}_{\alpha,\epsilon}=\lceil n(1-\alpha-\epsilon)\rceil$. This yields the first equality.
	In case $\alpha=0$ and $\epsilon<1$, then
	$n_\beta = n(1-\beta-\epsilon)$ and 	
	$\overline{n}_{\alpha,\epsilon}=\lceil n(1-\epsilon)\rceil$, which yields the second result.
We proceed similarly for the right-tail inequality.
\end{myproof}

The previous result shows that for the global concentration ($[\uu,\ou]=[0,1]$), we recover the classical DKW derivation.
Indeed, from Corollary~\ref{cor:specific}, we get the following expressions
\beqan
\Pr\bigg( \sup_{u \in [0,1]} U_n(u) -U(u) > 
\epsilon
\bigg) & =&\sum_{\ell=0}^{\lceil n(1\!-\!\epsilon)\rceil\!-\!1}\binom{n}{\ell}\epsilon\Big(\frac{\ell}{n}\!+\!\epsilon\Big)^{\ell\!-\!1}  \Big(1\!-\!\epsilon\!-\!\frac{\ell}{n}\Big)^{n\!-\!\ell}\\
\Pr\Big(\sup_{u\in [0,1]} U(u)-U_n(u) > \epsilon\Big)&=&\!\sum_{\ell=0}^{\lceil n(1-\epsilon) \rceil-1}\!\! 
\binom{n}{\ell}\epsilon
\bigg(\!\frac{\ell}{n}\!+\!\epsilon\!\bigg)^{\ell\!-\!1\!}\!\bigg(1\!-\! \frac{\ell}{n}\!-\!\epsilon\bigg)^{n\!-\!\ell}\,,
\eeqan
while, on the other hand, from \cite[eq.(50) p.9]{smirnov1944approximate}, we get the following equivalent expression
\beqan
\Pr\bigg( \sup_{u \in [0,1]} U_n(u) -U(u) > 
\epsilon
\bigg) & =& (1-\epsilon)^n + \sum_{\ell=\lfloor n\epsilon\rfloor+1}^{n-1} \binom{n}{\ell} \epsilon \bigg(\frac{\ell}{n}-\epsilon\bigg)^\ell \bigg(1- \frac{\ell}{n}+\epsilon\bigg)^{n-\ell-1}\\
&=& \sum_{\ell=\lfloor n\epsilon\rfloor+1}^{n} \binom{n}{\ell} \epsilon \bigg(\frac{\ell}{n}-\epsilon\bigg)^\ell \bigg(1- \frac{\ell}{n}+\epsilon\bigg)^{n-\ell-1}\\
&=&\sum_{\ell=0}^{n-\lfloor n\epsilon\rfloor-1} \binom{n}{\ell} \epsilon \bigg(1-\frac{\ell}{n}-\epsilon\bigg)^{n-\ell} \bigg(\frac{\ell}{n}+\epsilon\bigg)^{\ell-1}\,.
\eeqan
Now, Theorem~\ref{thm:concLeft},\ref{thm:concRight} and Corollary~\ref{cor:specific} provide a detailed control of the CDF concentration over arbitrary intervals $[\alpha,\beta]$ of $[0,1]$. As we detail in Section~\ref{sec:risk}, this is of special interest when $\alpha$ and $\beta$ are risk-levels and one is interested in functionals of the CDF such as the conditional value at risk or more generic risk measures, since in that case $\alpha$ and $\beta$ are known and specified by the practitioner.

\section{Numerical illustration of the bounds}\label{sec:xps}
In this section, we provide a numerical illustration of the concentration bounds provided in Theorem~\ref{thm:concLeft}
and  Theorem~\ref{thm:concRight}.  
This is made possible thanks to the fact the functions $\delta_{[\uu,\ou]}$ and $\delta_{[\uu,\ou]}$ are fully explicit although with fairly complicated expressions. We illustrate these functions as well as their inverse in this section in order to provide intuition and also to show that they can be easily computed. We thus humbly suggest the practitioner to make use of the quantity $\lepsilon_{[\uu,\ou]}(n,\delta)$ instead of the approximation  $\sqrt{\ln(1/\delta)/2n}$ suggested by DKW and Massart. Indeed, this approximation is primarily interesting for large values of $n$ in order to get an idea of the scaling of the bound. However, for small values of $n$, this approximation can be detrimental.

\begin{figure}[H]
	\centering
	\includeg{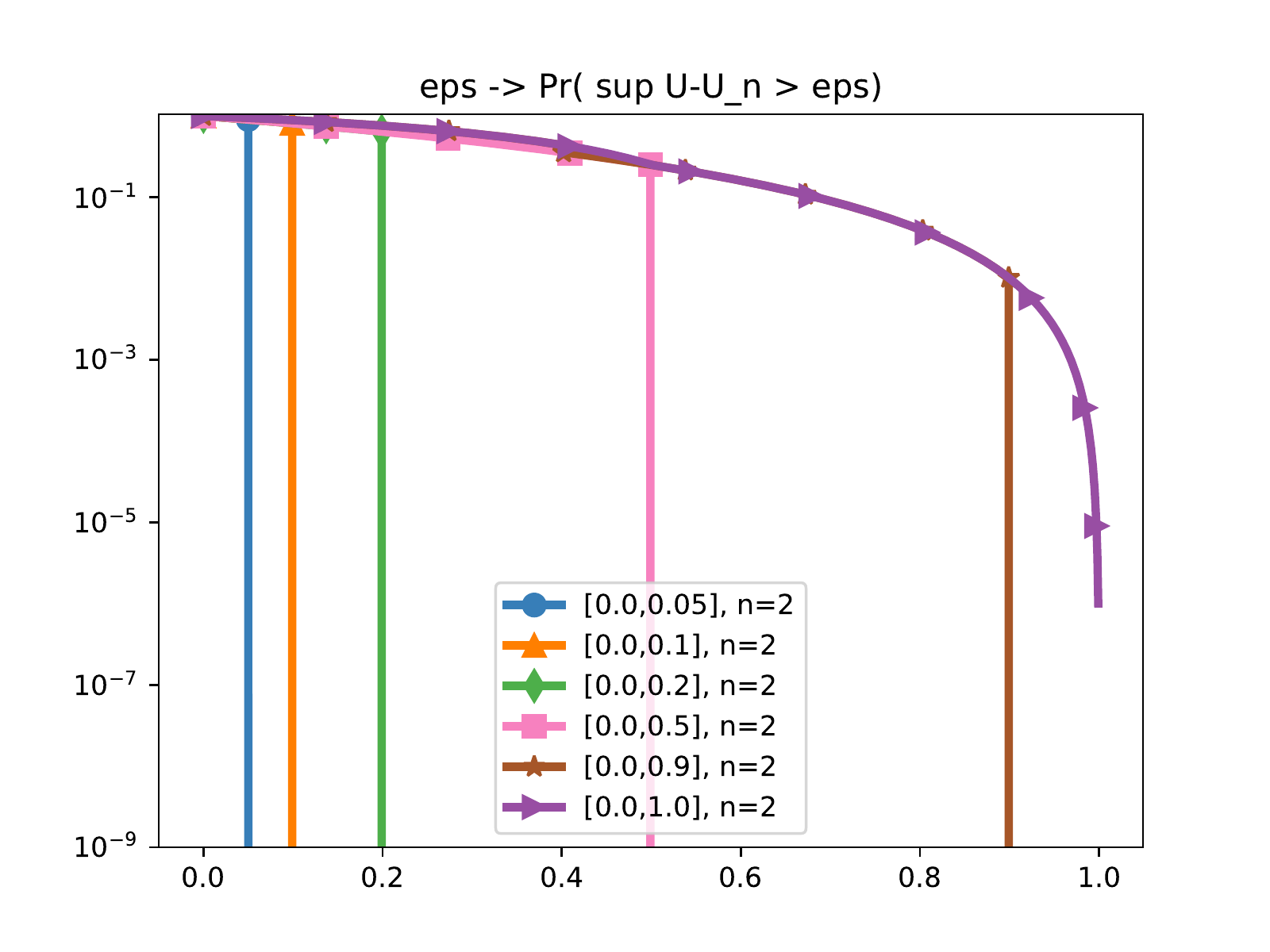}
	\includeg{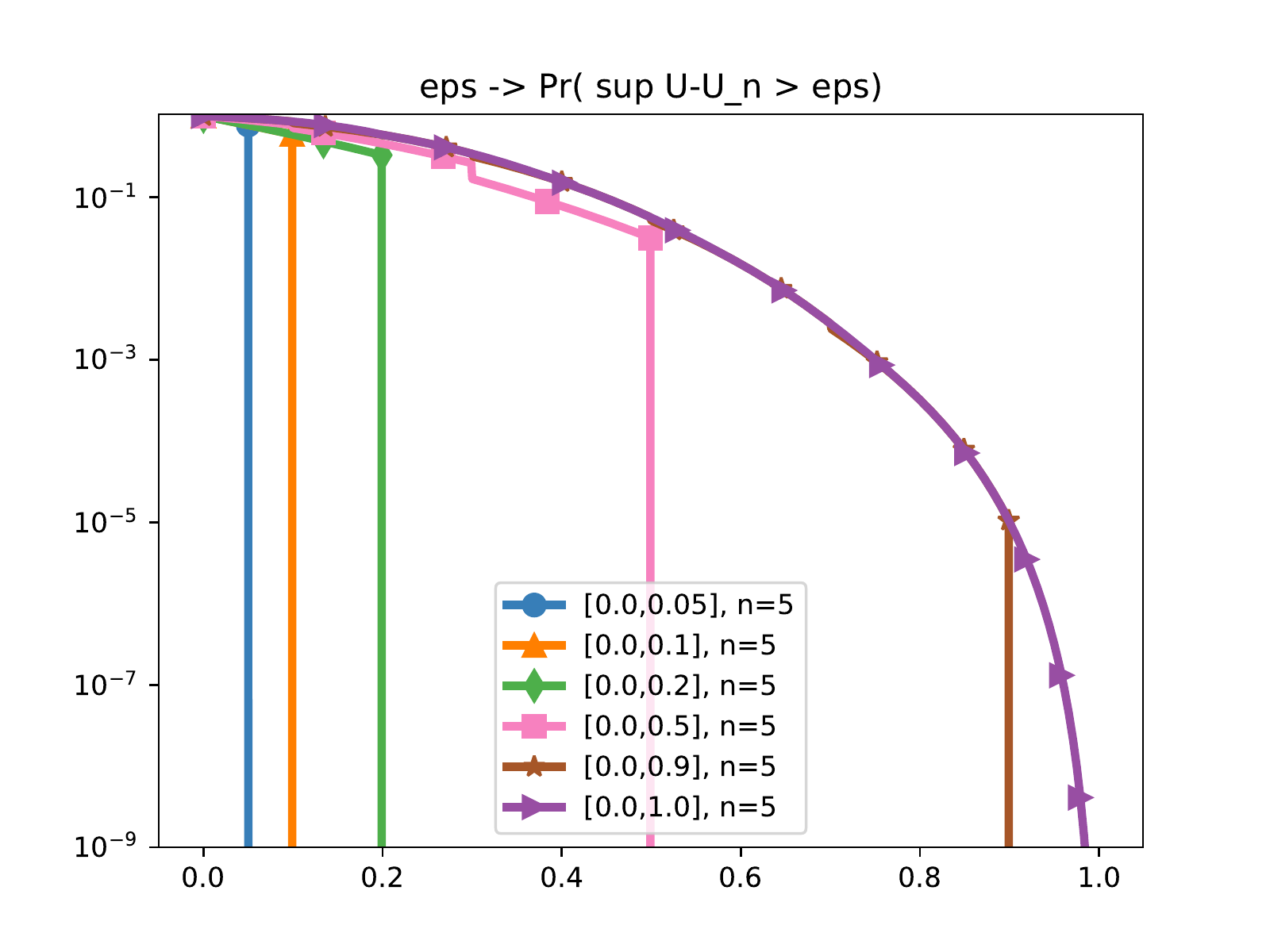}\\
	\includeg{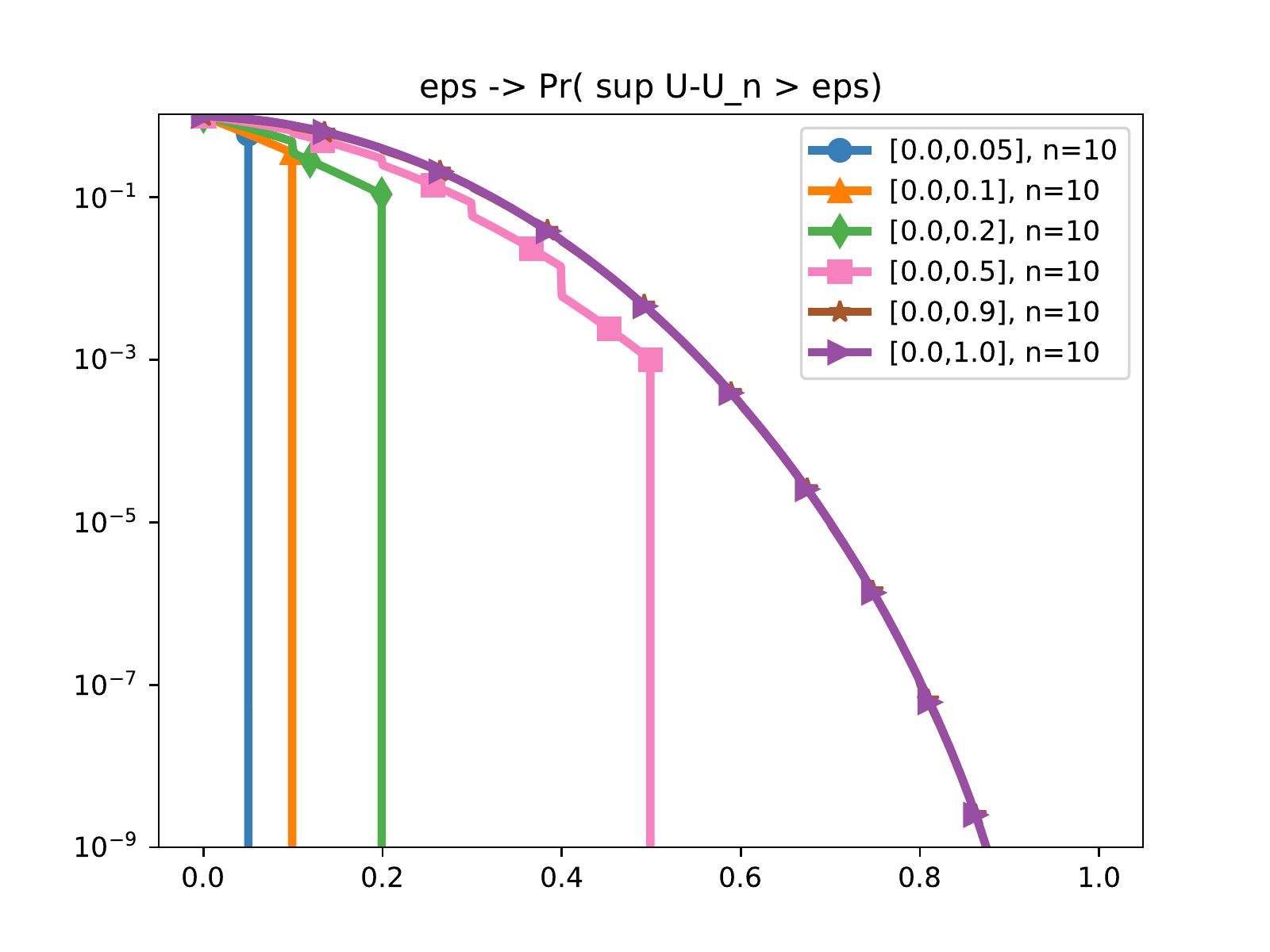}
	\includeg{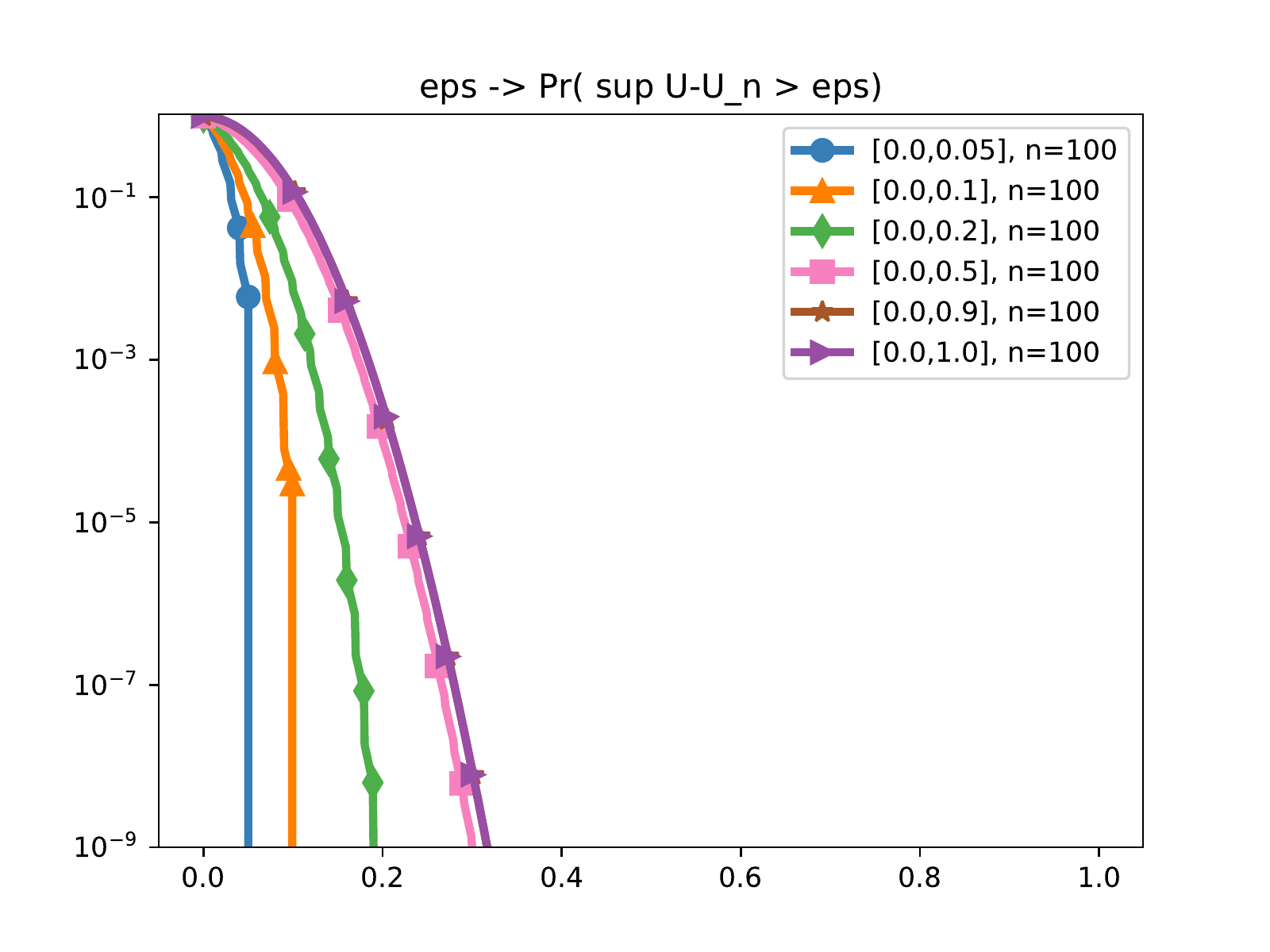}\\
	\includeg{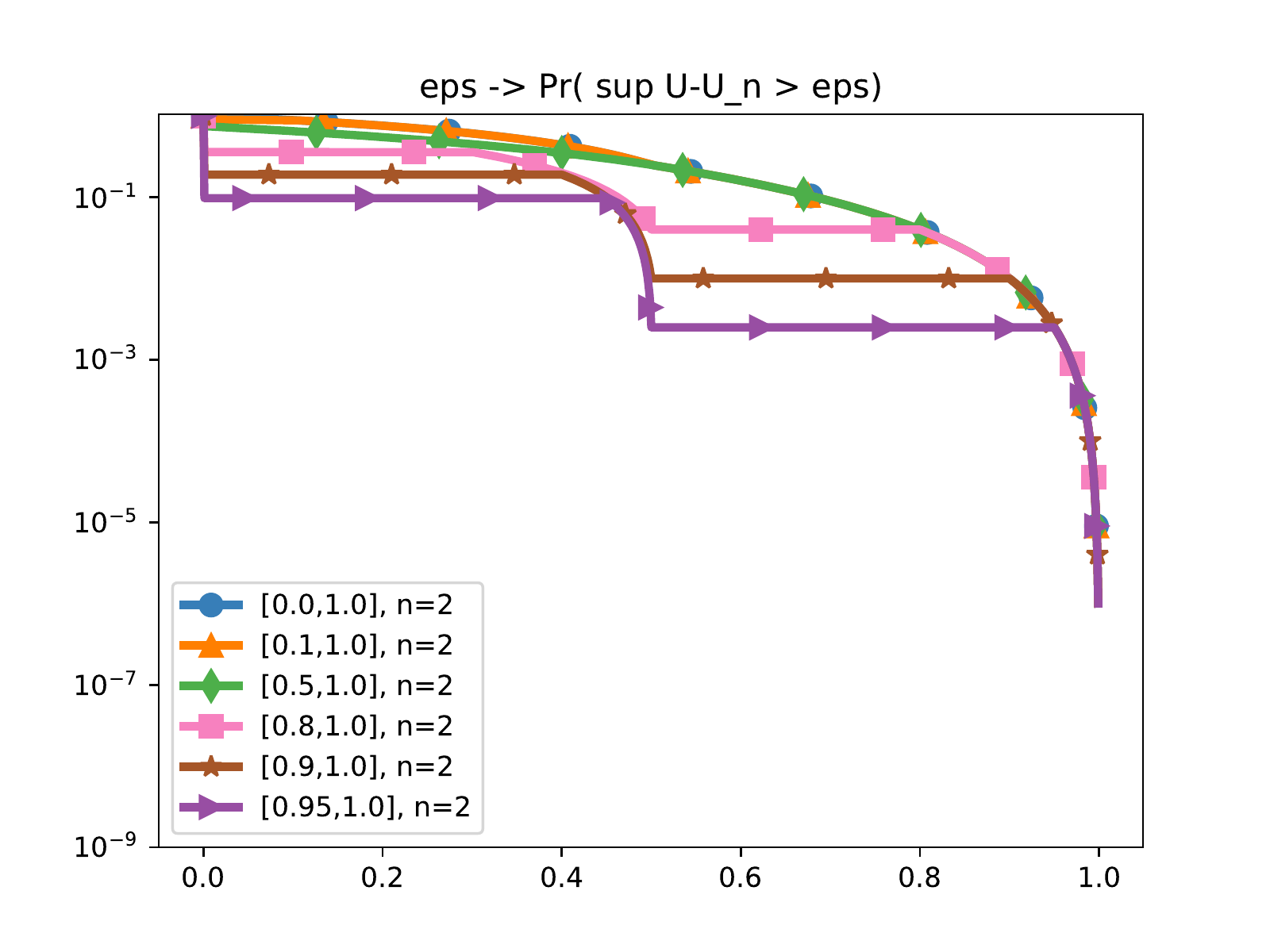}
	\includeg{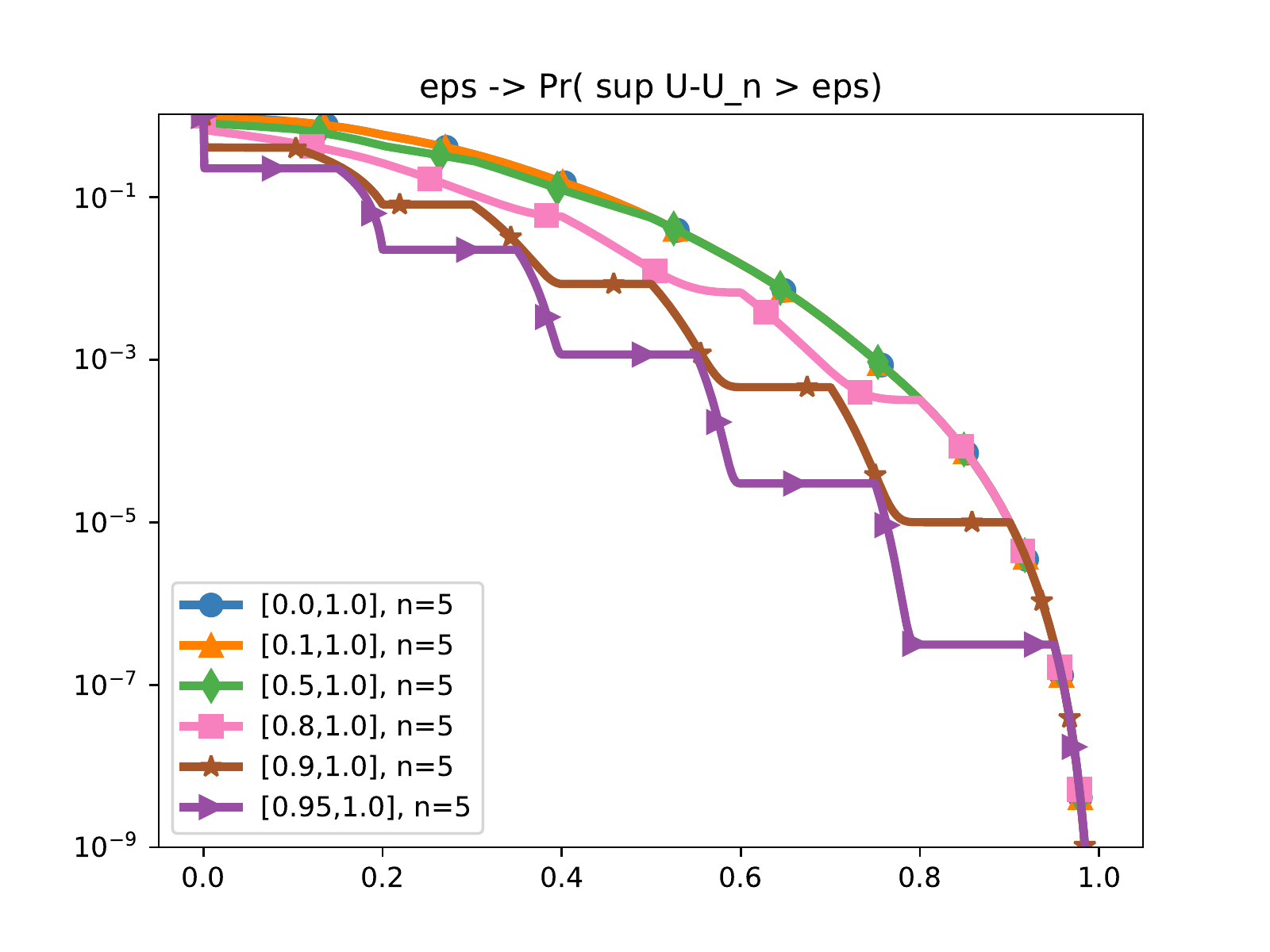}\\
	\includeg{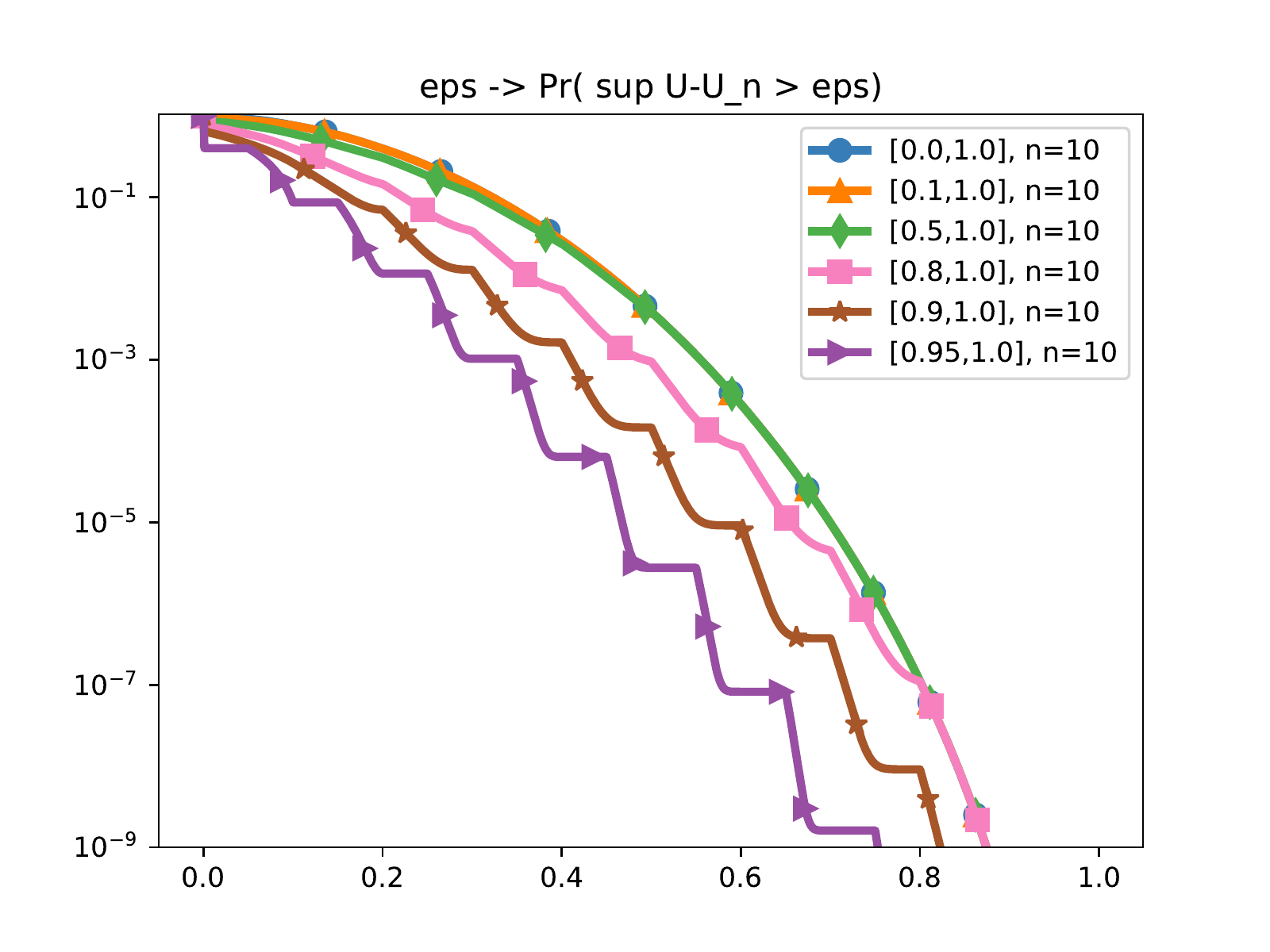}
	\includeg{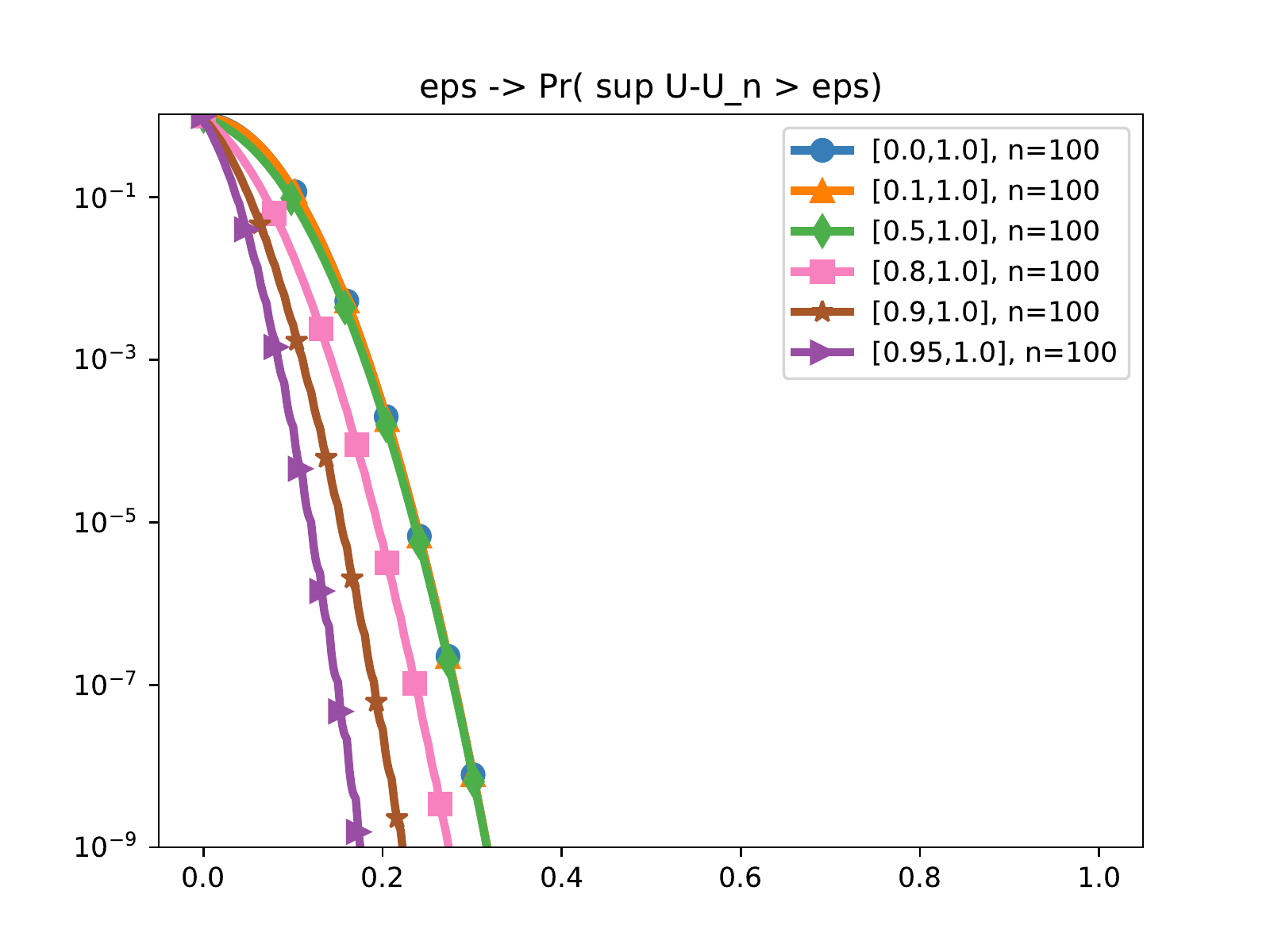}
	\caption{Plot of $\epsilon\to\tilde \delta_{[\uu,\ou]}(n,\epsilon)$ for various values of $n$ and 
		interval $[\uu,\ou]$.}
	\label{fig:delta1}
\end{figure}

In Figure~\ref{fig:delta0}, we plot $\epsilon\to\delta_{[\uu,\ou]}(n,\epsilon)$ for various intervals $[\uu,\ou]$, and in Figure~\ref{fig:delta1}, we plot
 $\epsilon\to\tilde\delta_{[\uu,\ou]}(n,\epsilon)$. The plots highlight  the non-trivial behavior of these functions, especially for small values of $n$, having plateaus, abrupt changes and non-linear behavior.  These functions become smoother as $n$ increases (which is expected). Let us note the striking impact of changing $[\uu,\ou]$ on the resulting function.
  In Appendix~\ref{app:MCMC}, we further provide in Figure~\ref{fig:MCMC1} and Figure~\ref{fig:MCMC2} a numerical comparison between the computation of the exact probabilities from Theorem~\ref{thm:concLeft}, and direct Monte Carlo simulations of the bound. The plots were obtained using  $M=10^4$ simulations (we consider that the accuracy of such plots is hence  good enough for values not less than $10^{-3}$), and perfectly match the theoretical bounds, as expected.

\begin{figure}[H]
	\centering
	\includeg{delta_00_005-01-02-05-09-10_n2}
	\includeg{delta_00_005-01-02-05-09-10_n5}\\
	\includeg{delta_00_005-01-02-05-09-10_n10}
	\includeg{delta_00_005-01-02-05-09-10_n100}\\
	\includeg{delta_00-01-05-08-09-095_10_n2}
	\includeg{delta_00-01-05-08-09-095_10_n5}\\
	\includeg{delta_00-01-05-08-09-095_10_n10}
	\includeg{delta_00-01-05-08-09-095_10_n100}
	\caption{Plot of $\epsilon\to\delta_{[\uu,\ou]}(n,\epsilon)$ for various values of $n$ and 
		interval $[\uu,\ou]$.}
	\label{fig:delta0}
\end{figure}

\begin{figure}
	\centering
	\includeg{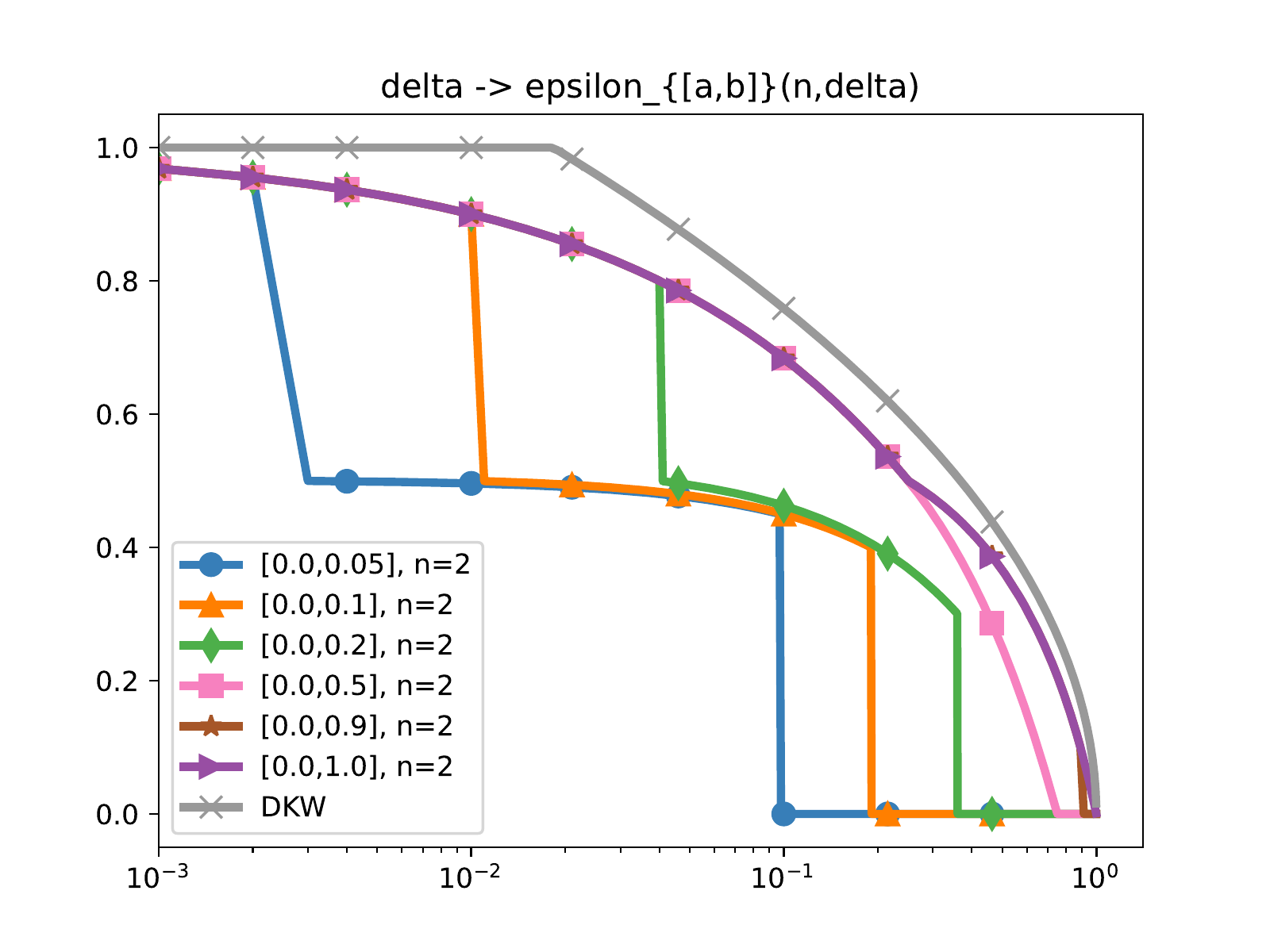}
	\includeg{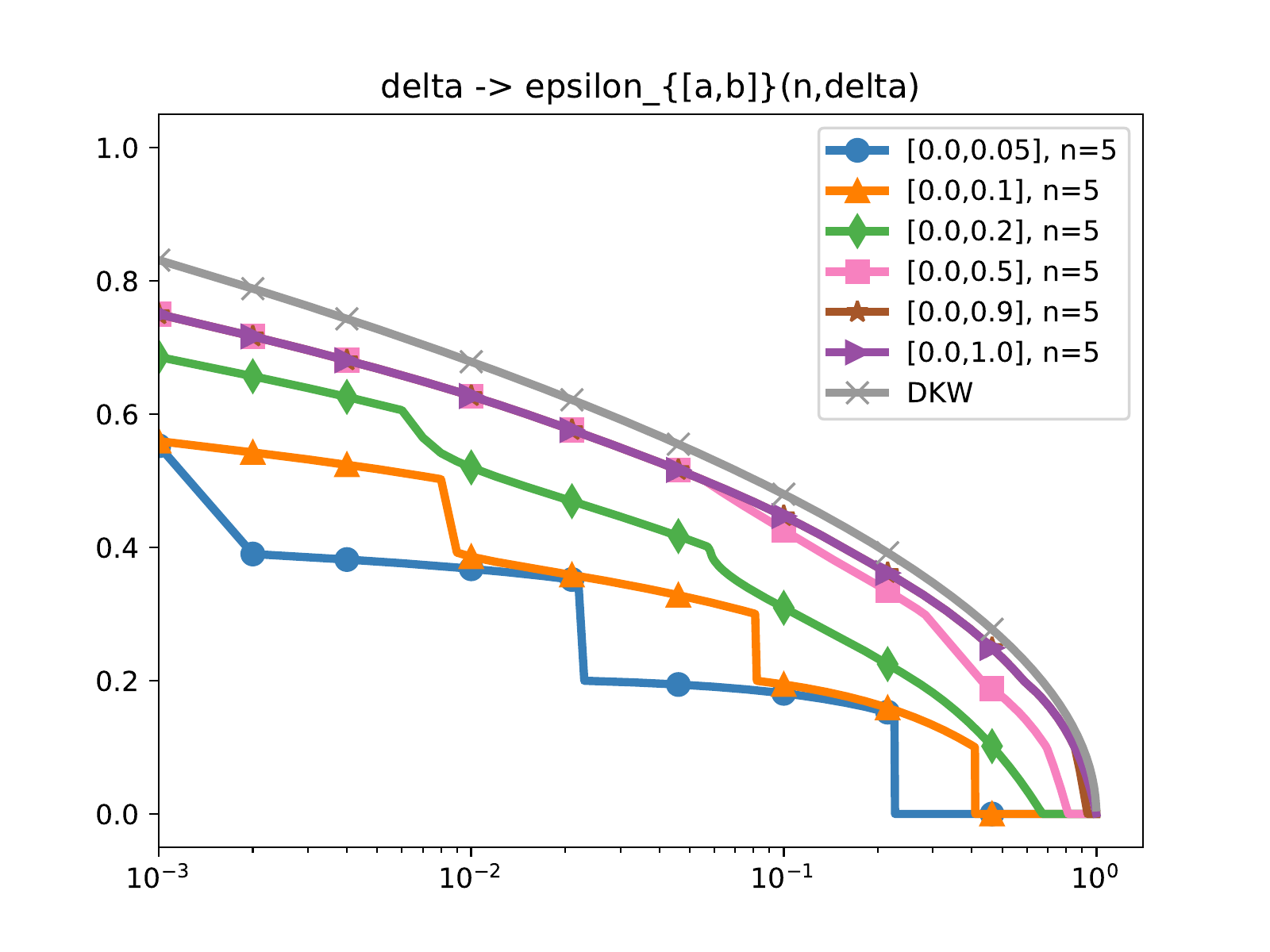}\\
	\includeg{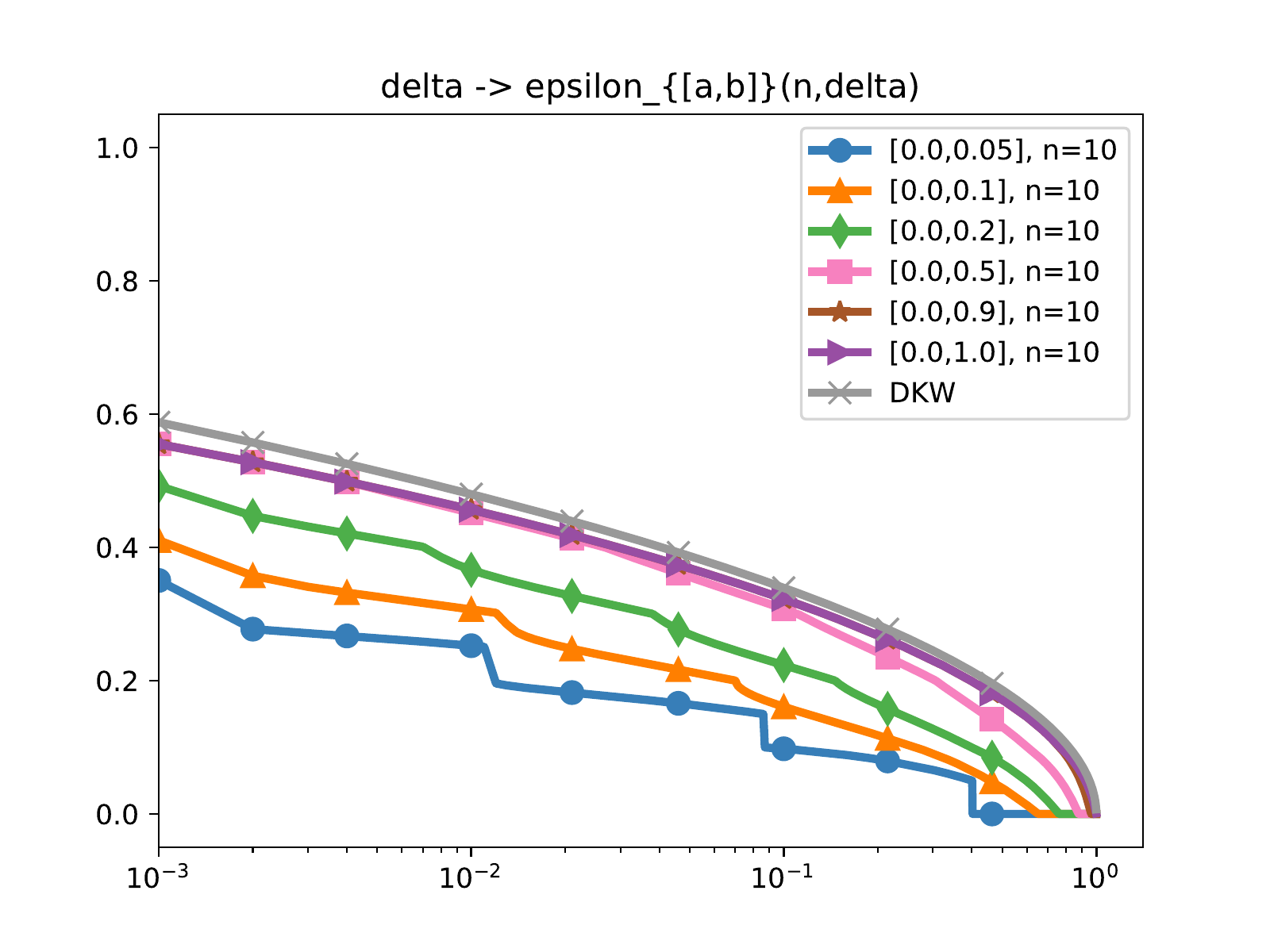}
	\includeg{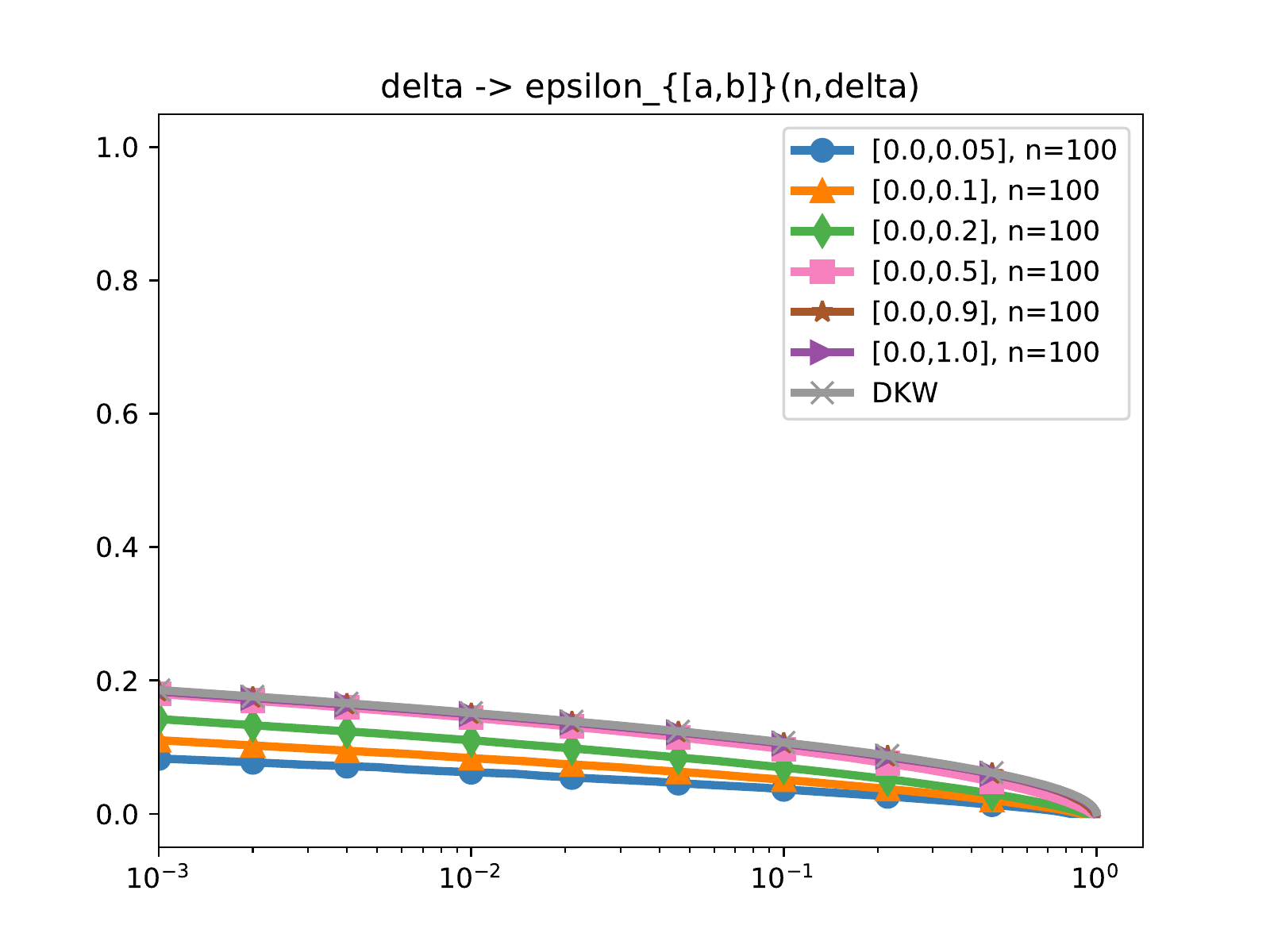}\\
	\includeg{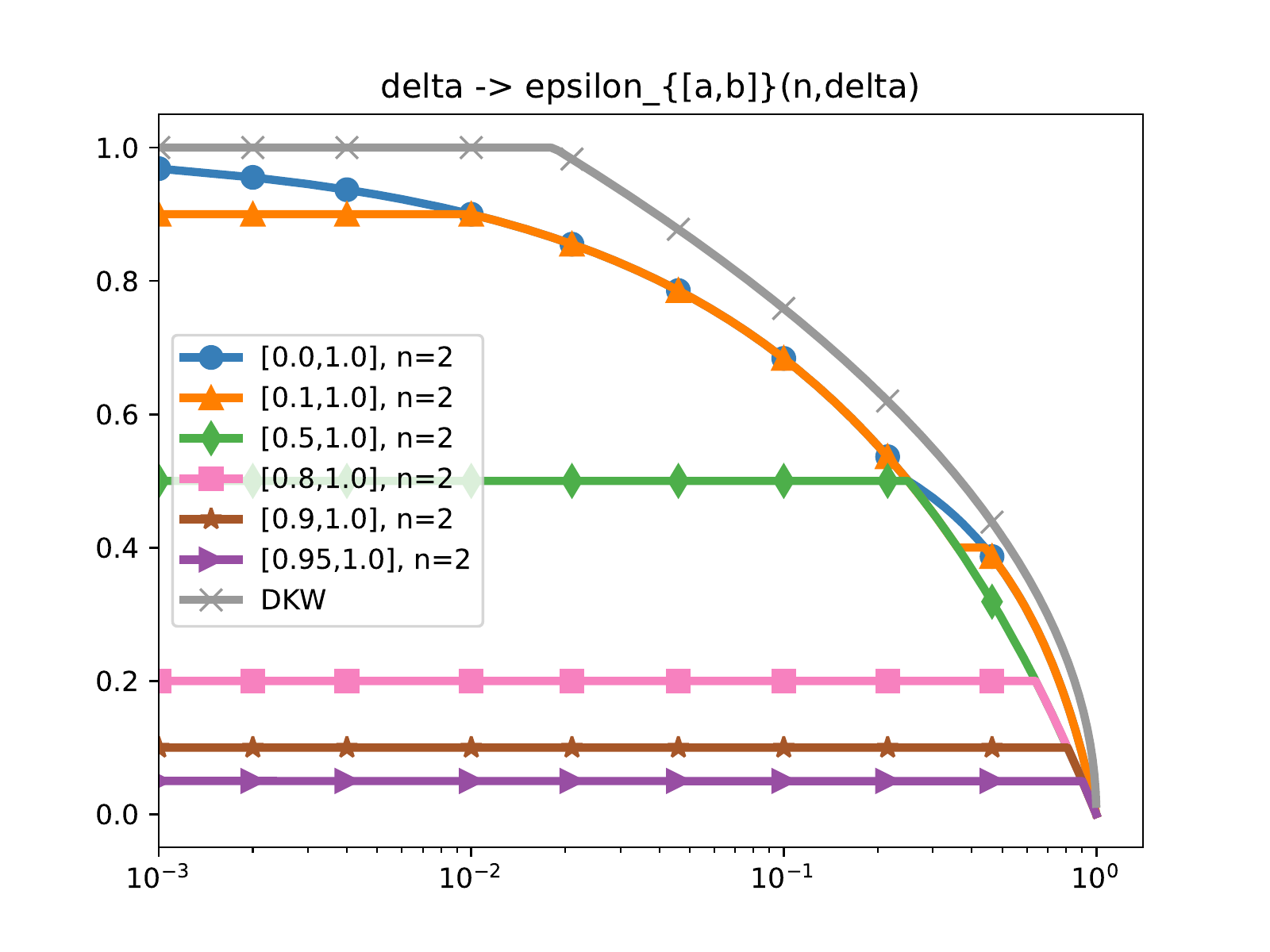}
	\includeg{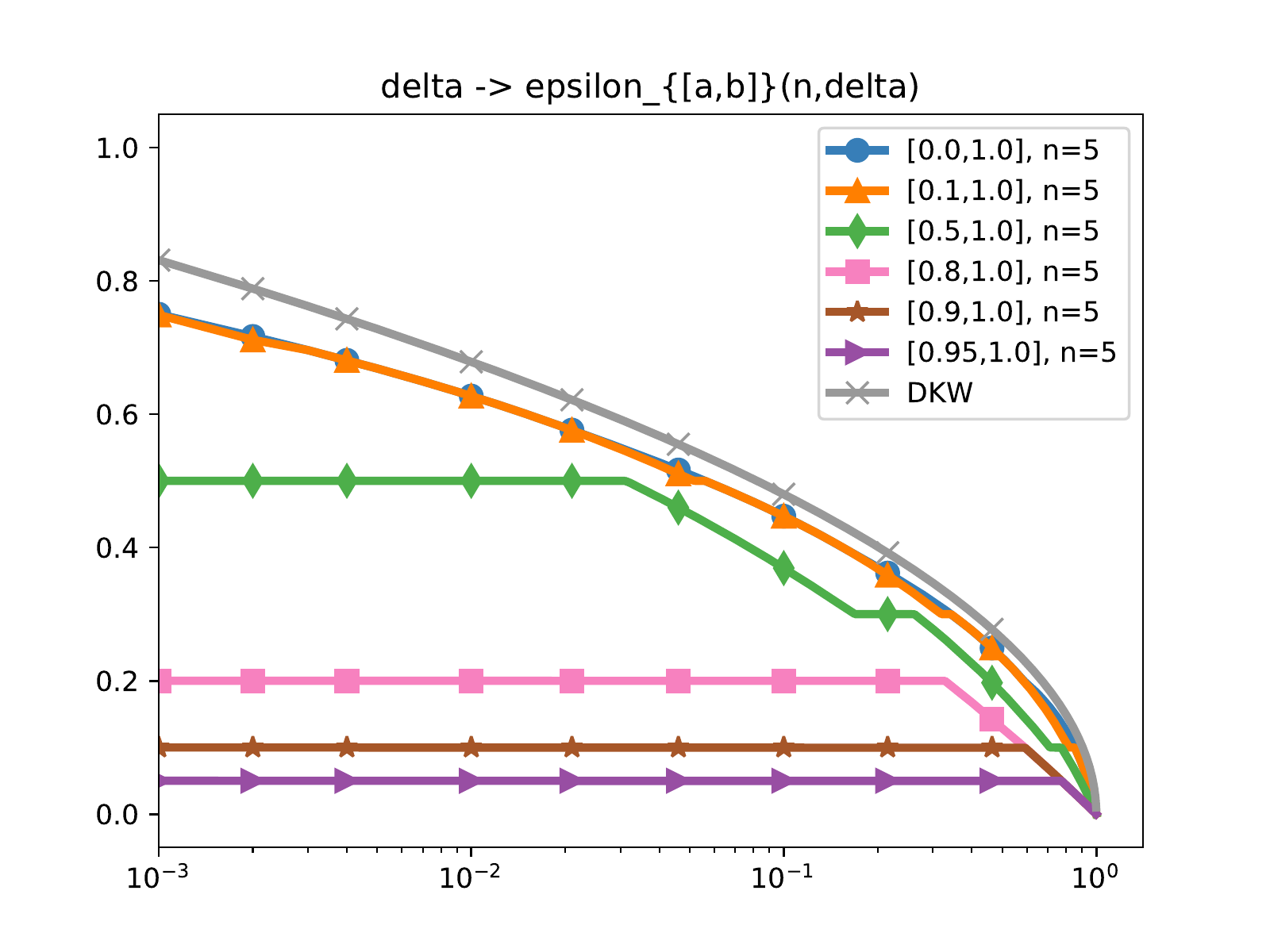}\\
	\includeg{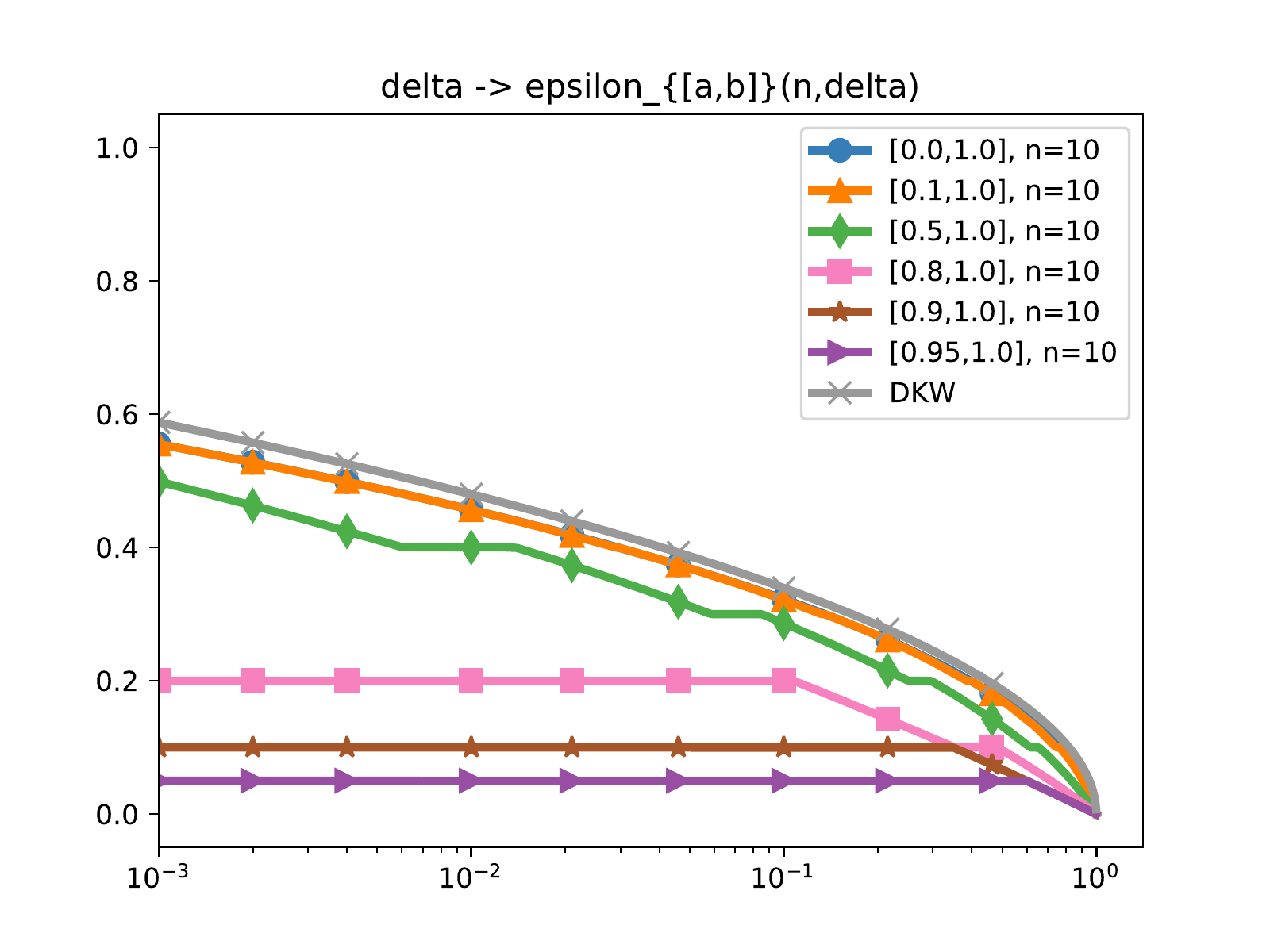}
	\includeg{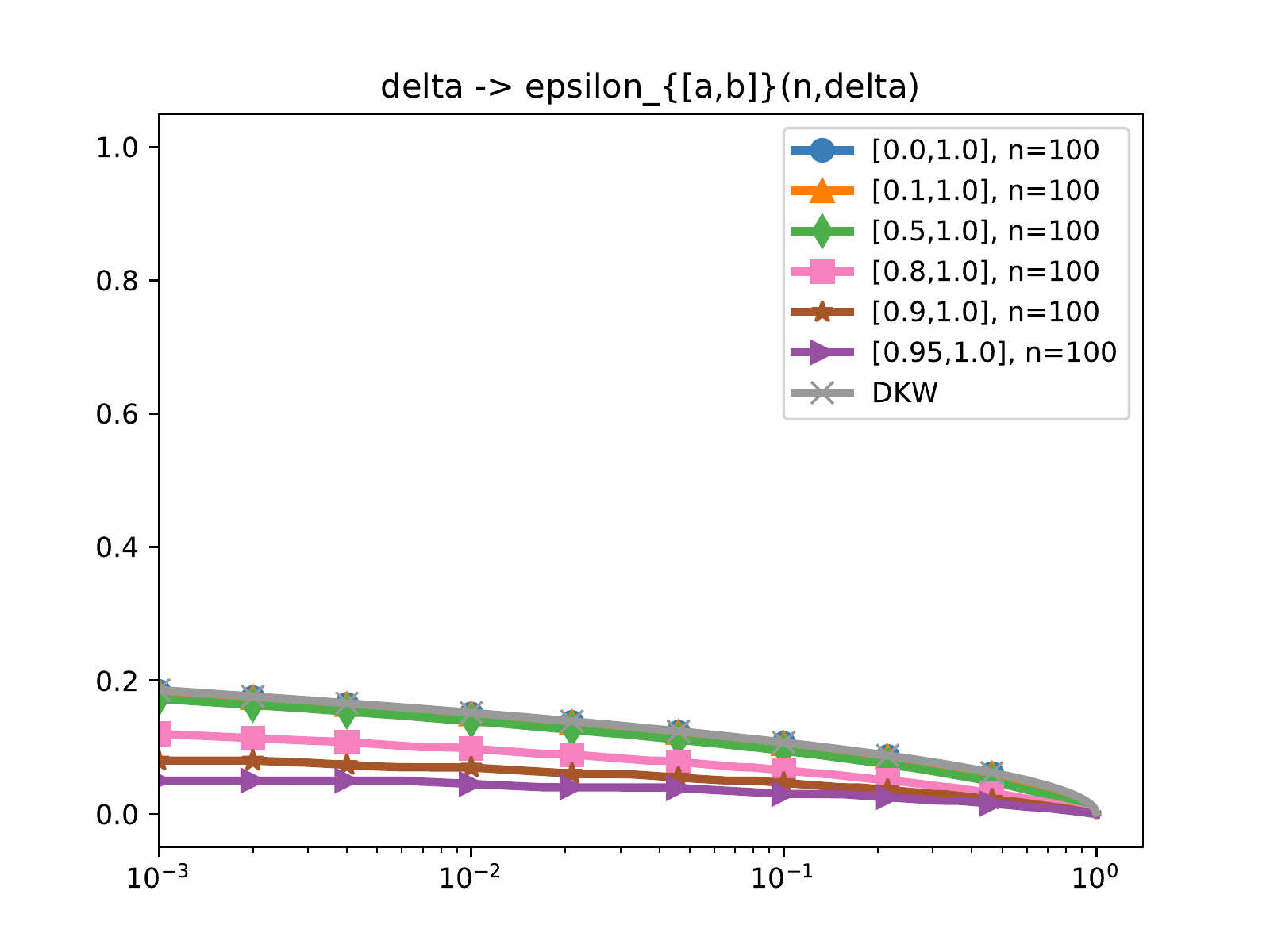}\\
	\caption{Plot of $\delta\to\epsilon_{[\uu,\ou]}(n,\delta)$ for various values of $n$ and 
		intervals $[\uu,\ou]$.}
	\label{fig:epsilon0}
\end{figure}

In  Figure~\ref{fig:epsilon0},  we plot the inverse function $\delta\to\lepsilon_{[\uu,\ou]}(n,\delta)$ as well as the 
explicit approximation given by Massart of the inverse probability function $\lepsilon_{[0,1]}$ for comparison.
More precisely, we plot $\epsilon(n,\delta)=\min\big\{\sqrt{\frac{\ln(1/\delta)}{2n}},1\}$, simply called ``DKW'' in the plots
(let us recall this function is only shown to control the deviations provided that $\delta \in[0,0.5)$, but this does not prevent us from plotting it for $[0,1]$).
From the perspective of confidence bands, we believe that this function is the most interesting to plot, as it shows the exact magnitude of the estimation errors, and as such are \textit{not} improvable.
The plots confirm that the Massart bound is a valid upper bound on the exact probability functions, as can be seen by comparing in each plot the curve called DKW with the one obtained for the full set $[0,1]$.
On the other hand, we also observe  that this convenient but simplified bound can be quite conservative when considering a local supremum as opposed to a supremum over $[0,1]$. 
Also, the plots confirm that a simple expression cannot be obtained to accurately describe the behavior of $\lepsilon_{[\uu,\ou]}$.
The simplified bound provided by Massart is in this sense an especially convenient tool to trade-off accuracy and simplicity of the expression. However, we recommend using the exact bounds out of tightness, especially when the number of considered samples $n$ is \textit{small}. Indeed in that case, the numerical cost of computing the exact bound should not be high.

\medskip
Interestingly, we remark that a computation of $\lepsilon_{[\uu,\ou]}$ and $\repsilon_{[\uu,\ou]}$ can be achieved numerically (hence, yielding the different plots) and readily translates into confidence bands over the local supremum deviation, offering tighter bounds for the practitioner.
The computation of $\lepsilon_{[\uu,\ou]}$ can further be done simply, e.g. by  dichotomous search, in an efficient way up to a desired precision $\eta$. 
For instance Figure~\ref{fig:epsilon0} has been obtained using precision $\eta=10^{-7}$. We believe this is a key point, as it indicates one can use such bounds in various applications, at a controlled computational cost. Note also that one can tabulate this function off-line, which is convenient for applications involving sequential decisions with increasingly many observations.

\medskip
For ease of use, we provide for the interested reader the implementation details in the gitlab repository that is available at
\verb|https://gitlab.inria.fr/omaillar/article-companion/-/tree/|

\noindent
\verb|master/2020-local-dkw|. In order to avoid some numerical instabilities, the code uses a simple trick to replace expressions such as $x^\ell$ by $\exp(\ell\ln(x))$. It also enables to reproduce all the plots displayed in the different figures of this article. Further, we provided a method to tabulate the inverse functions 
$\lepsilon_{[\uu,\ou]}$ and $\repsilon_{[\uu,\ou]}$ for any given sub-interval of $[0,1]$, number of observations $n$ and given probability level $\delta$, which we believe to be useful for the practitioner.

\newpage

\section{Application to the CVAR and other functionals of the CDF}\label{sec:risk}
In this section, we first introduce some background material and intuition about the risk measure known as the conditional value-at-risk ($\CVAR$).
For a probability distribution $\nu\in\cP(\cX)$ taking values in $\cX\subset\Real$, let $\Esp_\nu$ denote its expectation operator and $F:x\mapsto\Pr_\nu(X\leq x)$ its Cumulative Distribution Function (CDF).
Given $\alpha$, chosen by the practitioner, one may want to measure a risk using $\Esp_\nu[X |X<x_\alpha]$ for some $x_\alpha\in\cX$ such that $F(x_\alpha)=\alpha$ in case $X$ is a variable we want to maximize (reward), or $\Esp_\nu[X |X>x_\kappa]$ for some $x_\kappa\in\cX$ such that $F(x_\kappa)=\kappa$ in case $X$ is a variable we want to minimize (loss). Unfortunately, when $F$ is not continuously increasing, the points $x_\alpha,x_\kappa$
may not exist. This is typically the case when $\nu=\nu_n$ is the empirical distribution built from $n$ observation. 
One way to overcome this problem is to introduce
the Left Tail Conditional Expectation
$\LTCE_{\alpha}(\nu)=
\Esp_\nu\bigg[X \bigg| X \leq \ox_{\alpha}(\nu)
\bigg]$ together with the upper Value at Risk
$	\ox_{\alpha}(\nu)= \inf\{ x \in\Real : F(x)>\alpha\}$, or the Right Tail Conditional Expectation
$\RTCE_{\kappa}(\nu) = \Esp_\nu\bigg[X \bigg| X > \ux_{\kappa}(\nu) \bigg]$
together with the lower Value at Risk $\ux_{\kappa}(\nu)
= \inf\{ x \in\Real : F(x)\geq\kappa\}$.
Unfortunately, since in general $F(\ox_{\alpha}(\nu))\geq 
F(	\ux_{\alpha}(\nu)) \geq \alpha$ without equality,
this makes such quantities difficult to interpret.
The idea behind the classical definition of the \CVAR\ is hence to interpolate between the values of $F$ around $\alpha$, using the following\footnote{the superscript
	$r$ stands for rewards, and $\ell$ for losses} optimization problems (\cite{rockafellar2000optimization})
\beqan
\rCVAR_{\alpha}(\nu) = \sup_{x\in \Real}\bigg\{ \frac{1}{\alpha}\Esp[\min(X\!-\!x,0)]+x\bigg\}\,,\qquad
\lCVAR_{\kappa}(\nu) = \inf_{x\in \Real}\bigg\{x + \frac{1}{1-\kappa}\Esp[\max(X\!-\!x,0)]\bigg\}\,.
\eeqan
Let $X$ be a non negative random variable upper bounded by $\overline{x}$ for which $F$ is increasing and continuous. Then the following is known  (see e.g. \cite{acerbi2002spectral}, or \cite{thomas2019concentration}; we reproduce this result in Appendix~\ref{app:cvar} for completeness)
\begin{proposition}[\CVAR\ to CDF reduction]\label{prop:cvar}
	The quantity
	$\ux_{\kappa}(\nu)$ is a solution to the $\CVAR^\ell_\kappa$ optimization problem, and the following rewriting holds
	\beqan
	\lCVAR_{\kappa}(\nu) &\!=\!& \frac{1}{1\!-\!\kappa}\Esp_\nu\Big[X\indic{X\!>\! \ux_{\kappa}(\nu)}\Big] = \RTCE_{\kappa}(\nu)
	= \ox\!-\! \frac{1}{1\!-\!\kappa}\int_{0}^{\ox}(F(x)\!-\!\kappa)_+dx\,.\\
	\rCVAR_{\alpha}(\nu) &\!=\!& \frac{1}{\alpha}\Esp_\nu\Big[X\indic{X\!\leq\! \ox_{\alpha}(\nu)}\Big] = \LTCE_{\alpha}(\nu)
	=\frac{1}{\alpha}\int_{0}^{\ox} (\alpha-F(x))_{+}dx\,.
	\eeqan
\end{proposition}
In particular the \CVAR\ writes as a function of the CDF in the form
$\rCVAR_{\alpha}(\nu) = h^r\Big(\int_{\cX}g^r\big(F(x)\big)dx\Big)$, where $h^r$ and $g^r$ are monotonic functions. 
Further $g^r(\beta)=(\alpha-\beta)_+$ has support $[0,\alpha)$ that is a strict subset of $[0,1]$ for $\alpha<1$.
This property is actually not limited to the CVaR and allows to focus on controlling the deviations of the CDF to later control the risk-measure.
Indeed, one can then build the confidence bands
\beqan
\CVAR_{\alpha,n,\delta}^{r,\pm} = h^r\Big(\int_{\cX}g^r\big(F_n(x) \pm \epsilon_{[0,\alpha]}(n,\delta)\big)dx\Big)\,.
\eeqan

\paragraph{Functional of the CDF}
We now present a generalization of this procedure to other functionals of the CDF.
In the sequel, we let $\cS_J= \bigg\{  \alpha = (\alpha_0,\dots,\alpha_J), \alpha_0=0<\alpha_1< \dots< \alpha_J=1\}$ denote the set of increasing sequences
partitioning the interval $[0,1]$ into $J$ segments.
Further, for $\beta\in\cS_J$, we let $\cS_K(\beta) = \bigg\{\alpha \in \cS_K : \exists k_0<k_1,\dots<k_J \text{ s.t.} \alpha_{k_j}=\beta_j\bigg\}$. We now introduce a  definition for convenience.

\begin{definition}[Locally right-Lipschitz function]\label{def:locrL}
A function $\phi$ is \emph{locally lower-right-Lipschitz} if it satisfies
\beqan
\exists J\in\Nat_\star, \beta\in\cS_J, \ell\in\Real_+^J\,\,
\forall j\in\{1,\dots,J\}\forall y\in[\beta_{j-1},\beta_{j}],\quad \forall \epsilon>0, \frac{\phi(y)-\phi(y+\epsilon)}{\epsilon}\leq\ell_j\,.
\eeqan
A function $\phi$ is \emph{locally upper-left-Lipschitz} if it satisfies
\beqan
\exists J\in\Nat_\star, \beta\in\cS_J, \ell\in\Real_+^J\,\,
\forall j\in\{1,\dots,J\}\forall y\in[\beta_{j-1},\beta_{j}],\quad \forall \epsilon>0, \frac{\phi(y-\epsilon)-\phi(y)}{\epsilon}\leq\ell_j\,.
\eeqan
\end{definition}
For illustration, let us remark that $\lCVAR_{\kappa}(\nu) =  \int_\cX \phi(F(x))dx$ with $\phi(y)=1- \left(\frac{y-\kappa}{1-\kappa}\right)_{+}$.
This non-increasing function is locally upper-left Lipschitz with 
$J=2$, $\beta=(0,\kappa,1)$ and $\ell=(0,\frac{1}{1-\kappa})$.
On the other hand, $\rCVAR_{\alpha}(\nu) =  \int_\cX \phi(F(x))dx$ with $\phi(y)=\left(\frac{\alpha -y}{\alpha}\right)_{+}$.
This non-increasing  function is locally lower-right Lipschitz with 
$J=2$, $\beta=(0,\alpha,1)$ and $\ell=(\frac{1}{\alpha},0)$. This motivates the following result.

\begin{theorem}[Functional of continuous CDF deviation]\label{thm:functional} For a distribution on $\cX$ with continuous CDF $F$, let $\mu = \int_\cX \phi(F(x))dx$, where $\phi$ is assumed to be a known, non-increasing function. For any $K\in\Nat,\alpha\in\cS_k$ and all $\delta_k\in[0,1], k\in\{1,\dots,K\}$ it holds that
\beqan
\Pr\bigg( \mu < \sum_{k=0}^{K-1}\int_{\underline{x}_{\alpha_k}}^{\underline{x}_{\alpha_{k+1}}}  \phi\Big(\hat F_n(x) + \tilde \epsilon_{[\alpha_k,\alpha_{k+1}]}(n,\delta_k)\Big)\bigg)&\leq& \sum_{k=1}^K\delta_k\,,\\
\Pr\bigg( \mu > \sum_{k=0}^{K-1}\int_{\underline{x}_{\alpha_k}}^{\underline{x}_{\alpha_{k+1}}}  \phi\Big(\hat F_n(x) -  \epsilon_{[\alpha_k,\alpha_{k+1}]}(n,\delta_k)\Big)\bigg)&\leq& \sum_{k=1}^K\delta_k\,,
\eeqan 
where $\ux_{\kappa}= \inf\{ x \in\Real : F(x)\geq\kappa\}$.
Further, when $\phi$ is locally lower-right Lipschitz with known $(J,\beta,\ell)$ 
and the quantile function $\kappa\mapsto\ux_{\kappa}$ is $\gamma$-Lipschitz, then
if $\alpha \in \cS_K(\beta)$, the following holds
\beqan
\Pr\bigg( \mu < \int_\cX \phi(\hat F_n(x))dx - \sum_{j=1}^J\ell_j \sum_{k=k_{j-1}}^{k_{j}-1} \gamma(\alpha_{k+1}-\alpha_{k})\tilde \epsilon_{[\alpha_k,\alpha_{k+1}]}(n,\delta_k)\bigg)\leq \sum_{k=1}^K\delta_k\,.
\eeqan 
If instead $\phi$ is locally upper-left Lipschitz with known $(J,\beta,\ell)$ and $\alpha \in \cS_K(\beta)$, then
\beqan
\Pr\bigg( \mu > \int_\cX \phi(\hat F_n(x))dx + \sum_{j=1}^J\ell_j \sum_{k=k_{j-1}}^{k_{j}-1} \gamma(\alpha_{k+1}-\alpha_{k}) \epsilon_{[\alpha_k,\alpha_{k+1}]}(n,\delta_k)\bigg)\leq \sum_{k=1}^K\delta_k\,.
\eeqan 
\end{theorem}
We only stated the result for a non-increasing function $\phi$.
Alternative results for  a non-decreasing function 
and corresponding upper-right or lower-left assumptions can be derived too.

\begin{myproof}{of Theorem~\ref{thm:functional}}
Indeed, first, using the definition of $\mu$ an the monotony of $\phi$, it holds for all $\alpha \in \cS_K(\beta)$
\beqan
\mu & =& \sum_{k=0}^{K-1} \int_{\underline{x}_{\alpha_k}}^{\underline{x}_{\alpha_{k+1}}}  \phi\Big(\hat F_n(x) + F(x)- \hat F_n(x)\Big)dx\\
	&\geq&\sum_{k=0}^{K-1} \int_{\underline{x}_{\alpha_k}}^{\underline{x}_{\alpha_{k+1}}}  \phi\Big(\hat F_n(x) + \sup_{\underline{x}_{\alpha_k}\leq x \leq \underline{x}_{\alpha_{k+1}}}F(x)- \hat F_n(x)\Big)dx
\eeqan
Now, using that $\sup_{\underline{x}_{\alpha_k}\leq x \leq \underline{x}_{\alpha_{k+1}}}F(x)- \hat F_n(x) \stackrel{\cL}{=}\sup_{\alpha_k\leq u \leq \alpha_{k+1}}U(u)- \hat U_n(u)$, we get that 
\beqan
\mu&\stackrel{\cL}{\geq}&\sum_{k=0}^{K-1} \int_{\underline{x}_{\alpha_k}}^{\underline{x}_{\alpha_{k+1}}}  \phi\Big(\hat F_n(x) + \sup_{\alpha_k\leq u \leq \alpha_{k+1}}U(u)- \hat U_n(u)\Big)dx\\						
&\geq&\sum_{k=0}^{K-1}\int_{\underline{x}_{\alpha_k}}^{\underline{x}_{\alpha_{k+1}}}  \phi\Big(\hat F_n(x) + \tilde \epsilon_{[\alpha_k,\alpha_{k+1}]}(n,\delta_k)\Big)dx\,.
\eeqan
where the last inequality holds on an event of probability higher than $1- \sum_{k=1}^K\delta_k$.
Last, using the local Lipschitz property, we get
\beqan
\Pr\bigg( \mu < \int_\cX \phi\big(\hat F_n(x)\big) - \sum_{j=1}^J\ell_j \sum_{k=k_{j-1}}^{k_{j}-1} (\underline{x}_{\alpha_{k+1}}-\underline{x}_{\alpha_{k}})\tilde\epsilon_{[\alpha_k,\alpha_{k+1}]}(n,\delta_k)\bigg)\leq \sum_{k=1}^K\delta_k\,.
\eeqan 
We conclude using the assumption that $F$ has $\gamma$-Lipshitz quantile function.
\end{myproof}

\section{Time-uniform concentration inequalities}\label{sec:TU}

In this section, we now provide an extension of the previous result and focus on the number of samples $n$.
The previous result provide a confidence bound valid with high probability for each $n$.
In some situations, one way want to have a  high probability control valid \emph{simultaneously} for all $n$ in a given range, or even simultaneously all $n\in\Nat$. In order to derive such bounds, classical techniques consists of using (1) a union bound argument,
 (2) a geometric time-peeling argument together with Doob's maximal inequality for sub-martingales, or (3) a method of mixture (Laplace method) for specific distributions. These techniques lead to different bounds, the union bound technique being the simplest yet yielding the largest time-uniform  confidence bands.
 In the following, we provide a version of the geometric time-peeling argument for the control of the supremum CDF. One difficulty is that a Martingale cannot be easily built in this case, and hence we replace the use of Doob's maximal inequality with a weaker reflection inequality that can be traced back at least to \cite{james1975functional}.

We first show below a slight extension of 
\cite[Inequality 13.2.1]{shorack2009empirical} (the result from
\cite{shorack2009empirical} itself originates from, and slightly extends that of  \cite{james1975functional}). The proof of this result easily follows 
by looking at inequality (a) p. 513 in the proof of \cite[Inequality 13.2.1]{shorack2009empirical} and thus is not reproduced here.

\begin{lemma}[Reflection inequality]\label{lem:refinequality}
	Let $\eta\!>\!1, \lambda\!>\!0$ and then $\underline{u},\overline{u}\!\in\! [0,1]$.
	Let $c\!\in\!(0,1)$ be such that $\frac{(C-1)(1-c)^2\lambda^2}{C(\eta\!-\!1)}\geq \sup_{x\in[\uu,\ou]}x(1-x)$ for some $C>1$	Then,
	for all integers $n_1,n_2\!\in\!\Nat$ such that $n_2\!\leq n_1\!\eta$, it holds
	\beqan
	\Pr\bigg(  \max_{ n_1\leq n\leq n_2}  \sup_{x \in [\uu,\ou]} \sqrt{n}\Big(U_n(x) -U(x)\Big) > \lambda \bigg) &\leq& 
	C	\Pr\bigg( \sup_{x \in [\uu,\ou]} \sqrt{n_2}\Big(U_{n_2}(x) -U(x)\Big) >  \frac{c}{\sqrt{\eta}}\lambda\bigg)\,.\\
	\Pr\bigg(  \max_{ n_1\leq n\leq n_2}  \sup_{x \in [\uu,\ou]} \sqrt{n}\Big(U(x) -U_n(x)\Big)  > \lambda \bigg) &\leq& 
	C	\Pr\bigg( \sup_{x \in [\uu,\ou]} \sqrt{n_2}\Big(U(x)-U_{n_2}(x)\Big) >  \frac{c}{\sqrt{\eta}}\lambda\bigg)\,.\\
	\eeqan   
\end{lemma}

We deduce from this key result the following maximal inequality
\begin{corollary}[Maximal inequality]\label{cor:maxineq}
	Let $\eta\!>\!1$, $C\!>\!1$, $\underline{u},\overline{u}\!\in\![0,1]$ then $q \!=\! \sup_{x\in[\uu,\ou]}x(1-x)$.
	For any $n_1,n_2\!\in\!\Nat$ such that $n_2\!\leq\! \eta n_1$
	and for any $\epsilon>0$ such that 
	$\epsilon>\sqrt{\frac{Cq(\eta\!-\!1)}{(C-1)n_1}}$, it comes	
	\beqan
	\Pr\bigg(  \max_{ n_1\leq n\leq n_2}  \sup_{x \in [\uu,\ou]} U_n(x) \!-\!U(x)> \epsilon \bigg) &\leq&
	C\Pr\bigg( \sup_{x \in [\uu,\ou]} U_{n_2}(x) \!-\!U(x)
	> \sqrt{\frac{n_1}{n_2 \eta}}\Big(\epsilon  \!-\! \sqrt{\frac{Cq(\eta\!-\!1)}{(C-1)n_1}}\Big)\bigg)\,.\\	
	\Pr\bigg(  \max_{ n_1\leq n\leq n_2}  \sup_{x \in [\uu,\ou]} U(x) \!-\!U_n(x)> \epsilon \bigg) &\leq&
	C\Pr\bigg( \sup_{x \in [\uu,\ou]} U(x) \!-\!U_{n_2}(x)
	> \sqrt{\frac{n_1}{n_2 \eta}}\Big(\epsilon  \!-\! \sqrt{\frac{Cq(\eta\!-\!1)}{(C-1)n_1}}\Big)\bigg)\,.\\
	\eeqan	
\end{corollary}

\begin{myproof}{of Corollary~\ref{cor:maxineq}}
	First, we successively derive 
	\beqan
	\Pr\bigg(  \max_{ n_1\leq n\leq n_2}  \sup_{x \in [\uu,\ou]} U_n(x) -U(x)> \epsilon \bigg) &\leq&
	\Pr\bigg(  \max_{ n_1\leq n\leq n_2}  \sup_{x \in [\uu,\ou]} \frac{\sqrt{n}}{\sqrt{n_1}}\Big(U_n(x) -U(x)\Big)> \epsilon \bigg)\\
	&=&\Pr\bigg(  \max_{ n_1\leq n\leq n_2}  \sup_{x \in [\uu,\ou]} \sqrt{n}\Big(U_n(x) -U(x)\Big)> \sqrt{n_1}\epsilon \bigg)\\
	&\leq&C\Pr\bigg( \sup_{x \in [\uu,\ou]} \sqrt{n_2}\Big(U_{n_2}(x) -U(x)\Big)
	>  \sqrt{\frac{n_1}{\eta}}\epsilon c\bigg)\,.
	\eeqan
	Note that in the first line, we considered the event that
	$\sup_x U_n(x)-U(x)>0$. Indeed the complementary event does not intersect the event of interest.
	Now, we choose for $c\in(0,1)$ the maximal value such that 
	$\frac{n_1(C-1)(1-c)^2\epsilon^2}{C(\eta\!-\!1)}\geq \sup_{x\in[\uu,\ou]}x(1-x)$, that is 
	$c=1-\frac{\sqrt{C(\eta\!-\!1)q}}{\sqrt{(C-1)n_1}\epsilon}$, 
	provided that $\sqrt{\frac{C(\eta\!-\!1)q}{(C-1)n_1}}<\epsilon$.
	Reorganizing the terms yields the conclusion. We proceed similarly for the second inequality.
\end{myproof}

We are now ready to prove Theorem~\ref{thm:TuDKW}.
To this end, we combine Corollary~\ref{cor:maxineq} together with the local DKW inequality, on top of the standard geometric time-peeling technique.

\begin{theorem}[Time-Uniform local DKW inequality]\label{thm:TuDKW}
	Let $n\!\in\!\Nat$, and consider any random stopping time $N_n$ a.s. upper-bounded by $n$. 
	Let us introduce a function $\epsilon_{[\uu,\ou]}$ such that
	\beqan
	\forall \delta\!\in\!(0,1), \forall n\!\in\!\Nat,\quad
	\Pr\bigg( \sup_{x \in [\uu,\ou]} U_{n}(x) -U(x)> \epsilon_{[\uu,\ou]}\Big(n,\delta\Big)\bigg)\leq \delta\,.
	\eeqan 
	Then for all $\eta>1$, for all $\delta_n\in(0,1)$,  and $C>1$, it holds
	\beqan
	\Pr\bigg( \!\sup_{x \in [\uu,\ou]} U_{N_n}(x)\! -\!U(x)>
	\frac{1}{\sqrt{N_n\!-\! (\eta\!-\!1)}}\Big[\eta\sqrt{  N_n}\epsilon_{[\uu,\ou]}\bigg(\!N_n,\frac{\delta_n/C}{\big\lceil \frac{\ln(n)}{\ln(\eta)}\big\rceil}\!\bigg) + \sqrt{\!\frac{C}{C\!-\!1}q\eta(\eta\!-\!1)}\Big]	
	\bigg)\leq \delta_n\,.
	\eeqan
	Likewise, if $\tilde \epsilon_{[\uu,\ou]}$ controls 
	$\sup_{x \in [\uu,\ou]} U(x) -U_n(x)$, a similar inequality holds replacing  $U_n(x) -U(x)$ with $U(x) -U_n(x)$
	and $\epsilon_{[\uu,\ou]}$ with
	$\tilde \epsilon_{[\uu,\ou]}$.
\end{theorem}

\begin{corollary}[Time-Uniform global DKW inequality]\label{cor:TUDKW}
	In particular for $[\uu,\ou]=[0,1]$, 
	choosing $C=2$, and using that $q\leq 1/4$ it comes  $\forall\delta\in(0,0.5)$,
	\beqan
	\Pr\bigg( \sup_{x \in [0,1]} U_{N_n}(x) -U(x)>
	\frac{1}{\sqrt{2(N_n- (\eta\!-\!1))}}\Big[\eta\sqrt{ \ln\bigg(\bigg\lceil \frac{\ln(n)}{\ln(\eta)}\bigg\rceil\frac{2}{\delta}\bigg)} + \sqrt{\eta(\eta\!-\!1)}\Big]
	\bigg)\leq \delta\,.
	\eeqan
\end{corollary}
Corollary~\ref{cor:TUDKW} is stated for convenience, to show an explicit formula that can be used to control uniform deviations uniformly over time.
This result should be compared to the term $\sqrt{\ln(1/\delta)/2n}$ obtained for a single $n$ by application of Massart's inequality.
Note the $\log\log(n)$ scaling, compared to the $\log(n)$ term one would obtain from a simple union bound.

\begin{myproof}{of Theorem~\ref{thm:TuDKW}}
	Let us introduce $t_k = \lfloor\eta^k\rfloor$, for $k=0,\dots,K$ with $K = \big\lceil \frac{\ln(n)}{\ln(\eta)}\big\rceil$ (thus $n\leq t_K$), for some constant $\eta>1$.  Let  also $\epsilon$ be a non-increasing positive function on $\Real$, to be defined later.
	\beqan
	\lefteqn{
		\Pr\bigg( \sup_{x \in [\uu,\ou]} U_{N_n}(x) -U(x)> \epsilon(N_n) \bigg)}\\
	&\leq& \sum_{k=1}^K 
	\Pr\bigg( \exists t \in[t_{k\!-\!1},t_k \!-\!1] : \sup_{x \in [\uu,\ou]} \Big(U_{t}(x) -U(x)\Big)> \epsilon(t) \bigg)\\
	&\leq& \sum_{k=1}^K 
	\Pr\bigg( \max_ {t \in[t_{k\!-\!1},t_k \!-\!1]}\sup_{x \in [\uu,\ou]} \Big(U_{t}(x) -U(x)\Big)> \epsilon(t_k-1) \bigg)\\
	&\leq&  \sum_{k=1}^K 
	C\Pr\bigg( \sup_{x \in [\uu,\ou]} \Big(U_{t_k-1}(x) -U(x)\Big)> \sqrt{\frac{t_{k-1}}{(t_k -1) \eta}}\Big(\epsilon(t_k\!-\!1)  - \sqrt{\frac{Cq(\eta\!-\!1)}{(C\!-\!1)t_{k-1}}}\Big)\bigg)\\
	\eeqan
	In order to make use of local DKW inequality, let us choose the function $\epsilon$ such that 
	\beqan
	\sqrt{\frac{t_{k-1}}{(t_k -1) \eta}}\Big(\epsilon(t_k\!-\!1)  - \sqrt{\frac{Cq(\eta\!-\!1)}{(C\!-\!1)t_{k-1}}}\Big) \geq \epsilon_{[\uu,\ou]}(t_k\!-\!1,\delta_k)\,.
	\eeqan
	Indeed, this ensures by local DKW inequality that, provided that $\delta_k\in(0,1)$, 
	\beqan
	\Pr\bigg( \sup_{x \in [\uu,\ou]} \Big(U_{N_n}(x) -U(x)\Big)> \epsilon(N_n) \bigg) \leq \sum_{k=1}^K 
	C\delta_k\,,
	\eeqan
	and assuming $C>1$, we can further choose $\delta_k = \delta/(CK)$ for $\delta\in(0,1)$.
	That is, we want 
	\beqan
	\epsilon(t_k\!-\!1) \geq \frac{1}{\sqrt{t_{k-1}}}\bigg(	\sqrt{(t_k -1) \eta} \epsilon_{[\uu,\ou]}\Big(t_k\!-\!1,\frac{\delta}{CK}\Big) + \sqrt{\frac{C}{C\!-\!1}q(\eta\!-\!1)}\bigg)
	\eeqan
	We now use the fact that
	\beqan 
	t_{k}-1 &=& \lfloor\eta^k\rfloor-1\leq \eta^{k-1} \eta -1
	= \lfloor\eta^{k-1}\rfloor \eta + \Big(\eta^{k-1}-\lfloor\eta^{k-1}\rfloor\Big)\eta\!-\!1\\
	&\leq& t_{k-1} \eta + (\eta\!-\!1)\,.
	\eeqan
	Hence, we conclude the proof by choosing the function 
	\beqan
	\epsilon(t) &=& \sqrt{\frac{\eta}{t - (\eta\!-\!1)}}\Big(\sqrt{t \eta} \epsilon_{[\uu,\ou]}(t,\frac{\delta}{CK})  + \sqrt{\frac{C}{C\!-\!1}q(\eta\!-\!1)}\Big)\\
	&=&\frac{1}{\sqrt{t - (\eta\!-\!1)}}\Big(\eta\sqrt{t} \epsilon_{[\uu,\ou]}(t,\frac{\delta}{CK}) + \sqrt{\frac{C}{C\!-\!1}q\eta(\eta\!-\!1)}\Big)\,.
	\eeqan
\end{myproof}

\paragraph{An application to the cumulative control of sequential deviations}

Using Theorem~\ref{thm:TuDKW}, we deduce the following control, together  with tuning recommendations for the practitioner.
\begin{corollary}[Time uniform cumulative error control]\label{cor:TUcumulative}
	Let  $q = \sup_{x\in[\uu,\ou]}x(1-x)$.  Let $(\eta_t)_t$ be a decreasing sequence converging to $1$, and $(\delta_t)_t\subset(0,1)$. Then, it holds 
	\beqan
	\sum_{t=1}^T 
	\Pr\bigg( \sup_{x \in [\uu,\ou]} U_{N_t}(x) -U(x)&>&
	\frac{1}{\sqrt{N_t\!-\! (\eta_t\!-\!1)}}\Big[\eta_t\sqrt{  N_t}\epsilon_{[\uu,\ou]}\bigg(N_t, \delta_t\bigg) + \sqrt{2q\eta_t(\eta_t\!-\!1)}\Big]	
	\bigg)\\
& \leq& 2\sum_{t=1}^T \bigg\lceil\frac{\ln(t)}{\ln(\eta_t)}\bigg\rceil\delta_t,
	\eeqan
	In particular, provided that (a) $\ln(\eta_t)^{-1} = O(\ln^\alpha(t))$ for some $\alpha<1$, choosing $\delta_t$  such that $\ln(t)^2\delta_t=O(1/t)$ ensures
	that the right hand side term is $o(\ln(T))$. This quantity is also 
	\beqan
	o(\log(T)) &\text{ for (b)}&\eta_t= \frac{f(t)+1}{f(t)}\text{ and }\delta_t=\exp(-f(t))\text{ where}
	f(t)= \log(t)+\xi\log\log(t)\text{ with } \xi>2\,.\\	
	\leq 1 &\text{ for (c)}& \delta_t = \bigg\lceil\frac{\ln(t)}{\ln(\eta_t)}\bigg\rceil^{-1} \frac{1}{2g(t)} \text{ where } 
	\sum_{t=1}^\infty \frac{1}{g(t)}\leq 1 \text{ provided that }\delta_t\leq 1.
	\eeqan

\end{corollary}
\begin{remark}
	The tuning of $\eta_t$ using $f(t)$ is the one suggested in \cite{KLUCBJournal} for the tuning of the \texttt{KL-UCB} algorithm, and this is the one we suggest using in practice.
	In the last case c), the condition $\delta_t\leq 1$ constrains the choice of $\eta_t$ that cannot be too small. In particular, it cannot converge too fast towards $1$. Some classical choice for $g$ include $g(t)=3t^{3/2}$, $g(t)=t(t+1)$ or $g(t)=(t+1)\log^2(t+1)/\log(2)$. See Appendix~\ref{app:other} for other possible choices for $g$. 
	\end{remark}

\begin{remark}For comparison, note that a union bound argument yields the alternative bound
		\beqan
	\sum_{t=1}^T 
	\Pr\bigg( \sup_{x \in [\uu,\ou]} U_{N_t}(x) -U(x)> \epsilon_{[\uu,\ou]}(N_t,\delta_t)\bigg) \leq \sum_{t=1}^T t \delta_t\,,
	\eeqan
	which suggests choosing e.g. $\delta_t=(tg(t))^{-1}$. A classical choice is e.g. $\delta(t)=(t^2(t+1))^{-1}$.
	\end{remark}
\begin{myproof}{of Corollary~\ref{cor:TUcumulative}}
	We apply Theorem~ \ref{thm:TuDKW} with $C=2$ in order to get the claimed inequality.
	
	The first claim follows from the fact that $\bigg\lceil\frac{\ln(t)}{\ln(\eta_t)}\bigg\rceil\delta_t = O(t^{-1}\ln(t)^{1+\alpha-2})=o(t^{-1})$ under the considered assumption.
	
Regarding the second claim, for the choice $\eta_t= \frac{f(t)+1}{f(t)}$ and $\delta_t=\exp(-f(t))$ where
	$f(t)= \log(t)+\xi\log\log(t)$, we obtain that
	\beqan
	\sum_t \bigg\lceil\frac{\ln(t)}{\ln(\eta_t)}\bigg\rceil\delta_t &=& e\sum_t  \frac{1}{t\log^\xi(t)} \lceil \log^2(t) + \xi\log(t)\log\log(t)\rceil\\
	&\leq&e\sum_t \frac{1}{t\log^{\xi-2}(t)} + e\sum_t \frac{\log\log(t)}{t\log^{\xi-1}(t)} 
	+ e\sum_t \frac{1}{t\log^{\xi}(t)}\,.
	\eeqan
	Hence, provided that $\xi>2$, this sum is $o(\log(T))$.	  The last claim  is direct; note that the condition $\delta_t\leq 1$ constrains the choice of  admissible $\eta_t$.
\end{myproof}

\section*{Acknowledgement}
This work has been supported by CPER Nord-Pas-de-Calais/FEDER DATA Advanced data science
and technologies 2015-2020, the French Ministry of Higher Education and Research, Inria, the French Agence Nationale de la Recherche (ANR) under grant ANR-16-CE40-0002 (the BADASS project), the MEL, the I-Site
ULNE regarding project R-PILOTE-19-004-APPRENF, and the Inria A.Ex. SR4SG project.

\vspace{-4mm}

\bibliography{biblio.bib,biblioCVAR.bib}

\newpage
\appendix
\section{Proofs of the main results}\label{app:concentration}

\begin{myproof}{of Lemma~\ref{lem:concfirst}}

	We let $[n]=\{1,\dots,n\}$ for all $n\in\Nat_\star$.
	
	{\bf Left tail, step 1} Let first recall the following remark by Smirnov \cite[p.10]{smirnov1944approximate}, showing that if $u_{(1)},\dots,u_{(n)}$ denotes the order samples received from the uniform distribution, then 
	\beqa
	\Pr(\sup_{u\in [0,1]} U_n(u)-u \leq \epsilon)
	&=& \Pr\Big(\forall k \in [n],\quad U_n(u_{(k)}) - \epsilon \leq u_{(k)}\Big)\nonumber\\
	&=&\Pr\Big(\forall k \in [n],\quad k/n - \epsilon \leq u_{(k)}\Big)\nonumber\\
	&=& n! \!\int\!\dots\!\int\!\ind\bigg\{ 0\!\leq\! u_1 \!\leq\! \dots u_n \!\leq\! 1; \forall k, u_k \geq k/n-\epsilon\bigg\} du_1\dots du_n\,.\label{eqn:smirnov1}
	\eeqa
		When restricting the supremum to $[\alpha,\beta]$, this equality needs to be modified.
		First of all, it holds that
		\beqan
		\sup_{u\in [\alpha,\beta]} U_n(u)-u = \max\bigg\{U_n(v)-v: v \in \{\alpha\}\cup \{ u_{(1)},\dots,u_{(n)} \}\cap[\alpha,\beta]\bigg\}\,.
		\eeqan
		Hence, we deduce that 
		\beqan
		\lefteqn{\bigg\{\sup_{u\in [\alpha,\beta]} U_n(u)-u \leq \epsilon
		\bigg\}}\\ &=& 	\bigcap_{k \in[n]}	\bigg\{u_{(k)} \in[\alpha,\beta]\implies U_n(u_{(k)}) - \epsilon \leq u_{(k)} \bigg\} \cap \bigg\{ U_n(\alpha)\leq \epsilon+\alpha\bigg\}\\		
		&=& 	\bigcap_{k \in [n]}		\bigg\{u_{(k)} \in[\alpha,\beta]\implies k/n - \epsilon \leq u_{(k)} \bigg\} \cap \bigg\{\sum_{k=1}^n\indic{u_{(k)} \leq \alpha}\leq n(\epsilon+\alpha)\bigg\}\\
		 &=& 	\bigcap_{k \in [n]}		\bigg\{u_{(k)} \in[\alpha,\beta]\implies  k/n - \epsilon \leq u_{(k)} \bigg\} \cap \bigg\{ u_{(k)} \leq \alpha < u_{(k+1)} \implies k/n\leq \epsilon+\alpha\bigg\}\,,
		\eeqan
		where we introduced the term $u_{(n+1)}=1$ in the last line and used that $\epsilon+\alpha\geq 0$ to exclude the term $u_{(0)}=0$. Using the distribution of the order statistics, we deduce that
		\beqan
		\lefteqn{\Pr\Big(\sup_{u\in [\alpha,\beta]} U_n(u)-U(u) \leq \epsilon\Big)
		}\\	\! &=&\!\!n! \!\int\!\dots\!\int\!\ind\bigg\{ 0\!\leq\!  u_1 \!\leq \!\dots u_n \!\leq\! 1; \forall k\!\in\! [n],
	\begin{cases}
 \text{if } u_k \!\in\! [\alpha,\beta]& \text{ then }u_k \geq \frac{k}{n}-\epsilon\\
 \text{if } k\!>\!n(\epsilon\!+\!\alpha) & \text{ then } \alpha \!\notin\! [u_k,u_{k+1})\\
	\end{cases}\bigg\} du_1\dots du_n.
		\eeqan

	{\bf Left tail, step 2} Following \cite{smirnov1944approximate}, we introduce the notation $t_k=u_{n-k+1}$, constant $\gamma_k= (n-k+1)/n-\epsilon$ (non-negative for $k\leq n-\lfloor \epsilon n \rfloor$) as well as $\beta_k = \min(\gamma_k,\beta)$. We thus have  the following rewriting
	\beqan
	\bigg\{ t \in [0,1]:
	\text{if } t \in [\alpha,\beta] \text{ then }t \geq \gamma_k\bigg\} 
	&=& [0,\alpha]\cup[\min(\gamma_k,\beta),1]=[0,\alpha]\cup[\beta_k,1]\,,
	\eeqan
	which further reduces to $[0,1]$ when $\gamma_k\leq \alpha$.  
	We let $n_{\alpha,\epsilon}=n(1\!-\!\epsilon\!-\!\alpha)$, $\overline{n}_{\alpha,\epsilon}=\lceil n_{\alpha,\epsilon} \rceil$ and remark that  $\gamma_k\leq \alpha$ iff $k\geq n_{\alpha,\epsilon}+1$ . Also, $\gamma_k\leq \alpha$ as soon as
	$k>\overline{n}_{\alpha,\epsilon}$.
	Let us also note that $n-k+1>n(\epsilon+\alpha)$ iff $k<n_{\alpha,\epsilon}+1$, and that $\overline{n}_{\alpha,\epsilon}<n_{\alpha,\epsilon}+1$.
	This means in particular that if $k\geq \overline{n}_{\alpha,\epsilon}+1$, then both conditions in the integral vanish (so contribute to $1$ in the integral).
	\beqan
	\lefteqn{\Pr\Big(\sup_{u\in [\alpha,\beta]} U_n(u)-U(u) \leq \epsilon\Big)
	}\\	\! &=&\!\!n! \!\int\!\dots\!\int\!\!\ind\bigg\{ \!0\!\leq\! t_n \!\leq \!\dots t_1 \!\leq\! 1; \forall k\!\in\![n], t_k \!\in\! [0,\alpha]\!\cup\![\beta_k,1],\forall k\!\leq\! \overline{n}_{\alpha,\epsilon}, \alpha \!\notin\! [t_k,t_{k-1})\!\bigg\} dt_n\dots dt_1\\
	&=&\!\!n!\! \!\int\!\dots\!\int\!\!\mathbb{I}\bigg\{ \!0\!\leq \!t_{\overline{n}_{\alpha,\epsilon}} \!\leq\! \dots t_1\! \leq\! 1;\\
	&& \qquad\qquad\forall k\!\leq\! \overline{n}_{\alpha,\epsilon}, t_k \!\in\! [0,\alpha]\!\cup\![\beta_k,\!1] \text{ and }\alpha \!\notin\! [t_k,t_{k-1})\!\!\bigg\} J_{n\!-\!\overline{n}_{\alpha,\epsilon}}\!(t_{\overline{n}_{\alpha,\epsilon}}) dt_{\overline{n}_{\alpha,\epsilon}\!}\!\dots dt_1\,,
	\eeqan 	where we integrated out all terms for $k> \overline{n}_{\alpha,\epsilon}$ into  the short-hand notation $J_k(x)=\frac{x^k}{k!}$.
	In the integral, we note that 
	if $t_k\leq \alpha$, then so must be all terms $t_{k'}$ for $k'\geq k$.
	
	We now proceed with integration. Starting with $t_1$, we see that if $t_1\leq \alpha$, then this implies $\alpha \in [t_1,t_0]$.
	Hence, the corresponding terms are $0$, and it remains to integrate $t_1$ on $(\alpha,1]$, that is on $[\beta_1,1]$.
	
Regarding $t_2$, if $t_2\leq \alpha<t_1$, this contradicts
$\alpha \notin [t_2,t_1]$, hence it remains it remains to integrate $t_2$ on $(\alpha,1]$, that is on $[\beta_2,1]$.
Proceeding similarly, for all $k\leq \overline{n}_{\alpha,\epsilon}$ we obtain that 
\beqan
\Pr\Big(\sup_{u\in [\alpha,\beta]} U_n(u)-U(u) \leq \epsilon\Big) 
&=&
	n! \int_{\beta_1}^1\!\int_{\beta_2}^{t_1}
\!\dots  \!\int_{\beta_{\overline{n}_{\alpha,\epsilon}}}^{t_{\overline{n}_{\alpha,\epsilon}-1}}\!\!J_{n-\overline{n}_{\alpha,\epsilon}}(t_{\overline{n}_{\alpha,\epsilon}})dt_{\overline{n}_{\alpha,\epsilon}}\!\dots\! dt_1\!\,.\\
\eeqan
	
	In order to compute the multiple integral, similarly to \cite{smirnov1944approximate},  we make use of the following variant of the Taylor expansion 
	\beqan
	f(x) &=& f(a_1) +  \sum_{\ell=1}^{n-k-1}f^{(\ell)}(a_{\ell+1})\int_{a_1}^{x}\int_{a_2}^{t_1}\dots\int_{a_k}^{t_{\ell-1}}dt_1\dots dt_\ell\\
	&& + \int_{a_1}^{x}\int_{a_2}^{t_1}\dots\int_{a_{n-k}}^{t_{n-k-1}}f^{(n-k)}(t_{n-k})dt_1\dots dt_{n-k}\,,
	\eeqan
	which, using the $I_k$ notation, yields the following form
	\beqan
	\int_{a_1}^{x}\int_{a_2}^{t_1}\dots\int_{a_{n-k}}^{t_{n-k-1}}f^{(n-k)}(t_{n-k})dt_1\dots dt_{n-k} = f(x)- f(a_1)-\sum_{\ell=1}^{n-k-1}f^{(\ell)}(a_{\ell+1})I_\ell(x;a_1,\dots,a_\ell)\,.
	\eeqan
	
	This is applied to the function $f(x)=x^n$, $k=n-\overline{n}_{\alpha,\epsilon}$. Indeed, we then get
	$f^{(n-k)}(x) =\frac{n!x^k}{k!}= \frac{n!x^{n-\overline{n}_{\alpha,\epsilon}}}{(n-\overline{n}_{\alpha,\epsilon})!}=n!J_{n-\overline{n}_{\alpha,\epsilon}}(x)$. This in turns
yields
	\beqan
	\Pr\Big(\sup_{u\in [\alpha,\beta]} U_n(u)-U(u) \leq \epsilon\Big)
	&=&1- \beta_{1}^{n}-\sum_{\ell=1}^{\overline{n}_{\alpha,\epsilon}-1} \frac{n!}{(n-\ell)!}\beta_{\ell+1}^{n-\ell}I_{\ell}(1;\beta_1,\dots \beta_\ell) \\
	&=&1 -\sum_{\ell=0}^{\overline{n}_{\alpha,\epsilon}-1} \binom{n}{\ell}\beta_{\ell+1}^{n-\ell}\ell!I_{\ell}(1;\beta_1,\dots \beta_\ell) \,,
	\eeqan
	using the convention that  $I_0(x;\emptyset)=1$. 	This completes the proof regarding the left tail concentration.
	
	{\bf Right tail} 	We proceed similarly for the right tail. 
	First, using our notation we note that 
		\beqan
	\sup_{u\in [\alpha,\beta]} U(u)-U_n(u) = \max\bigg\{\lim_{u\to v;u<v} U(u)-U_n(u): v \in  \{\beta\}\cup\{u_{(1)},\dots,u_{(n)}\}\cap(\alpha,\beta]\bigg\}
	\eeqan
			To be more precise, we let $(\eta_k)_{k\in[n]}>0$ be  arbitrary small constants. We also let $\eta_0<\min_{k\in[n]}\eta_k$ and define
			$\overline{\eta}= \max_{k\in[n]}\eta_k$. We further introduce, for each $k\!\in\![n]$, $u^-_{(k)}$ such that 
		$u_{(k)} - \eta = u^-_{(k)} < u_{(k)}$, and $\beta^{-}$ such that 	$\beta-\eta_0 = \beta^{-} < \beta$. Then, we introduce the notation
		\beqan
\sup^{\eta}_{u\in [\alpha,\beta]} U(u)-U_n(u) = \max\bigg\{U(v)-U_n(v): v \in  \{\beta^{-}\}\cup\{u_{(1)}^{-},\dots,u_{(n)}^{-}\}\cap(\alpha,\beta]\bigg\}\,.
\eeqan		
Before proceeding, we note that $\forall n\!\in\!\Nat,\lim_{\overline{\eta}\to0}\Pr(\min_{k\in[n]} u_{(k)} - u_{(k-1)}>\eta_k)= 1$
Indeed, it holds
\beqan
\Pr\big(\min_{k\in[n]} u_{(k)} - u_{(k-1)}\leq \eta_k\big) &=& \Pr\big(\exists {k\in[n]},\,\, u_{(k)} - u_{(k-1)}\leq \eta_k\big)\\
&\leq& \Pr\big(\exists {i,j\in[n], i<j},\,\, |X_i-X_j|\leq\overline{\eta}\big)\\
&	\leq&  \frac{n(n-1)}{2} \Pr(|X_1-X_2|\leq\overline{\eta}) = n(n-1)\overline{\eta}\,.
\eeqan
In the following, we use a construction similar to that of \cite{kolmogorov1956skorokhod} for Skorokhod convergence.
Note that under the event that $\Omega_n = \big\{\min_{k\in[n]} u_{(k)} - u_{(k-1)}> \eta_k\big\}$ (where $u_{(0)}=0$) we have the following rewriting
	\beqan
	\lefteqn{\bigg\{\sup^{\eta}_{u\in [\alpha,\beta]} u-U_n(u) \leq \epsilon
		\bigg\}}\\ &=& 	\bigcap_{k \in[n]}	\bigg\{u_{(k)}^{-} \in[\alpha,\beta]\implies u_{(k)}^{-} - \epsilon \leq U_n(u_{(k)}^{-}) \bigg\} \cap \bigg\{ \beta^{-} -\epsilon\leq U_n(\beta^{-})\bigg\}\\		
	&=& 	\bigcap_{k \in [n]}		\bigg\{u_{(k)}^{-} \in[\alpha,\beta]\implies u_{(k)}^{-} - \epsilon \leq (k-1)/n  \bigg\} \cap \bigg\{\sum_{k=1}^n\indic{u_{(k)} \leq \beta^{-}}\geq n(\beta^{-}-\epsilon)\bigg\}\\
	&=& 	\bigcap_{k \in [n]}	\!	\bigg\{u_{(k)}^{-} \!\in\![\alpha,\beta]\!\implies\!  u_{(k)}^{-} \!-\! \epsilon \leq (k\!-\!1)/n \bigg\} \cap \bigg\{ u_{(k\!-\!1)}\! \leq\! \beta^{-} \!<\! u_{(k)} \!\implies\! (k\!-\!1)/n\geq \beta^{-}\!-\!\epsilon\bigg\}\,.
	\eeqan
	In the last line,  we used that $\beta^-=\beta-\eta_0$ and $\eta_0<\min_{k\in[n]}\eta_k$ to rewrite $\sum_{k=1}^n\indic{u_{(k)} \leq \beta^{-}}$, in terms of $u_{(k)}\! \leq\! \beta^{-} \!<\! u_{(k\!+\!1)}$. Then we shifted $k$ by $1$, and used the fact that $\beta\leq 1$ implies $1\geq \beta^{-}-\epsilon$ in order to exclude the term $u_{(n+1)}=1$. 

	We let $\tilde \beta=1-\beta, \tilde \alpha=1-\alpha$, $\tau_k=1-u_k$
	and introduce for all $k$ the constant $\rho_k= 1-\epsilon-k/n$  (non-negative for all $k\leq n(1-\epsilon)$ as well as
	$\tilde \alpha_k = \min(\rho_{k-1},\tilde \alpha)$. We also let $\tilde \beta^{+}=1-\tilde\beta^{-}$. Finally, we 
		let $\tau_{k}<_k \tau_{k-1}$ if and only if $\tau_k+\eta_k<\tau_{k-1}$.
	Using the distribution of the order statistics together with these notations,we then naturally study the quantity $\lim_{\overline{\eta}\to0}
	\Pr^{\eta}\Big(\sup_{u\in [\alpha,\beta]} U(u)-U_n(u) \leq \epsilon\Big)$ where
	\beqan
	\lefteqn{\Pr^\eta\Big(\sup_{u\in [\alpha,\beta]} U(u)-U_n(u) \leq \epsilon\Big)=\Pr\Big(\sup^\eta_{u\in [\alpha,\beta]} U(u)-U_n(u) \leq \epsilon \cap \Omega_n\Big)
	}\\	\! &\!\!=&\!\!n!\!\! \int\!\!\!\dots\!\!\int\!\!\ind\bigg\{\!\Omega_n\!\cap \!0\leq \! u_1\! \leq\! \dots u_n\! \leq\!\! 1; \forall k\!\in\! [n],\!
	\begin{cases}
		\text{if } u_k^{-} \!\in\! [\alpha,\beta]& \!\!\text{then }u_k^{-} \!\leq\! \frac{k-1}{n}\!+\!\epsilon\\
		\text{if } k\!-\!1\!<\!n(\beta^{-}\!-\!\epsilon) & \!\!\text{then } \beta^{-} \!\!\notin\! [u_{k-1},u_{k})\\
	\end{cases}\!\!\bigg\} du_1\dots du_n.\\
&\!\!=&\!\!n!\! \!\int\!\! \! \dots\!\! \int\!\! \ind\bigg\{ \!0\!\leq\!  \tau_n \!<_n \!\dots \tau_1 \!<_1 \!1; \forall k\!\in\! [n],
\begin{cases}
	\text{if } \tau_k^{+} \!\in\! [\tilde \beta,\tilde \alpha]& \text{ then }\tau_k^{+}\! \geq \! \rho_{k-1}\\
	\text{if } k\!-\!1<n(\beta^{-}\!-\!\epsilon) & \text{ then } \tilde \beta^{+} \!\notin\! [\tau_{k},\tau_{k-1})\\
\end{cases}\!\bigg\} d\tau_n\dots d\tau_1.
	\eeqan

	Now, $[0,\tilde \beta]\cup [\min(\rho_{k-1},\tilde \alpha),1]$ reduces to $[0,1]$ when
	$\rho_{k-1}\leq \tilde \beta$.  	
	We let $n_{\beta,\epsilon}=n(\beta-\epsilon)$, $\overline{n}_{\beta,\epsilon}=\lceil n_{\beta,\epsilon} \rceil$ and remark that  $\rho_{k-1}\leq \tilde \beta$ iff $k-1\geq n_{\beta,\epsilon}$. 	We first deal with the case when $n_{\beta,\epsilon}\in\Nat$.
	 In this situation, provided that $\eta_0$ is sufficiently small, then $\overline{n}_{\beta^{-},\epsilon}=\overline{n}_{\beta,\epsilon}=n_{\beta,\epsilon}$  and also, $k-1<n_{\beta^{-},\epsilon}$ iff $k\leq \overline{n}_{\beta^{-},\epsilon}$.
	 If $k> \overline{n}_{\beta,\epsilon}$, then $\rho_{k-1}\leq \tilde \beta$ and both restrictions disappear in the integral.	 
	In the general situation when $n_{\beta,\epsilon}\notin\Nat$, then 
	 $\overline{n}_{\beta,\epsilon}>n_{\beta,\epsilon}$ and 
	 $\overline{n}_{\beta^{-},\epsilon}=\overline{n}_{\beta,\epsilon}$.
	  Also, $k-1<n_{\beta^{-},\epsilon}$ iff $k\leq \overline{n}_{\beta^{-},\epsilon}$.
	  If $k> \overline{n}_{\beta,\epsilon}$  then $\rho_{k-1}\leq \tilde \beta$ and again restrictions disappear in the integral.
	We deduce that provided that $\overline{\eta}$ is sufficiently small, 
	\beqan
	\lefteqn{\Pr^\eta\Big(\sup_{u\in [\alpha,\beta]} U(u)-U_n(u) \leq \epsilon\Big)
	}\\\!\!	 &\!\!=\!\!&\!\!n! \!\int\!\dots\!\int\!\!\ind\bigg\{ 0 \!\leq\! \tau_n\! <_n\! \dots\! \tau_1\!<_1\! 1; \forall k\!\in\![n], \tau_k \!\in\! [0,\tilde \beta]\cup [\tilde \alpha_k,1],
\forall k\!\leq\!\overline{n}_{\beta^{-},\epsilon}, \tilde \beta^{+} \notin [\tau_{k},\tau_{k-1})  \bigg\} d\tau_n\dots d\tau_1\\
	\!\!&\!\!=\!\!&\!\!n!\!\! \int\!\!\dots\!\!\int\!\!\mathbb{I}\bigg\{ 0\!\leq\! \tau_{\overline{n}_{\beta,\epsilon}}\!\! <_{{\overline{n}_{\beta,\epsilon}}}\! \dots\! \tau_1\! <_1\! 1;\\
	&&\qquad\qquad \forall k \!\leq\! \overline{n}_{\beta,\epsilon},  \tau_k \!\in\! [0,\tilde \beta]\!\cup\![\tilde \alpha_k,1] \text{ and } 	\tilde \beta^{+} \!\!\notin\! [\tau_{k},\!\tau_{k-1})\! \bigg\} J ^{\eta}_{n-\overline{n}_\beta}(\!\tau_{\overline{n}_{\beta,\epsilon}}\!) d\tau_{\overline{n}_{\beta,\epsilon}}\!\dots\! d\tau_1,
	\eeqan 	where we integrated out all terms for $k> \overline{n}_{\beta,\epsilon}$ in  the term
	$J^{\eta}_m(x)$, that satisfies $\lim_{\overline{\eta}\to0}J^{\eta}_m(x)= \frac{x^m}{m!}$.

		We now proceed with integration. Starting with $\tau_1$, we see that if $\tau_1\leq \tilde \beta$, then this implies $\tilde \beta \in [\tau_1,\tau_0)$.
		The case when $\tilde \beta^+\geq \tau_0=1$, that is $\beta^{-}\leq 0$ is excluded by the assumption that $\beta>0$.
		Hence, this in turns implies $\tilde \beta^{+} \in [\tau_1,\tau_0)$, provided that $\eta_0<\tau_0-\tilde\beta=\beta$.
		Since this event is excluded by the indicator function, the corresponding terms are $0$, and it remains to integrate $\tau_1$ on $(\tilde \beta,1]$, that is on $[\tilde \alpha_1,1]$.	
	Regarding $\tau_2$, if $\tau_2\leq \tilde \beta^{+}<\tau_1$, this contradicts
	$\tilde \beta^{+} \notin [\tau_2,\tau_1)$, hence it remains to integrate $\tau_2$ on $(\tilde \beta^{+},1]$, that is on $[\tilde \alpha_2,1]$.
	We proceed similarly for all $k\leq \overline{n}_{\alpha,\epsilon}$. We obtain that for $\overline{\eta}$ sufficiently small, 
	\beqan
	\lefteqn{
	\Pr^\eta\Big(\sup_{u\in [\alpha,\beta]} U(u)-U_n(u) \leq \epsilon\Big)}\\
	&=&
	n! \int_{\tilde\alpha_1}^1\!\int_{\tilde\alpha_2}^{\tau_1}
	\!\dots  \!\int_{\tilde\alpha_{\overline{n}_{\beta,\epsilon}}}^{\tau_{\overline{n}_{\beta,\epsilon}-1}}\!\!\ind\bigg\{ 0\!\leq\! \tau_{\overline{n}_{\beta,\epsilon}}\!\! <_{{\overline{n}_{\beta,\epsilon}}}\! \dots\! \tau_1\! <_1\! 1\bigg\}J^\eta_{n-\overline{n}_{\beta,\epsilon}}(\tau_{\overline{n}_{\beta,\epsilon}})d\tau_{\overline{n}_{\beta,\epsilon}}\!\dots\! d\tau_1\!\,.
	\eeqan

Now, we remark that $\lim_{\overline{\eta}\to0}\ind\bigg\{ 0\!\leq\! \tau_{\overline{n}_{\beta,\epsilon}}\!\! <_{{\overline{n}_{\beta,\epsilon}}}\! \dots\! \tau_1\! <_1\! 1\bigg\}= \ind\bigg\{ 0\!\leq\! \tau_{\overline{n}_{\beta,\epsilon}}\!\! \leq \dots\! \tau_1\! \leq 1\bigg\}$, and so
\beqan
\lim_{\overline{\eta}\to0}\Pr^\eta\Big(\sup_{u\in [\alpha,\beta]} U(u)-U_n(u) \leq \epsilon\Big)
&=&
n! \int_{\tilde\alpha_1}^1\!\int_{\tilde\alpha_2}^{\tau_1}
\!\dots  \!\int_{\tilde\alpha_{\overline{n}_{\beta,\epsilon}}}^{\tau_{\overline{n}_{\beta,\epsilon}-1}}\!\!J_{n-\overline{n}_{\beta,\epsilon}}(\tau_{\overline{n}_{\beta,\epsilon}})d\tau_{\overline{n}_{\beta,\epsilon}}\!\dots\! d\tau_1\!\,.
\eeqan

In order to compute the multiple integral, we resort to a Taylor expansion as for the Left tail, and deduce that 
	\beqan
\lim_{\overline{\eta}\to0}	\Pr^{\eta}\Big(\sup_{u\in [\alpha,\beta]} U(u)-U_n(u) \leq \epsilon\Big)
	&=&1 -\sum_{\ell=0}^{\overline{n}_{\beta,\epsilon}-1} \binom{n}{\ell}\tilde \alpha_{\ell+1}^{n-\ell}\ell!I_{\ell}(1;\tilde \alpha_1,\dots \tilde \alpha_\ell) \,.
	\eeqan
	
It remains to note that $\lim_{\overline{\eta}\to0} \Pr\Big(\sup^\eta_{u\in [\alpha,\beta]} U(u)-U_n(u) \leq \epsilon \cap \Omega_n^c\Big) \leq \lim_{\overline{\eta}\to0} \Pr\Big(\Omega_n^c\Big) =0$, and thus 
\beqan
\lim_{\overline{\eta}\to0} \Pr\Big(\sup^\eta_{u\in [\alpha,\beta]} U(u)-U_n(u)\Big)&=&\lim_{\overline{\eta}\to0} \Pr\Big(\sup^\eta_{u\in [\alpha,\beta]} U(u)-U_n(u)\cap \Omega_n\Big)\\
&=&\lim_{\overline{\eta}\to0}\Pr^\eta\Big(\sup_{u\in [\alpha,\beta]} U(u)-U_n(u) \leq \epsilon\Big)\,.
\eeqan
This shows that the limit of $\Pr\Big(\sup^\eta_{u\in [\alpha,\beta]} U(u)-U_n(u)\Big)$ indeed exists and hence gives the value of
$\Pr\Big(\sup_{u\in [\alpha,\beta]} U(u)-U_n(u)\Big)$.

\end{myproof}

\begin{myproof}{of Theorem~\ref{thm:concLeft}}
	We now compute for $\ell\geq 1$ the quantity
	\beqan
	I_\ell(x;\beta_1,\dots,\beta_\ell) = \int_{\beta_1}^{x}\int_{\beta_2}^{t_1}\dots\int_{\beta_\ell}^{t_{\ell-1}}dt_\ell\dots dt_1\,.
	\eeqan	
	We further let $n_\beta=n(1-\beta-\epsilon)$, $ \underline{n}_\beta=\lfloor n_\beta \rfloor$ and note that
	$\beta_k =\min((n-k+1)/n-\epsilon,\beta)$ is equal to $\beta$
	iff $k\leq n_\beta+1$. Also, $\beta_k = \beta$ as soon as $k\leq \underline{n}_\beta+1$.
	Last, $\gamma_k = (n-k+1)/n-\epsilon$.
	
	{\bf Case 1} When $n_\beta<0$, then $\beta_k = \gamma_k$ for all $k\geq 1$.
	In this case, since $\gamma_k-\gamma_{k-1} = -1/n$, we deduce that
	\beqan
	I_\ell(1;\beta_1,\dots, \beta_\ell) &=& 
	I_\ell(1;\gamma_1,\dots, \gamma_\ell)\\
	&=& \frac{(1-\gamma_1)(1-\gamma_{1}+\ell/n)^{\ell-1}}{\ell!}\,,
	\eeqan
	and hence since $\gamma_1 = 1-\epsilon$,  and $\gamma_{\ell+1} = 1-\ell/n-\epsilon$,
	\beqan
	\Pr\Big(\sup_{u\in [\alpha,\beta]} U_n(u)-U(u) < \epsilon\Big)
	&=& 1 -\sum_{\ell=0}^{\overline{n}_{\alpha,\epsilon}-1} \binom{n}{\ell}\bigg( \Big(1-\ell/n-\epsilon\Big)^{n-\ell}-\alpha^{n-\ell}\bigg)\epsilon(\epsilon+\ell/n)^{\ell-1}\,.
	\eeqan
	
	{\bf Case 2} We now consider the general case when  $ \underline{n}_\beta\geq0$.
	For instance if $n_\beta\geq0$ but $ \underline{n}_\beta=0$ (that is, $0\leq n_\beta<1$), then, we deduce that $\beta_k = \gamma_k$ for all $k\geq 2$, while $\beta_1 = \beta$.		Hence, we deduce that
	\beqan
	I_\ell(1;\beta_1,\dots, \beta_\ell) =
	\begin{cases}
		\int_{\beta}^{1}dt_1 = (1-\beta) &\text{if } \ell = 1\\
		\int_{\beta}^{1}I_{\ell-1}(t_{1};\gamma_{2},\dots,\gamma_{\ell})dt_1 &\text{if } \ell >1\,.
	\end{cases}	
	\eeqan
	Likewise, when $ \underline{n}_\beta=1$, then we deduce that $\beta_k = \gamma_k$ for all $k\geq 3$, while $\beta_1 = \beta_2 = \beta$, and so
	\beqan
	I_\ell(1;\beta_1,\dots, \beta_\ell) =
	\begin{cases}
		\int_{\beta}^{1}\int_{\beta}^{t_1}\dots\int_{\beta}^{t_{\ell-1}}dt_\ell\dots dt_1 = \frac{(1-\beta)^\ell}{\ell!} &\text{if } \ell \leq 2\\
		\int_{\beta}^{1}\int_{\beta}^{t_1}I_{\ell-2}(t_{ \underline{n}_\beta+1};\gamma_{ \underline{n}_\beta+2},\dots,\gamma_{\ell})dt_{2}dt_1 &\text{if } \ell > 2\,.
	\end{cases}	
	\eeqan
	More generally, for a generic $ \underline{n}_\beta\geq 0$, we deduce (using the convention that $t_0=1$) that
	\beqan
	I_\ell(1;\beta_1,\dots, \beta_\ell) =
	\begin{cases}
		\int_{\beta}^{1}\int_{\beta}^{t_1}\dots\int_{\beta}^{t_{\ell-1}}dt_\ell\dots dt_1 = \frac{(1-\beta)^\ell}{\ell!} &\text{if } \ell \leq  \underline{n}_\beta+1\\
		\int_{\beta}^{1}\int_{\beta}^{t_1}\dots\int_{\beta}^{t_{ \underline{n}_\beta}}I_{\ell- \underline{n}_\beta-1}(t_{ \underline{n}_\beta+1};\gamma_{ \underline{n}_\beta+2},\dots,\gamma_{\ell})dt_{ \underline{n}_\beta+1}\dots dt_1 &\text{if } \ell >  \underline{n}_\beta+1\,.
	\end{cases}	
	\eeqan
	Further, since $\gamma_\ell-\gamma_{\ell-1}=-1/n$ for all $\ell$, and introducing $\ell_\beta = \ell- \underline{n}_\beta-1$ we also have
	(see \cite{smirnov1944approximate})
	\beqan
	I_{\ell_\beta}(t;\gamma_{ \underline{n}_\beta+2},\dots,\gamma_{\ell}) &=& \frac{(t-\gamma_{ \underline{n}_\beta+2})(t-\gamma_{ \underline{n}_\beta+2}+\ell_\beta/n)^{\ell_\beta-1}}{\ell_\beta!}\\
	&=&\frac{(t-\beta+C_{\ell_\beta})^{\ell_\beta}}{\ell_\beta!}
	- \frac{1}{n}\frac{(t-\beta+C_{\ell_\beta})^{\ell_\beta-1}}{(\ell_\beta-1)!}\,,
	\eeqan
	where in the second line, we also introduced
	$C_{\ell_\beta}= \beta-\gamma_{ \underline{n}_\beta+2}+\ell_\beta/n = (\beta+\epsilon)- (n-\ell)/n$. In particular, $C_{\ell_\beta}= (\ell-n_\beta)/n > 0$ for $ \ell >  \underline{n}_\beta+1$.  From this expression, we deduce that if $\ell >  \underline{n}_\beta+1$, then
	\beqan
	I_{\ell}(1;\beta_1,\dots, \beta_\ell)
	&=&
	\int_{\beta}^{1}\int_{\beta}^{t_1}\dots\int_{\beta}^{t_{ \underline{n}_\beta}} \frac{(t_{ \underline{n}_\beta+1}-\beta+C_{\ell_\beta})^{\ell_\beta}}{\ell_\beta!}dt_{ \underline{n}_\beta+1}\dots dt_1\\
	&&  - \frac{1}{n} 
	\int_{\beta}^{1}\int_{\beta}^{t_1}\dots\int_{\beta}^{t_{ \underline{n}_\beta}} \frac{(t_{ \underline{n}_\beta+1}-\beta+C_{\ell_\beta})^{\ell_\beta-1}}{(\ell_\beta-1)!}dt_{ \underline{n}_\beta+1}\dots dt_1\,.
	\eeqan 
	In order to compute both terms, we use the following inequality for given $k,\ell,C$,
	\beqan
	\lefteqn{
		\int_{\beta}^{1}\int_{\beta}^{t_1}\dots\int_{\beta}^{t_{k}} \frac{(t_{k+1}-\beta+C)^{\ell}}{\ell!}dt_{k+1}\dots dt_1}\\ &=&		
	\int_{\beta}^{1}\int_{\beta}^{t_1}\dots\int_{\beta}^{t_{k-1}} \frac{(t_{k}-\beta+C)^{\ell+1}}{(\ell+1)!}dt_{k}\dots dt_1- \frac{C^{\ell+1}}{(\ell+1)!}\underbrace{
		\int_{\beta}^{1}\!\int_{\beta}^{t_1}\!\dots\!\int_{\beta}^{t_{k-1}} dt_{k}\dots dt_1}_{B_k}\\
	&=&\int_{\beta}^{1}\int_{\beta}^{t_1}\dots\int_{\beta}^{t_{k-j-1}} \frac{(t_{k-j}\!-\!\beta\!+\!C)^{\ell+j+1}}{(\ell\!+\!j\!+\!1)!}dt_{k-j}\dots dt_1
	\!-\! \frac{C^{\ell+j+1}}{(\ell\!+\!j\!+\!1)!}B_{k-j}-\!\dots- \!\frac{C^{\ell+1}}{(\ell\!+\!1)!}B_k\\
	&=& \frac{(1-\beta+C)^{\ell+k+1}}{(\ell+k+1)!} - \sum_{j=0}^{k}\frac{C^{\ell+j+1}}{(\ell+j+1)!}\frac{(1-\beta)^{k-j}}{(k-j)!}\,.
	\eeqan
	Hence, we deduce that
	\beqan
	I_{\ell}(1;\beta_1,\dots, \beta_\ell)&=&
	\frac{(1-\beta+C_{\ell_\beta})^{\ell_\beta+ \underline{n}_\beta+1}}{(\ell_\beta+ \underline{n}_\beta+1)!} - \frac{1}{n}
	\frac{(1-\beta+C_{\ell_\beta})^{\ell_\beta+ \underline{n}_\beta}}{(\ell_\beta+ \underline{n}_\beta)!}\\
	&&- \sum_{j=0}^{ \underline{n}_\beta}\frac{C_{\ell_\beta}^{\ell_\beta+j+1}}{(\ell_\beta+j+1)!}\frac{(1-\beta)^{ \underline{n}_\beta-j}}{( \underline{n}_\beta-j)!}
	+\frac{1}{n} \sum_{j=0}^{ \underline{n}_\beta}\frac{C_{\ell_\beta}^{\ell_\beta+j}}{(\ell_\beta+j)!}\frac{(1-\beta)^{ \underline{n}_\beta-j}}{( \underline{n}_\beta-j)!}\\
	&=&\bigg[\frac{1-\beta+C_{\ell_\beta}}{\ell} - \frac{1}{n}\bigg]
	\frac{(1-\beta+C_{\ell_\beta})^{\ell-1}}{(\ell-1)!}\\
	&&- \sum_{j=0}^{ \underline{n}_\beta}\bigg[\frac{C_{\ell_\beta}}{\ell- \underline{n}_\beta+j}-\frac{1}{n}\bigg] \frac{C_{\ell_\beta}^{\ell- \underline{n}_\beta+j-1}}{(\ell- \underline{n}_\beta+j-1)!}\frac{(1-\beta)^{ \underline{n}_\beta-j}}{( \underline{n}_\beta-j)!}\\
	\eeqan
	After reorganizing the terms, and remarking that 
	$C_{\ell_\beta}= b-1+\epsilon + \ell/n = (\ell-n_\beta)/n$, we obtain that then if $\ell >  \underline{n}_\beta+1$, then
	\beqan
	I_{\ell}(1;\beta_1,\dots, \beta_\ell)&=&
	\frac{1}{n}\bigg[\frac{\ell\!+\!n\epsilon}{\ell} \!-\! 1\bigg]
	\frac{(\ell/n\!+\!\epsilon)^{\ell\!-\!1}}{(\ell\!-\!1)!}- \sum_{j=0}^{ \underline{n}_\beta}\frac{1}{n}\bigg[\frac{\ell\!-\!n_\beta}{\ell\!-\! \underline{n}_\beta\!+\!j}\!-\!1\bigg] \frac{((\ell\!-\!n_\beta)/n)^{\ell\!-\! \underline{n}_\beta\!+\!j\!-\!1}}{(\ell\!-\! \underline{n}_\beta\!+\!j\!-\!1)!}\frac{(1\!-\!\beta)^{ \underline{n}_\beta\!-\!j}}{( \underline{n}_\beta\!-\!j)!}\\
	&=&\epsilon
	\frac{(\ell/n\!+\!\epsilon)^{\ell\!-\!1}}{\ell!}+\frac{1}{n(\ell-1)!} \sum_{j=0}^{ \underline{n}_\beta}\frac{n_\beta\!-\!j}{\ell\!-\!j}\binom{\ell\!-\!1}{j}\bigg(\frac{\ell\!-\!n_\beta}{n}\bigg)^{\ell\!-\!j\!-\!1}(1\!-\!\beta)^{j}\\
	&=&\epsilon
	\frac{(\ell/n\!+\!\epsilon)^{\ell\!-\!1}}{\ell!}+\frac{1}{\ell!} \sum_{j=0}^{ \underline{n}_\beta}\frac{n_\beta\!-\!j}{n}\binom{\ell}{j}\bigg(\frac{\ell\!-\!n_\beta}{n}\bigg)^{\ell\!-\!j\!-\!1}(1\!-\!\beta)^{j}\,.
	\eeqan
	Combining all steps together,
	we deduce that if $ \underline{n}_\beta\geq0$ and $\lfloor n_\beta\rfloor+1\leq \overline{n}_{\alpha,\epsilon}-1$, then
	\beqan
	\lefteqn{
		\Pr\Big(\sup_{u\in [\alpha,\beta]} U_n(u)-U(u) < \epsilon\Big)=	1 - 
		\sum_{\ell=0}^{ \underline{n}_\beta+1}\binom{n}{\ell}\beta_{\ell+1}^{n-\ell}(1-\beta)^\ell} \\
	&&- \sum_{\ell= \underline{n}_\beta+2}^{\overline{n}_{\alpha,\epsilon}-1} 
	\binom{n}{\ell}\bigg(1\!-\! \frac{\ell}{n}\!-\!\epsilon\bigg)^{n-\ell}\epsilon
	\bigg(\!\frac{\ell}{n}\!+\!\epsilon\!\bigg)^{\ell\!-\!1}\\
	&&- \sum_{\ell= \underline{n}_\beta+2}^{\overline{n}_{\alpha,\epsilon}-1}
	\binom{n}{\ell}\bigg(1\!-\! \frac{\ell}{n}\!-\!\epsilon\bigg)^{n-\ell} \sum_{j=0}^{ \underline{n}_\beta}\bigg[\frac{n_\beta\!-\!j}{n}\bigg]\binom{\ell}{j}\bigg(\frac{\ell\!-\!n_\beta}{n}\bigg)^{\ell\!-\!j\!-\!1}(1\!-\!\beta)^{j}
	\eeqan
	where 	$\beta_k =\min((n-k+1)/n-\epsilon,b)$,	
	$\overline{n}_{\alpha,\epsilon}=\lceil n(1-\alpha-\epsilon)\rceil$
	$ \underline{n}_\beta=\lfloor n(1-\beta-\epsilon) \rfloor$.
	Introducing the term $m_\beta = \min\{\lfloor n_\beta\rfloor+1,\overline{n}_{\alpha,\epsilon}-1\}$, 
	we get more generally when $n_\beta>0$,
	\beqan
	\lefteqn{
		\Pr\Big(\sup_{u\in [\alpha,\beta]} U_n(u)-U(u) < \epsilon\Big)=	1 - 
		\sum_{\ell=0}^{m_\beta}\binom{n}{\ell}\Big(\!\min\Big\{1\!-\!\frac{\ell}{n}\!-\!\epsilon,\beta\Big\}\!\Big)^{n\!-\!\ell}(1-\beta)^\ell} \\
	&&-\!\!\sum_{\ell=m_\beta+1}^{\overline{n}_{\alpha,\epsilon}-1}\!\! 
	\binom{n}{\ell}\!\bigg(1\!-\! \frac{\ell}{n}\!-\!\epsilon\bigg)^{n\!-\!\ell}\!\bigg[\epsilon
	\bigg(\!\frac{\ell}{n}\!+\!\epsilon\!\bigg)^{\ell\!-\!1\!}\!\!+\! \sum_{j=0}^{m_\beta-1}\!\bigg[\frac{n_\beta\!-\!j}{n}\bigg]\binom{\ell}{j}\bigg(\frac{\ell\!-\!n_\beta}{n}\bigg)^{\ell\!-\!j\!-\!1\!}\!(1\!-\!\beta)^{j}\bigg]\,.
	\eeqan
	
\end{myproof}

\begin{myproof}{of Lemma~ \ref{thm:concRight}}
	
	We let $\tilde \beta=1-\beta, \tilde \alpha=1-\alpha$, $\tau_k=1-u_k$
	and consider for all $k$ the constant $\rho_k= 1-\epsilon-k/n$  (non-negative for all $k\leq n(1-\epsilon)$.
	We recall that 	$\tilde \alpha_k = \min(\rho_{k-1},\tilde \alpha)$.
	We let $n_{\beta,\epsilon}=n(\beta-\epsilon)$, $\overline{n}_{\beta,\epsilon}=\lceil n_{\beta,\epsilon} \rceil$ and remark that  $\rho_k\leq \tilde \beta$ iff $k\geq n_{\beta,\epsilon}$. Also, $\rho_k\leq \tilde \beta$ as soon as
	$k-1> \overline{n}_{\beta,\epsilon}$.		
	
	We now compute the quantity
	\beqan
	I_\ell(x;\tilde \alpha_1,\dots,\tilde \alpha_\ell) = \int_{\tilde \alpha_1}^{x}\int_{\tilde \alpha_2}^{t_1}\dots\int_{\tilde \alpha_\ell}^{t_{\ell-1}}dt_\ell\dots dt_1\,,
	\eeqan	
	where  $\tilde \alpha_k = \min(\rho_{k-1},\tilde \alpha)$.
	We further let $n_{\tilde \alpha}=n(\alpha-\epsilon)$, $\underline{n}_\alpha=\lfloor n_\alpha \rfloor$ and note that
	$\tilde \alpha_k =\min(1-\epsilon-(k-1)/n,\tilde \alpha)$ is equal to $\tilde \alpha$
	iff $k\leq n_{\tilde \alpha}+1$. 
	Also, $\tilde \alpha_k = \tilde \alpha$ as soon as $k\leq\underline{n}_{\tilde \alpha}+1$.

	{\bf Case 1} When $n_{\tilde \alpha}<0$, then $\tilde \alpha_k = \rho_{k-1}$ for all $k\geq 1$.
In this case, since $\rho_k-\rho_{k-1} = -1/n$, we deduce that
\beqan
I_\ell(1;\tilde \alpha_1,\dots, \tilde \alpha_\ell) &=& 
I_\ell(1;\rho_0,\dots, \rho_{\ell-1})\\
&=& \frac{(1-\rho_{0})(1-\rho_{0}+\ell/n)^{\ell-1}}{\ell!}\,,
\eeqan
and hence since $\rho_0 = 1-\epsilon$,  and $\rho_{\ell} = 1-\ell/n-\epsilon$, it comes
\beqan
\Pr\Big(\sup_{u\in [\alpha,\beta]} U(u)-U_n(u) \geq \epsilon\Big)
&=& \sum_{\ell=0}^{\overline{n}_{\beta,\epsilon}-1} \binom{n}{\ell} \Big(1-\frac{\ell}{n}-\epsilon\Big)^{n-\ell}\epsilon(\epsilon+\frac{\ell}{n})^{\ell-1}\,.
\eeqan

{\bf Case 2} We now consider the general case when  $ \underline{n}_{\tilde \alpha}\geq0$.
For instance if $n_{\tilde \alpha}\geq0$ but $ \underline{n}_{\tilde \alpha}=0$ (that is, $0\leq n_{\tilde \alpha}<1$), then, we deduce that $\tilde \alpha_k = \rho_{k-1}$ for all $k\geq 2$, while $\tilde \alpha_1 = \tilde \alpha$.		Hence, we deduce that
\beqan
I_\ell(1;\tilde \alpha_1,\dots, \tilde \alpha_\ell) =
\begin{cases}
	\int_{\tilde \alpha}^{1}dt_1 = (1-\tilde \alpha) &\text{if } \ell = 1\\
	\int_{\tilde \alpha}^{1}I_{\ell-1}(t_{1};\rho_{1},\dots,\rho_{\ell-1})dt_1 &\text{if } \ell >1\,.
\end{cases}	
\eeqan
Likewise, when $ \underline{n}_{\tilde \alpha}=1$, then we deduce that $\tilde \alpha_k = \rho_{k-1}$ for all $k\geq 3$, while $\tilde \alpha_1 = \tilde \alpha_2 = \tilde \alpha$, and so
\beqan
I_\ell(1;\tilde \alpha_1,\dots, \tilde \alpha_\ell) =
\begin{cases}
	\int_{\tilde \alpha}^{1}\int_{\tilde \alpha}^{t_1}\dots\int_{\tilde \alpha}^{t_{\ell-1}}dt_\ell\dots dt_1 = \frac{(1-\tilde \alpha)^\ell}{\ell!} &\text{if } \ell \leq 2\\
	\int_{\tilde \alpha}^{1}\int_{\tilde \alpha}^{t_1}I_{\ell-2}(t_{ \underline{n}_{\tilde \alpha}+1};\rho_{ \underline{n}_{\tilde \alpha}+1},\dots,\rho_{\ell-1})dt_{2}dt_1 &\text{if } \ell > 2\,.
\end{cases}	
\eeqan
More generally, for a generic $ \underline{n}_{\tilde \alpha}\geq 0$, we deduce (using the convention that $t_0=1$) that
\beqan
I_\ell(1;\tilde \alpha_1,\dots, \tilde \alpha_\ell) =
\begin{cases}
	\int_{\tilde \alpha}^{1}\int_{\tilde \alpha}^{t_1}\dots\int_{\tilde \alpha}^{t_{\ell-1}}dt_\ell\dots dt_1 = \frac{(1-\tilde \alpha)^\ell}{\ell!} &\text{if } \ell \leq  \underline{n}_{\tilde \alpha}+1\\
	\int_{\tilde \alpha}^{1}\int_{\tilde \alpha}^{t_1}\dots\int_{\tilde \alpha}^{t_{ \underline{n}_{\tilde \alpha}}}I_{\ell- \underline{n}_{\tilde \alpha}-1}(t_{ \underline{n}_{\tilde \alpha}+1};\rho_{ \underline{n}_{\tilde \alpha}+1},\dots,\rho_{\ell-1})dt_{ \underline{n}_{\tilde \alpha}+1}\dots dt_1 &\text{if } \ell >  \underline{n}_{\tilde \alpha}+1\,.
\end{cases}	
\eeqan
Further, since $\rho_\ell-\rho_{\ell-1}=-1/n$ for all $\ell$, and introducing $\ell_{\tilde \alpha} = \ell- \underline{n}_{\tilde \alpha}-1$ we also have
(see \cite{smirnov1944approximate})
\beqan
I_{\ell_{\tilde \alpha}}(t;\rho_{ \underline{n}_{\tilde \alpha}+1},\dots,\rho_{\ell-1}) &=& \frac{(t-\rho_{ \underline{n}_{\tilde \alpha}+1})(t-\rho_{ \underline{n}_{\tilde \alpha}+1}+\ell_{\tilde \alpha}/n)^{\ell_{\tilde \alpha}-1}}{\ell_{\tilde \alpha}!}\\
&=&\frac{(t-{\tilde \alpha}+C_{\ell_{\tilde \alpha}})^{\ell_{\tilde \alpha}}}{\ell_{\tilde \alpha}!}
- \frac{1}{n}\frac{(t-{\tilde \alpha}+C_{\ell_{\tilde \alpha}})^{\ell_{\tilde \alpha}-1}}{(\ell_{\tilde \alpha}-1)!}\,,
\eeqan
where in the second line, we also introduced
$C_{\ell_{\tilde \alpha}}= {\tilde \alpha}-\rho_{ \underline{n}_{\tilde \alpha}+1}+\ell_{\tilde \alpha}/n = (1-\alpha+\epsilon)- (n-\ell)/n$.

In particular, it holds that $C_{\ell_{\tilde \alpha}}= (\ell-n_{\tilde \alpha})/n > 0$ for $ \ell >  \underline{n}_{\tilde \alpha}+1$.  From this expression, we deduce that if $\ell >  \underline{n}_{\tilde \alpha}+1$, then
\beqan
I_{\ell}(1;{\tilde \alpha}_1,\dots, {\tilde \alpha}_\ell)
&=&
\int_{{\tilde \alpha}}^{1}\int_{{\tilde \alpha}}^{t_1}\dots\int_{{\tilde \alpha}}^{t_{ \underline{n}_{\tilde \alpha}}} \frac{(t_{ \underline{n}_{\tilde \alpha}+1}-{\tilde \alpha}+C_{\ell_{\tilde \alpha}})^{\ell_{\tilde \alpha}}}{\ell_{\tilde \alpha}!}dt_{ \underline{n}_{\tilde \alpha}+1}\dots dt_1\\
&&  - \frac{1}{n} 
\int_{{\tilde \alpha}}^{1}\int_{{\tilde \alpha}}^{t_1}\dots\int_{{\tilde \alpha}}^{t_{ \underline{n}_{\tilde \alpha}}} \frac{(t_{ \underline{n}_{\tilde \alpha}+1}-{\tilde \alpha}+C_{\ell_{\tilde \alpha}})^{\ell_{\tilde \alpha}-1}}{(\ell_{\tilde \alpha}-1)!}dt_{ \underline{n}_{\tilde \alpha}+1}\dots dt_1\,.
\eeqan 
Hence, we deduce that
\beqan
I_{\ell}(1;{\tilde \alpha}_1,\dots, {\tilde \alpha}_\ell)&=&
\frac{(1-{\tilde \alpha}+C_{\ell_{\tilde \alpha}})^{\ell_{\tilde \alpha}+ \underline{n}_{\tilde \alpha}+1}}{(\ell_{\tilde \alpha}+ \underline{n}_{\tilde \alpha}+1)!} - \frac{1}{n}
\frac{(1-{\tilde \alpha}+C_{\ell_{\tilde \alpha}})^{\ell_{\tilde \alpha}+ \underline{n}_{\tilde \alpha}}}{(\ell_{\tilde \alpha}+ \underline{n}_{\tilde \alpha})!}\\
&&- \sum_{j=0}^{ \underline{n}_{\tilde \alpha}}\frac{C_{\ell_{\tilde \alpha}}^{\ell_{\tilde \alpha}+j+1}}{(\ell_{\tilde \alpha}+j+1)!}\frac{(1-{\tilde \alpha})^{ \underline{n}_{\tilde \alpha}-j}}{( \underline{n}_{\tilde \alpha}-j)!}
+\frac{1}{n} \sum_{j=0}^{ \underline{n}_{\tilde \alpha}}\frac{C_{\ell_{\tilde \alpha}}^{\ell_{\tilde \alpha}+j}}{(\ell_{\tilde \alpha}+j)!}\frac{(1-{\tilde \alpha})^{ \underline{n}_{\tilde \alpha}-j}}{( \underline{n}_{\tilde \alpha}-j)!}\\
&=&\bigg[\frac{1-{\tilde \alpha}+C_{\ell_{\tilde \alpha}}}{\ell} - \frac{1}{n}\bigg]
\frac{(1-{\tilde \alpha}+C_{\ell_{\tilde \alpha}})^{\ell-1}}{(\ell-1)!}\\
&&- \sum_{j=0}^{ \underline{n}_{\tilde \alpha}}\bigg[\frac{C_{\ell_{\tilde \alpha}}}{\ell- \underline{n}_{\tilde \alpha}+j}-\frac{1}{n}\bigg] \frac{C_{\ell_{\tilde \alpha}}^{\ell- \underline{n}_{\tilde \alpha}+j-1}}{(\ell- \underline{n}_{\tilde \alpha}+j-1)!}\frac{(1-{\tilde \alpha})^{ \underline{n}_{\tilde \alpha}-j}}{( \underline{n}_{\tilde \alpha}-j)!}\\
\eeqan
where we used that $\ell_{\tilde \alpha}+ \underline{n}_{\tilde \alpha}+1 = \ell$.
After reorganizing the terms, and remarking that 
$C_{\ell_{\tilde \alpha}}= \epsilon-\alpha + \ell/n = (\ell-n_{\tilde \alpha})/n$, we obtain that then if $\ell >  \underline{n}_{\tilde \alpha}+1$, then
\beqan
I_{\ell}(1;{\tilde \alpha}_1,\dots, {\tilde \alpha}_\ell)\!\!&\!\!=\!&\!
\!\frac{1}{n}\bigg[\frac{\ell\!+\!n\epsilon}{\ell} \!-\! 1\bigg]
\frac{(\ell/n\!+\!\epsilon)^{\ell\!-\!1}}{(\ell\!-\!1)!}\!- \!
\sum_{j=0}^{ \underline{n}_{\tilde \alpha}}\frac{1}{n}\bigg[\frac{\ell\!-\!n_{\tilde \alpha}}{\ell\!-\! \underline{n}_{\tilde \alpha}\!+\!j}\!-\!1\!\bigg] \!\frac{((\ell\!-\!n_{\tilde \alpha})/n)^{\ell\!-\! \underline{n}_{\tilde \alpha}\!+\!j\!-\!1}}{(\ell\!-\! \underline{n}_{\tilde \alpha}\!+\!j\!-\!1)!}\frac{(1\!-\!{\tilde \alpha})^{ \underline{n}_{\tilde \alpha}\!-\!j}}{( \underline{n}_{\tilde \alpha}\!-\!j)!}\\
&=&\epsilon
\frac{(\ell/n\!+\!\epsilon)^{\ell\!-\!1}}{\ell!}+	
\frac{1}{n}\sum_{j=0}^{ \underline{n}_{\tilde \alpha}}\bigg(n_{\tilde \alpha}\!-\!j\bigg)\bigg(\frac{\ell\!-\!n_{\tilde \alpha}}{n}\bigg)^{\ell\!-\!j\!-\!1}(1\!-\!{\tilde \alpha})^{ j} \frac{1}{\ell!}\binom{\ell}{j}
\\	
&=&\epsilon
\frac{(\ell/n\!+\!\epsilon)^{\ell\!-\!1}}{\ell!}+\frac{1}{\ell!} \sum_{j=0}^{ \underline{n}_{\tilde \alpha}}\frac{n_{\tilde \alpha}\!-\!j}{n}\binom{\ell}{j}\bigg(\frac{\ell\!-\!n_{\tilde \alpha}}{n}\bigg)^{\ell\!-\!j\!-\!1}(1\!-\!{\tilde \alpha})^{j}\,.
\eeqan

Combining all steps together,
we deduce that if $\underline{n}_{\tilde \alpha}\geq0$ and $\lfloor n_{\tilde \alpha}\rfloor+1\leq \overline{n}_{\beta,\epsilon}-1$, then
\beqan
\lefteqn{
	\Pr\Big(\sup_{u\in [\alpha,\beta]} U_n(u)-U(u) < \epsilon\Big)=	1 - 
	\sum_{\ell=0}^{ \underline{n}_{\tilde \alpha}+1}\binom{n}{\ell}\tilde \alpha_{\ell+1}^{n-\ell}(1-{\tilde \alpha})^\ell} \\
&&- \sum_{\ell= \underline{n}_{\tilde \alpha}+2}^{\overline{n}_{\beta,\epsilon}-1} 
\binom{n}{\ell}\bigg(1\!-\! \frac{\ell}{n}\!-\!\epsilon\bigg)^{n-\ell}\epsilon
\bigg(\!\frac{\ell}{n}\!+\!\epsilon\!\bigg)^{\ell\!-\!1}\\
&&- \sum_{\ell= \underline{n}_{\tilde \alpha}+2}^{\overline{n}_{\beta,\epsilon}-1}
\binom{n}{\ell}\bigg(1\!-\! \frac{\ell}{n}\!-\!\epsilon\bigg)^{n-\ell} \sum_{j=0}^{ \underline{n}_{\tilde \alpha}}\bigg[\frac{n_{\tilde \alpha}\!-\!j}{n}\bigg]\binom{\ell}{j}\bigg(\frac{\ell\!-\!n_{\tilde \alpha}}{n}\bigg)^{\ell\!-\!j\!-\!1}(1\!-\!{\tilde \alpha})^{j}
\eeqan
where 	${\tilde \alpha}_k =\min((n-(k-1))/n-\epsilon,\tilde \alpha)$,	
$\overline{n}_{\beta,\epsilon}=\lceil n(\beta-\epsilon)\rceil$
$ \underline{n}_{\tilde \alpha}=\lfloor n(\alpha-\epsilon) \rfloor$.
Introducing the term $m_\beta = \min\{\lfloor n_{\tilde \alpha}\rfloor+1,\overline{n}_{\beta,\epsilon}-1\}$, 
we get more generally when $n_{\tilde \alpha}>0$,
\beqan
\lefteqn{
	\Pr\Big(\sup_{u\in [\alpha,\beta]} U_n(u)-U(u) < \epsilon\Big)=	1 - 
	\sum_{\ell=0}^{m_\beta}\binom{n}{\ell}\Big(\!\min\Big\{1\!-\!\frac{\ell}{n}\!-\!\epsilon,\tilde \alpha\Big\}\!\Big)^{n\!-\!\ell}(1-{\tilde \alpha})^\ell} \\
&&-\!\!\sum_{\ell=m_\beta+1}^{\overline{n}_{\beta,\epsilon}-1}\!\! 
\binom{n}{\ell}\!\bigg(1\!-\! \frac{\ell}{n}\!-\!\epsilon\bigg)^{n\!-\!\ell}\bigg[\epsilon
\bigg(\!\frac{\ell}{n}\!+\!\epsilon\!\bigg)^{\ell\!-\!1\!}\!\!+\! \sum_{j=0}^{m_\beta-1}\!\bigg[\frac{n_{\tilde \alpha}\!-\!j}{n}\bigg]\binom{\ell}{j}\bigg(\frac{\ell\!-\!n_{\tilde \alpha}}{n}\bigg)^{\ell\!-\!j\!-\!1\!}\!(1\!-\!{\tilde \alpha})^{j}\bigg]\,.
\eeqan

\end{myproof}

\section{Monte Carlo simulations of the confidence bounds}\label{app:MCMC}

\vspace{-0.5cm}
\begin{figure}[H]
	\centering
	\includegraphics[width=0.43\linewidth, trim = {0.9cm 0.8cm 1cm 0.8cm}, clip]{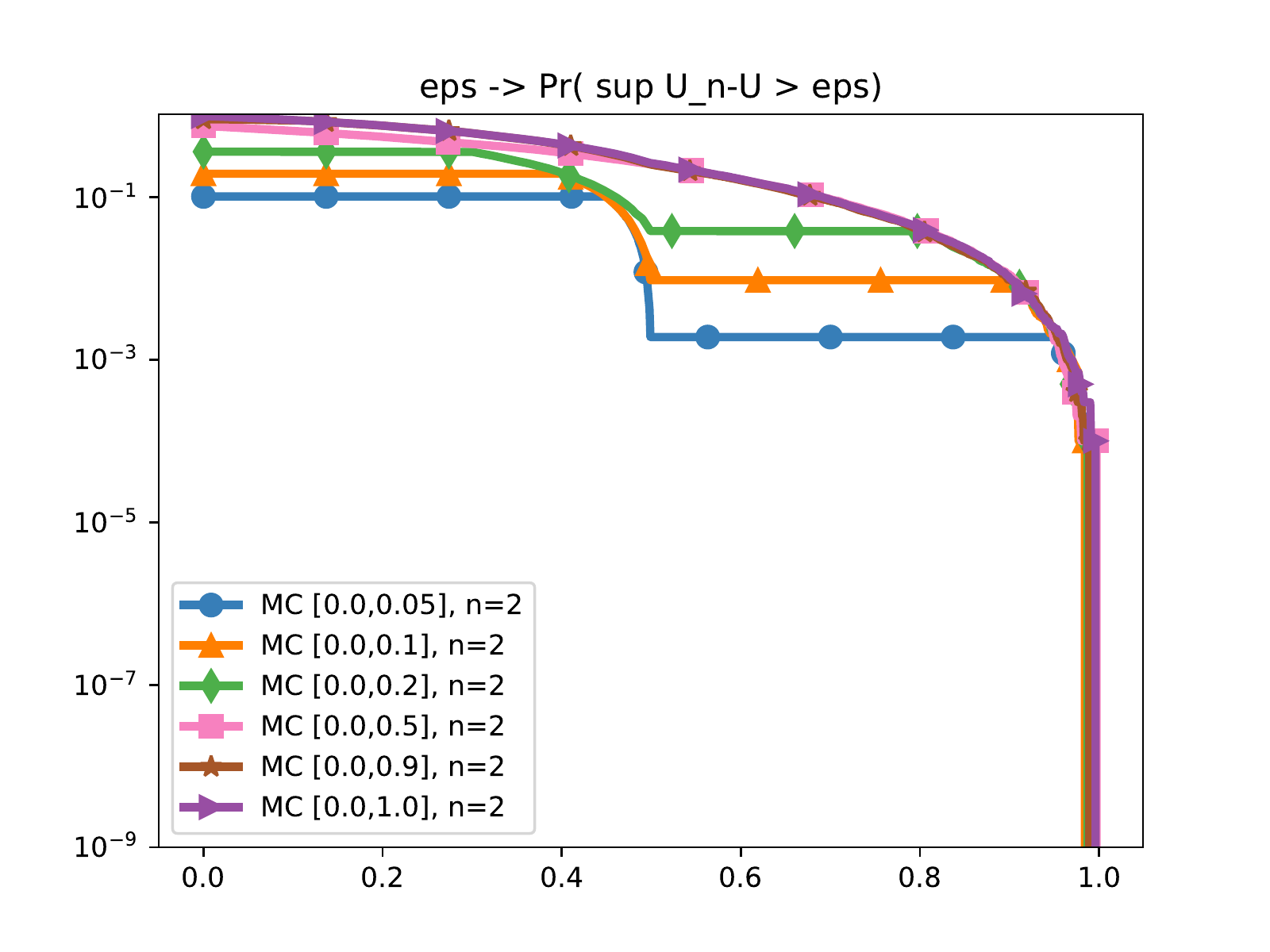}
	\includegraphics[width=0.43\linewidth, trim = {0.9cm 0.8cm 1cm 0.8cm}, clip]{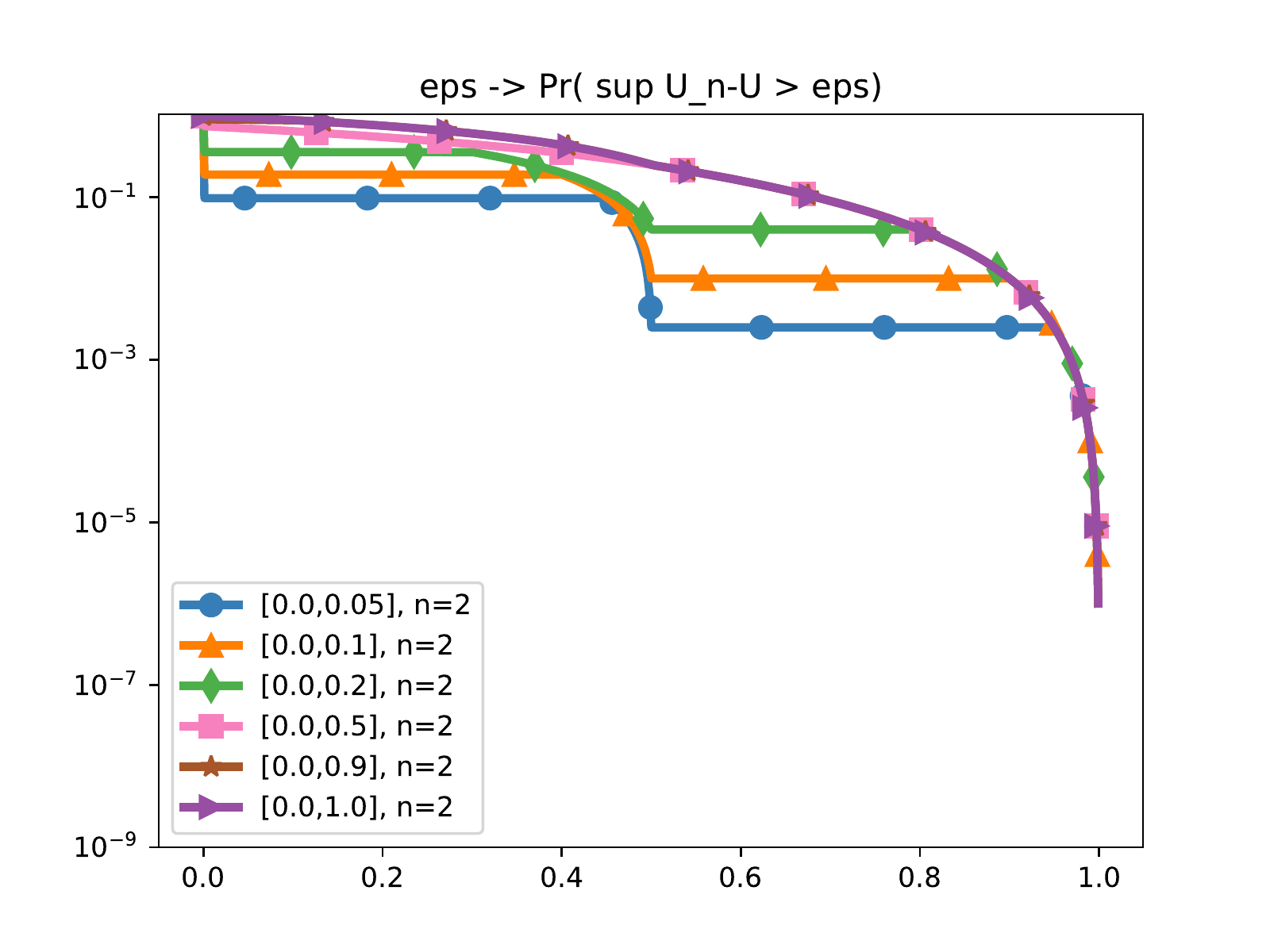}\\
	\includegraphics[width=0.43\linewidth, trim = {0.9cm 0.8cm 1cm 1.4cm}, clip]{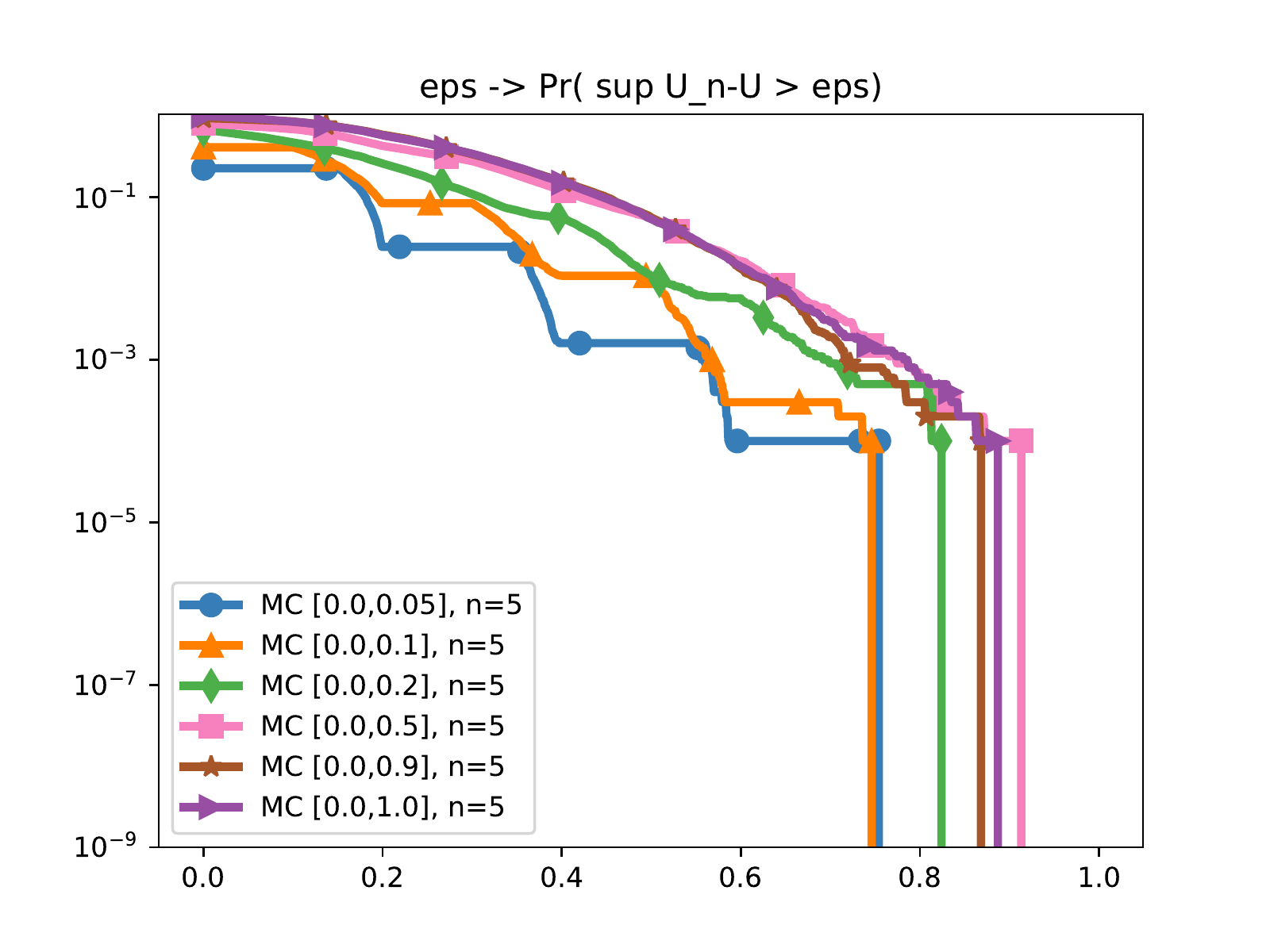}
	\includegraphics[width=0.43\linewidth, trim = {0.9cm 0.8cm 1cm 1.4cm}, clip]{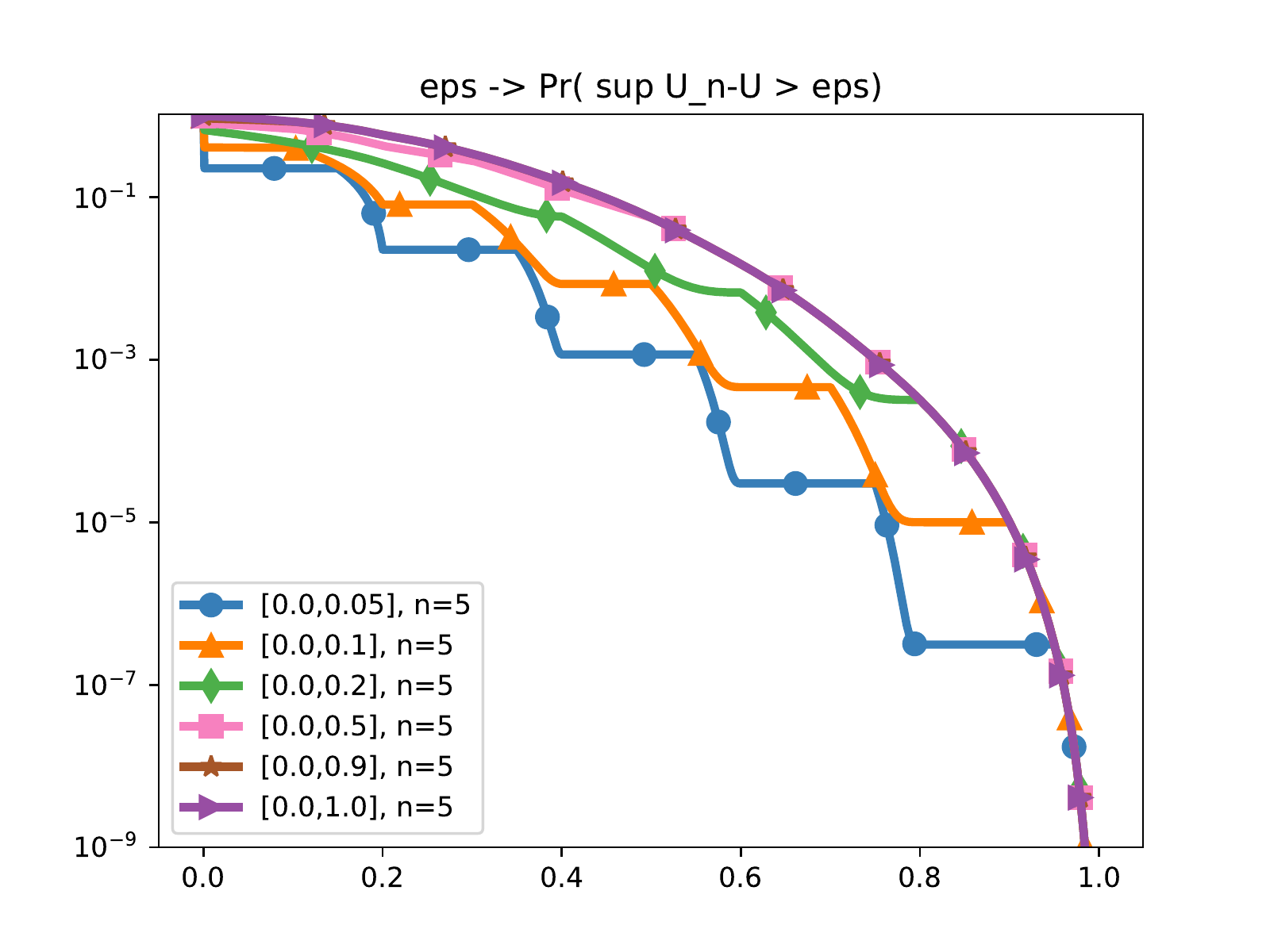}\\
	\includegraphics[width=0.43\linewidth, trim = {0.9cm 0.8cm 1cm 1.4cm}, clip]{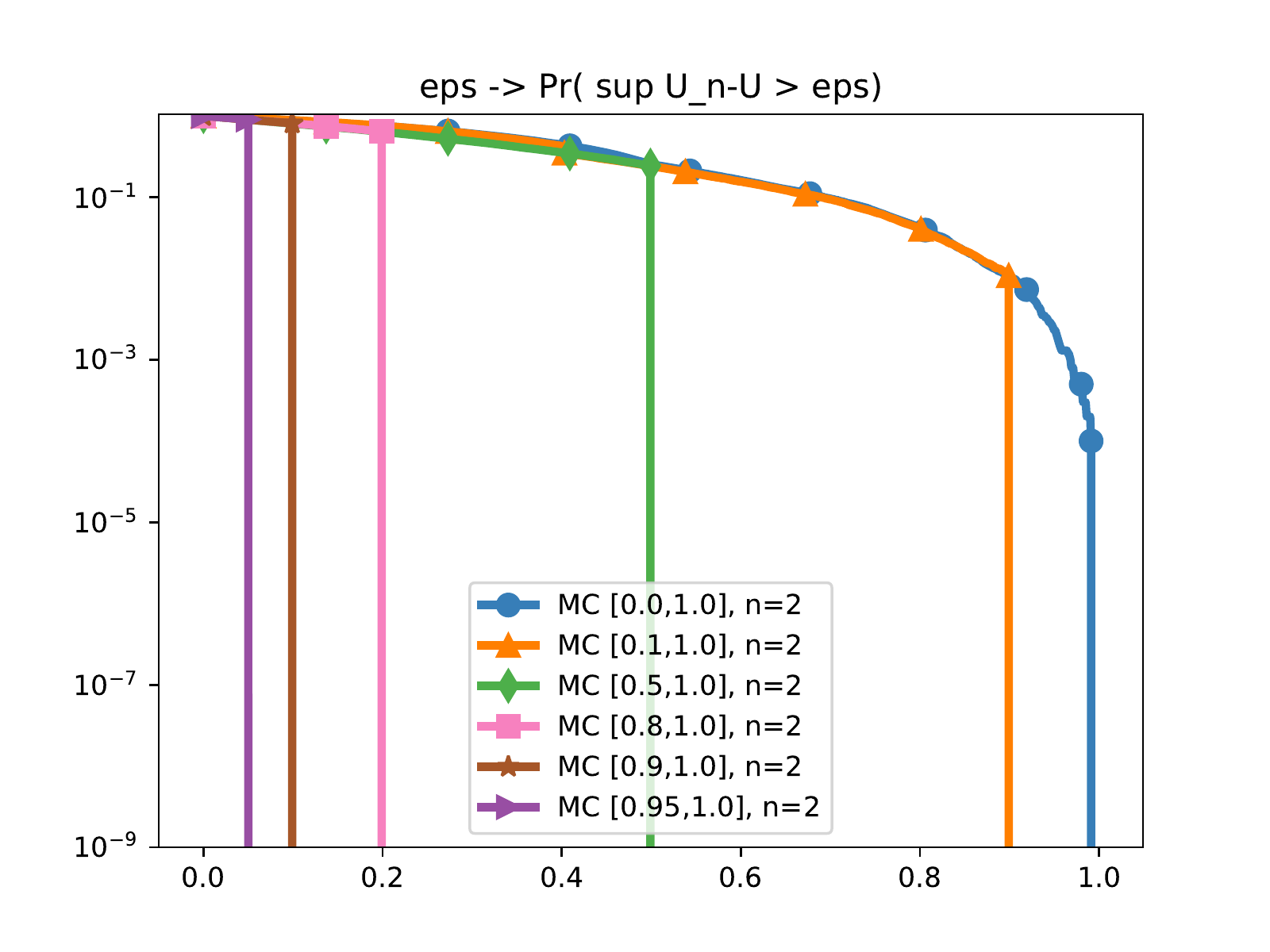}
	\includegraphics[width=0.43\linewidth, trim = {0.9cm 0.8cm 1cm 1.4cm}, clip]{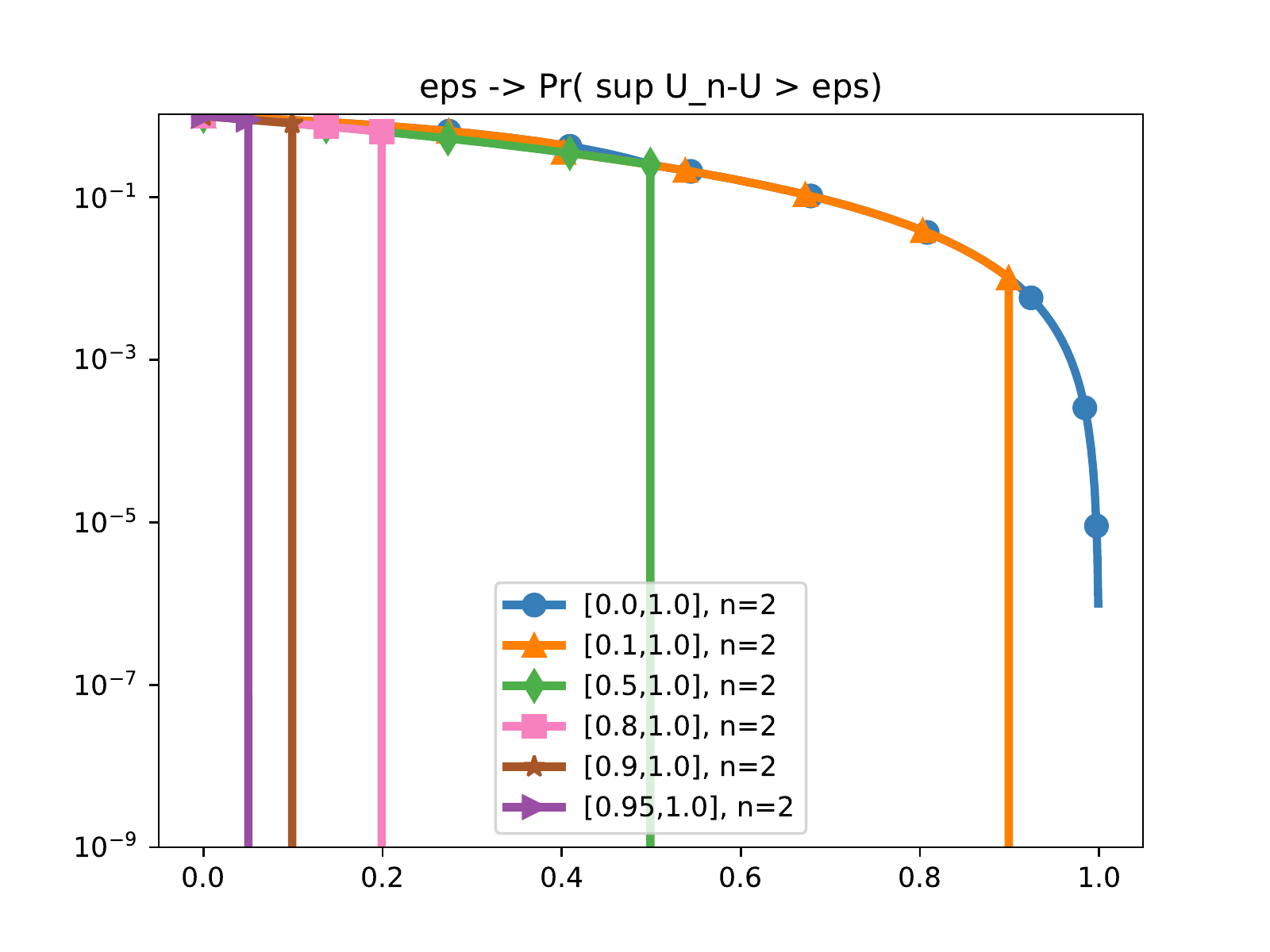}\\
	\includegraphics[width=0.43\linewidth, trim = {0.9cm 0.8cm 1cm 1.4cm}, clip]{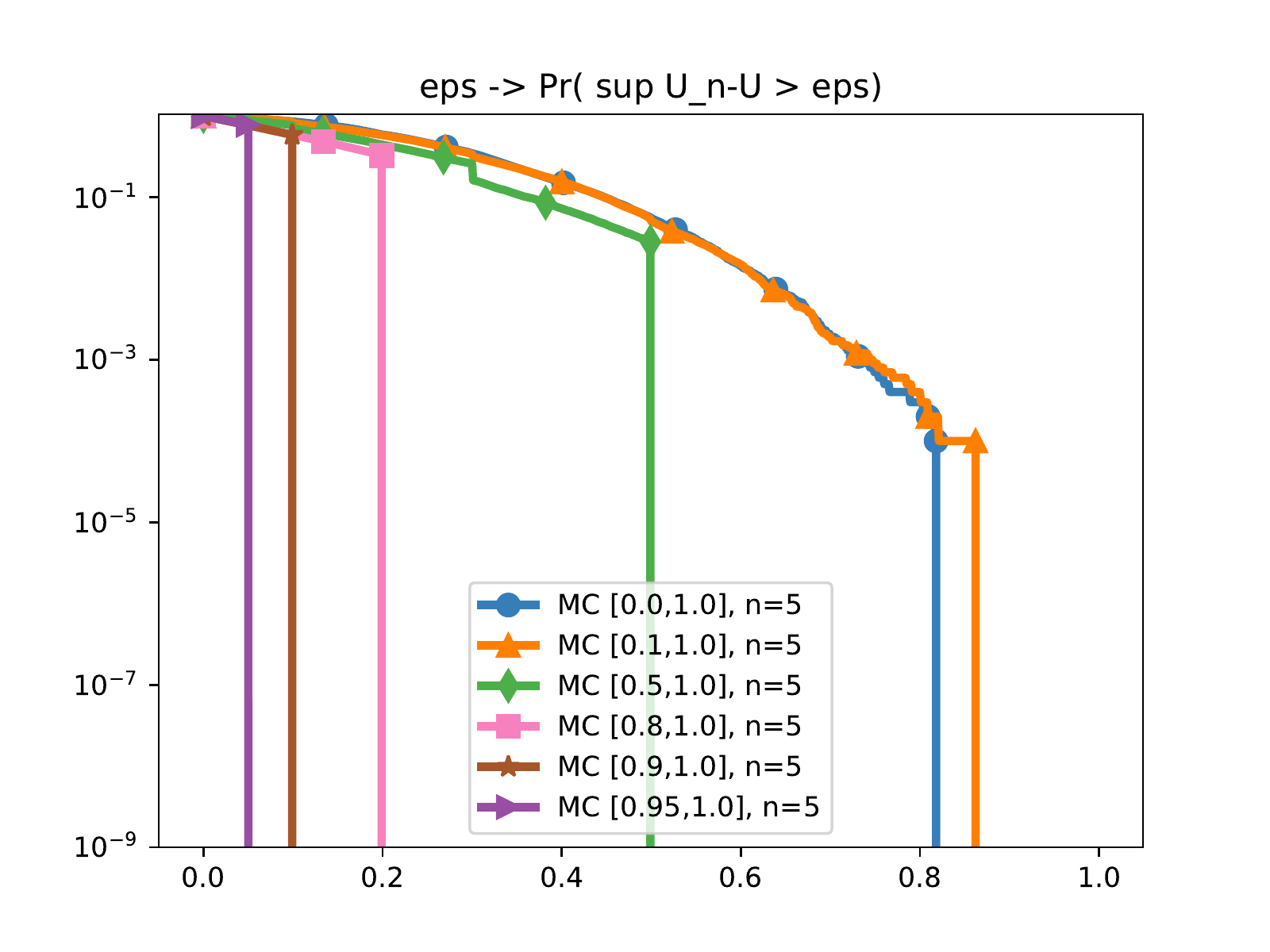}
	\includegraphics[width=0.43\linewidth, trim = {0.9cm 0.8cm 1cm 1.4cm}, clip]{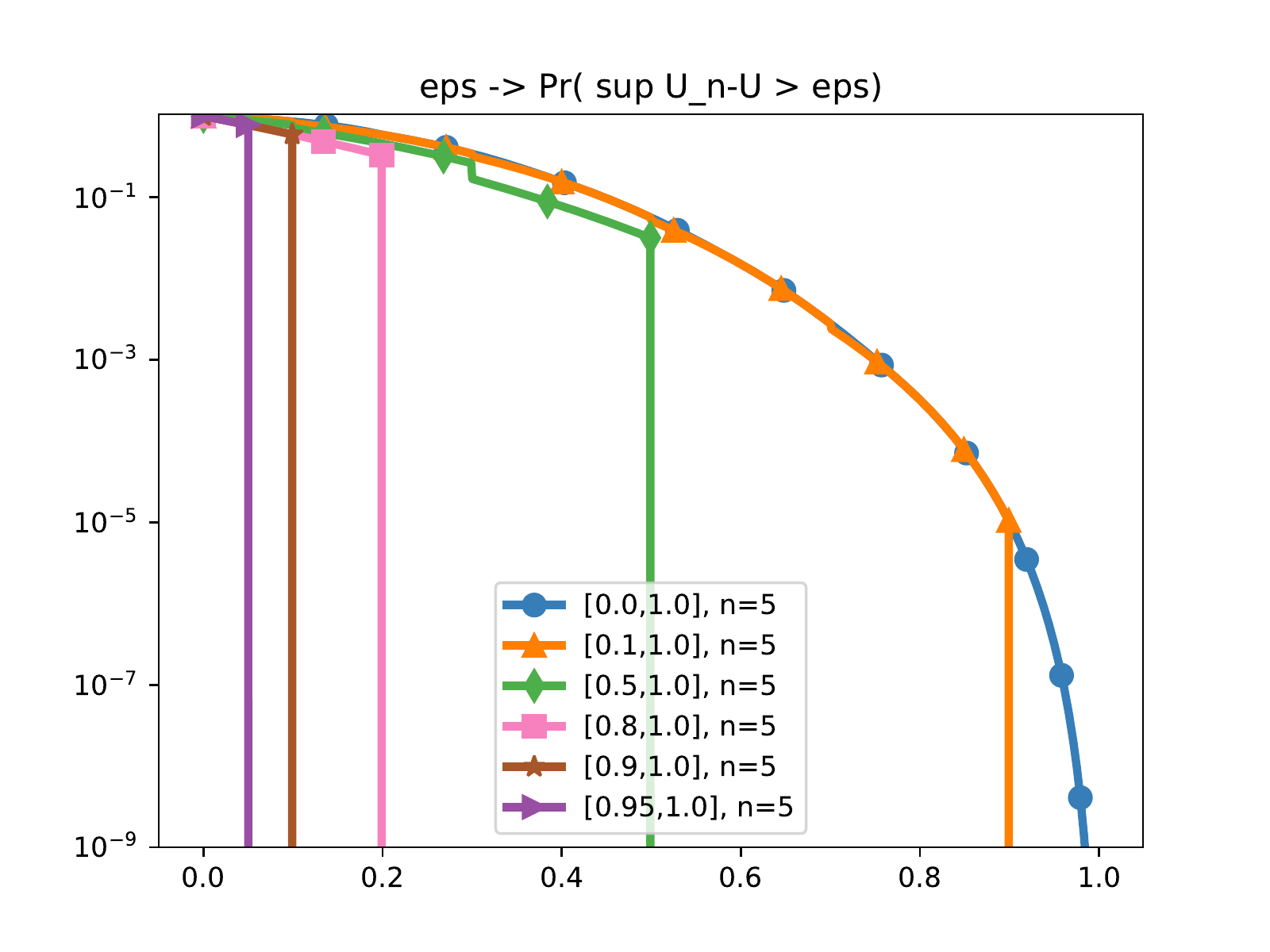}
	\vspace{-0.5cm}
	\caption{MCMC (left) versus Exact (Right) plot of $\epsilon\to\delta_{[\uu,\ou]}(n,\epsilon)$ for various values of $n$ and 
		interval $[\uu,\ou]$ build from $M=10^4$ replicate.}
	\label{fig:MCMC1}
\end{figure}

\begin{figure}[H]
	\centering
	\includegraphics[width=0.43\linewidth, trim = {0.9cm 0.8cm 1cm 0.8cm}, clip]{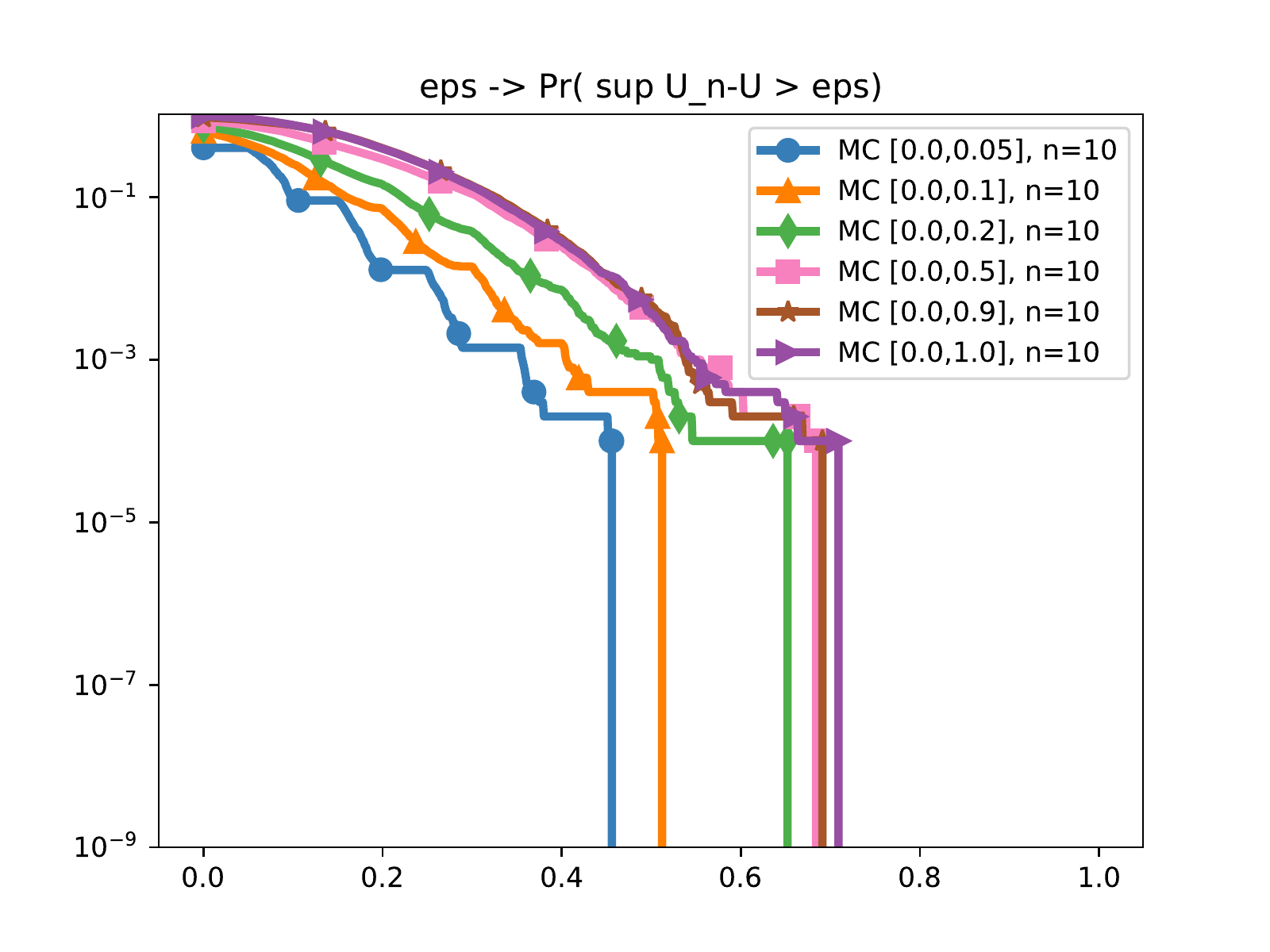}
	\includegraphics[width=0.43\linewidth, trim = {0.9cm 0.8cm 1cm 0.8cm}, clip]{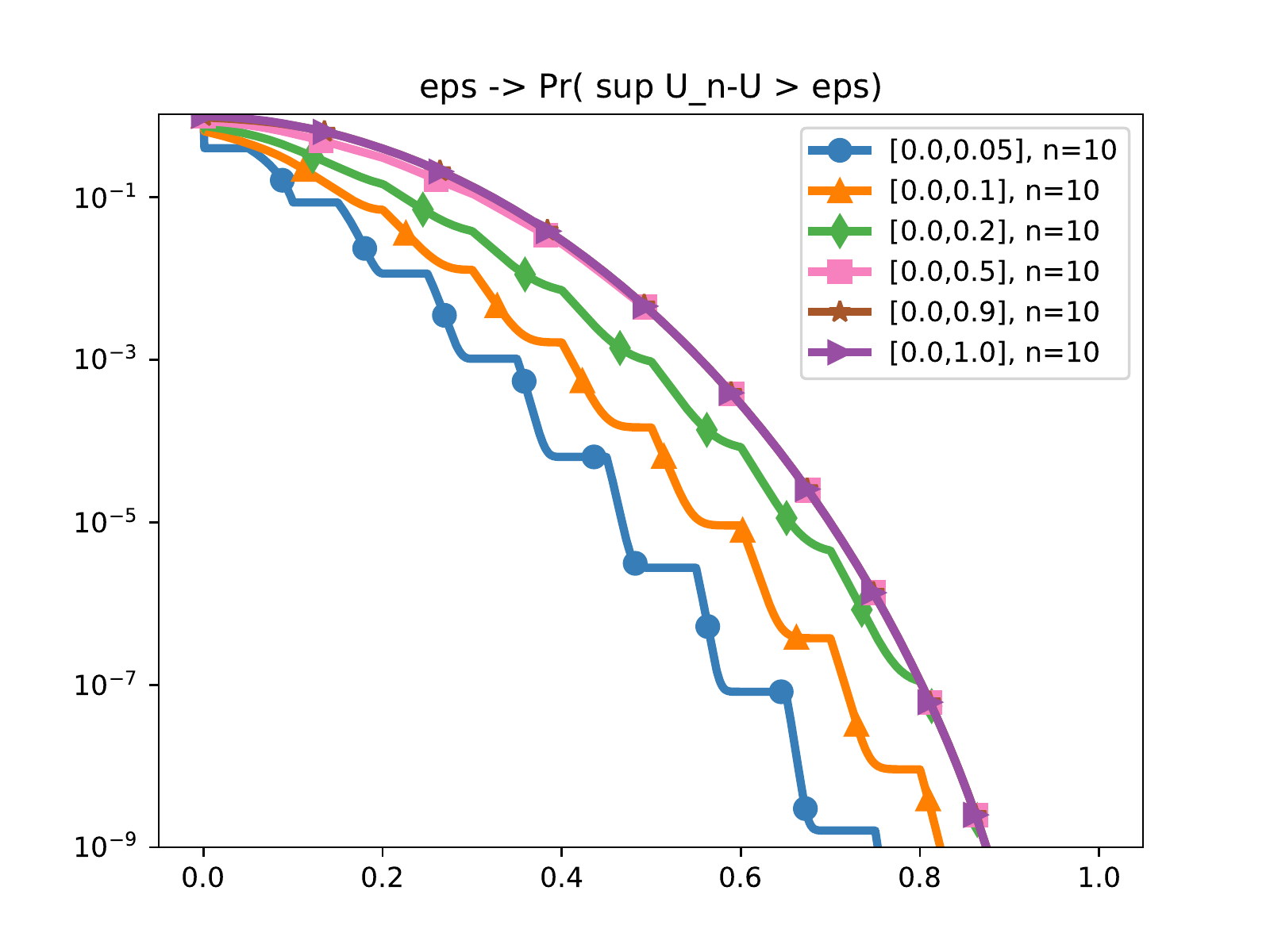}\\
	\includegraphics[width=0.43\linewidth, trim = {0.9cm 0.8cm 1cm 0.8cm}, clip]{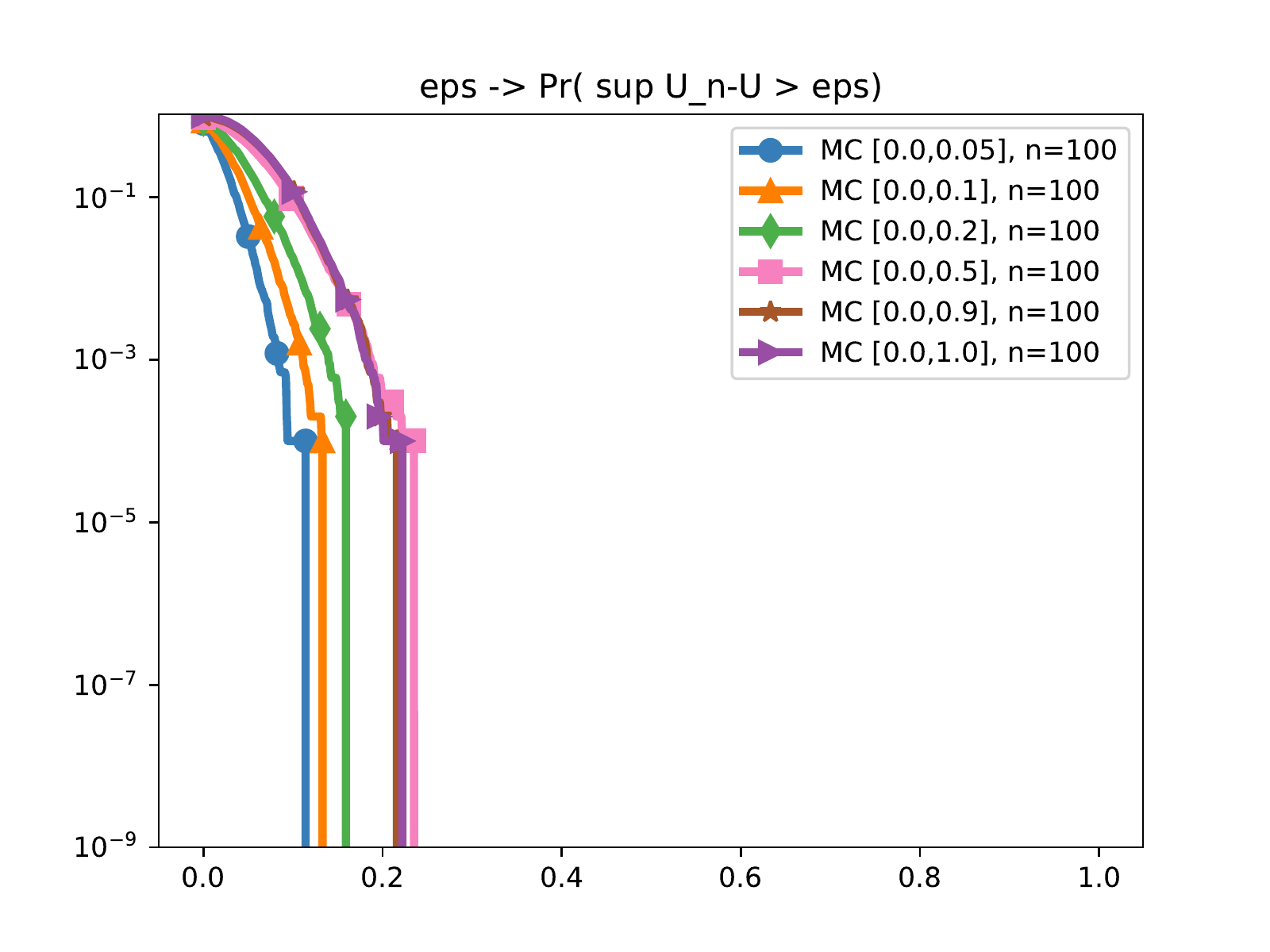}
	\includegraphics[width=0.43\linewidth, trim = {0.9cm 0.8cm 1cm 0.8cm}, clip]{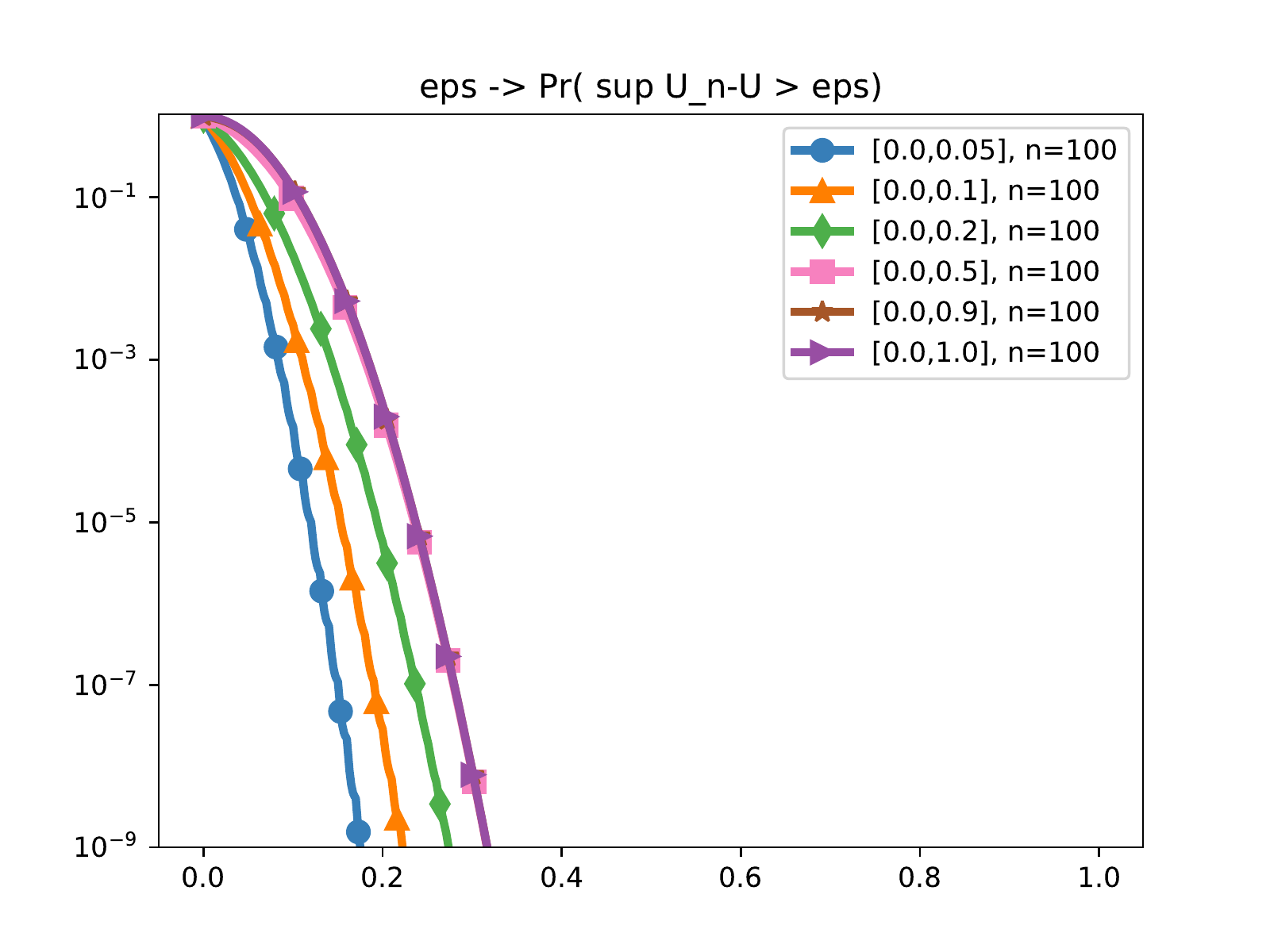}\\
	\includegraphics[width=0.43\linewidth, trim = {0.9cm 0.8cm 1cm 0.8cm}, clip]{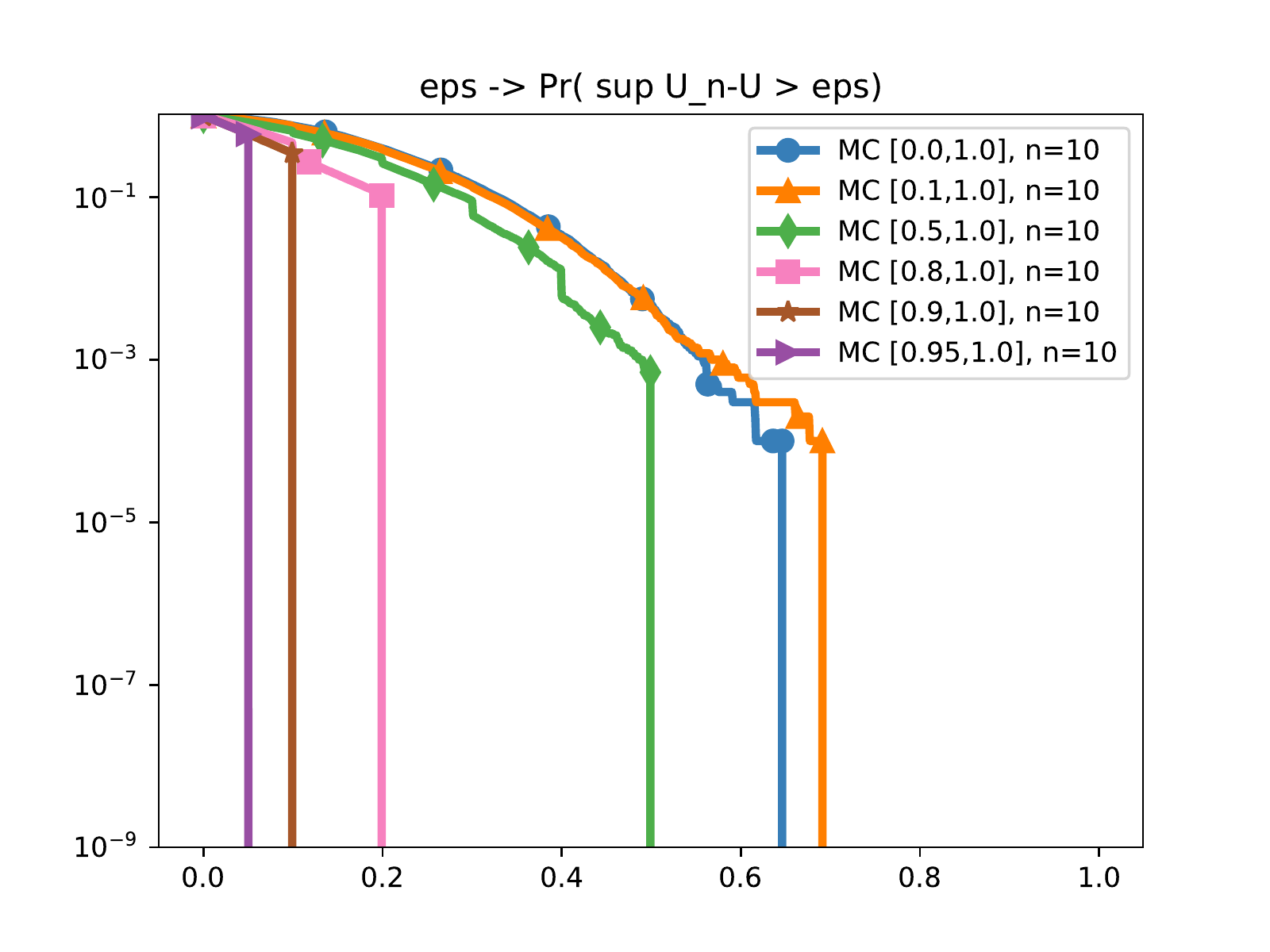}
	\includegraphics[width=0.43\linewidth, trim = {0.9cm 0.8cm 1cm 0.8cm}, clip]{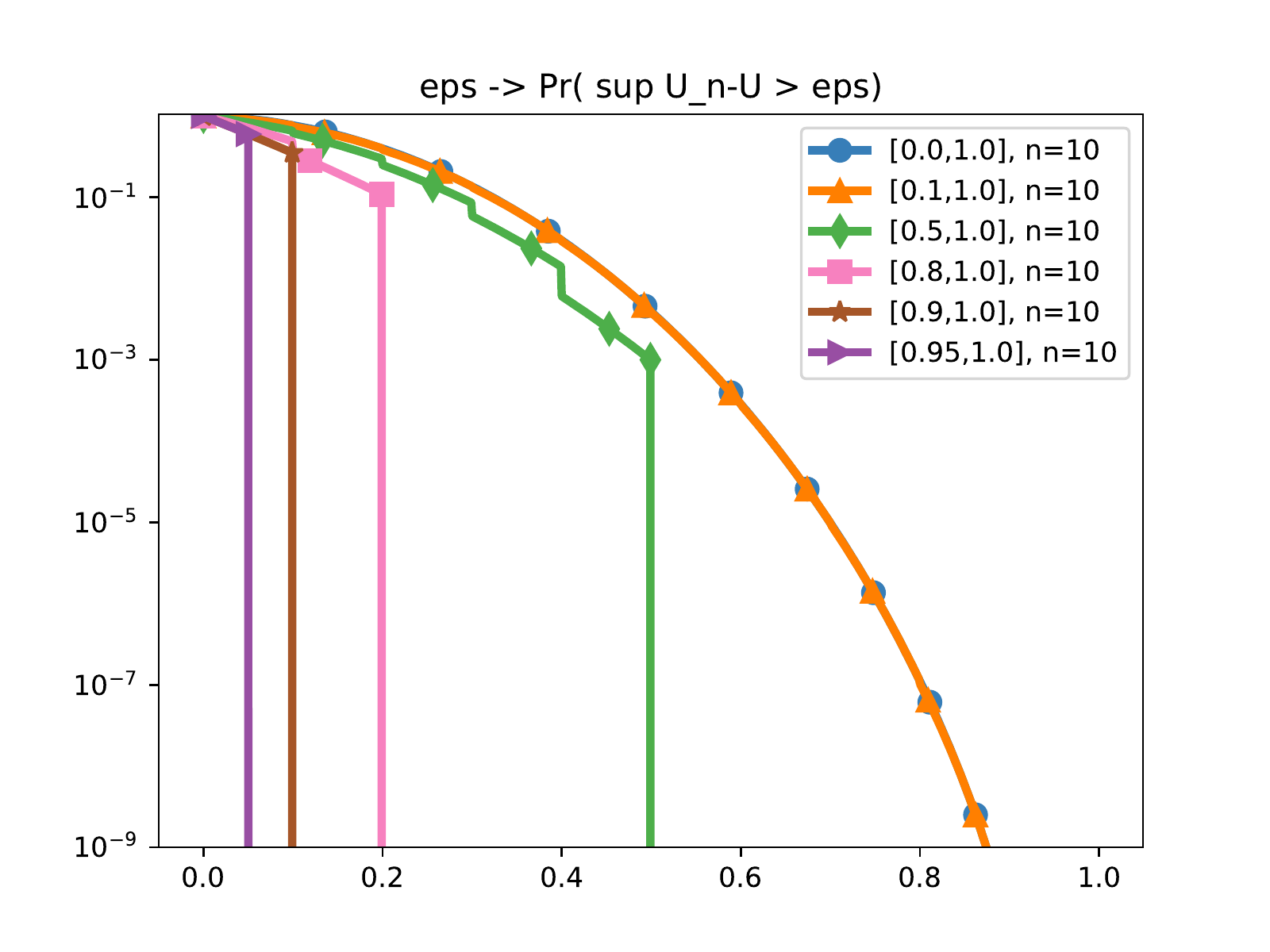}\\
	\includegraphics[width=0.43\linewidth, trim = {0.9cm 0.8cm 1cm 0.8cm}, clip]{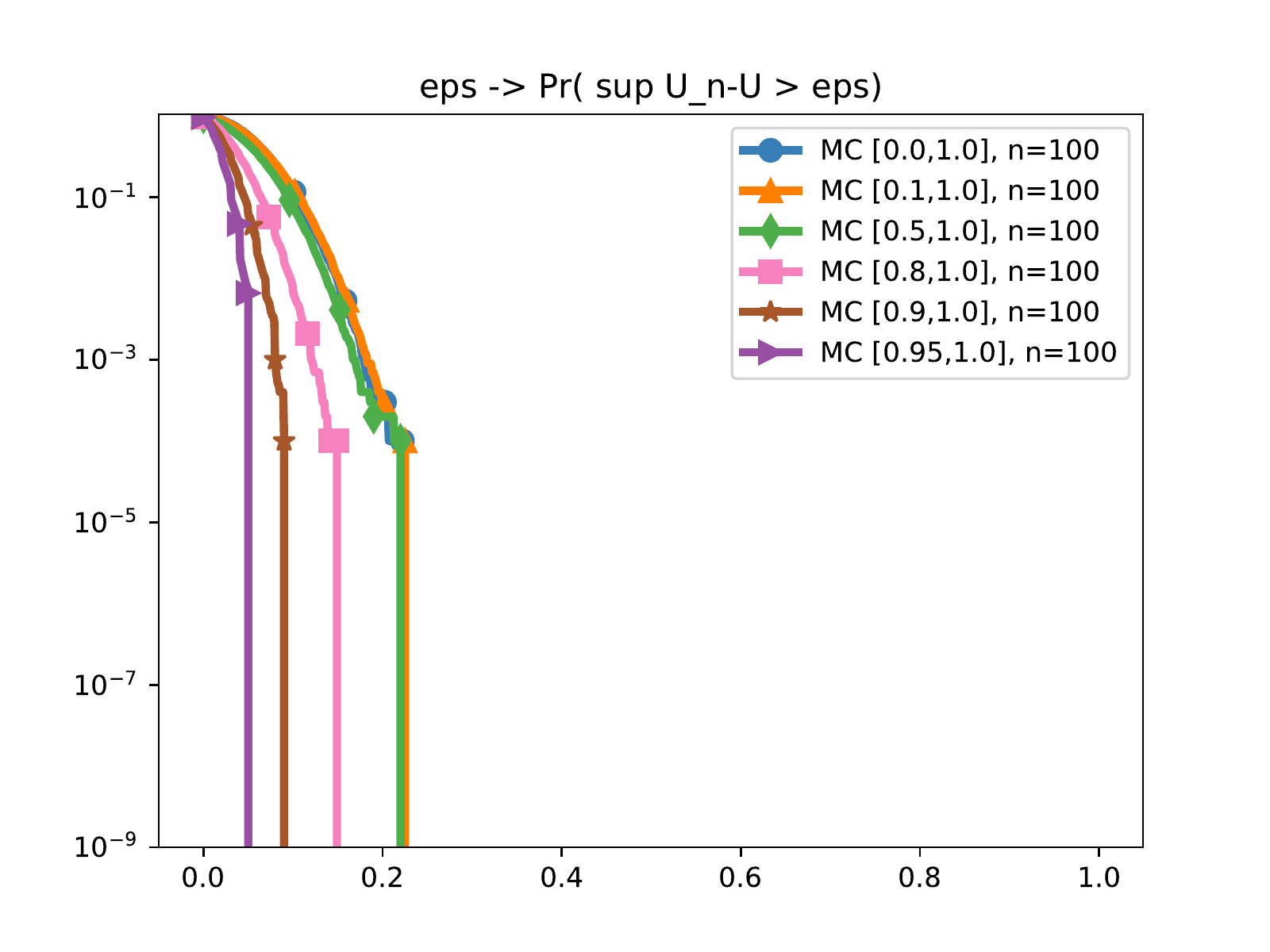}
	\includegraphics[width=0.43\linewidth, trim = {0.9cm 0.8cm 1cm 0.8cm}, clip]{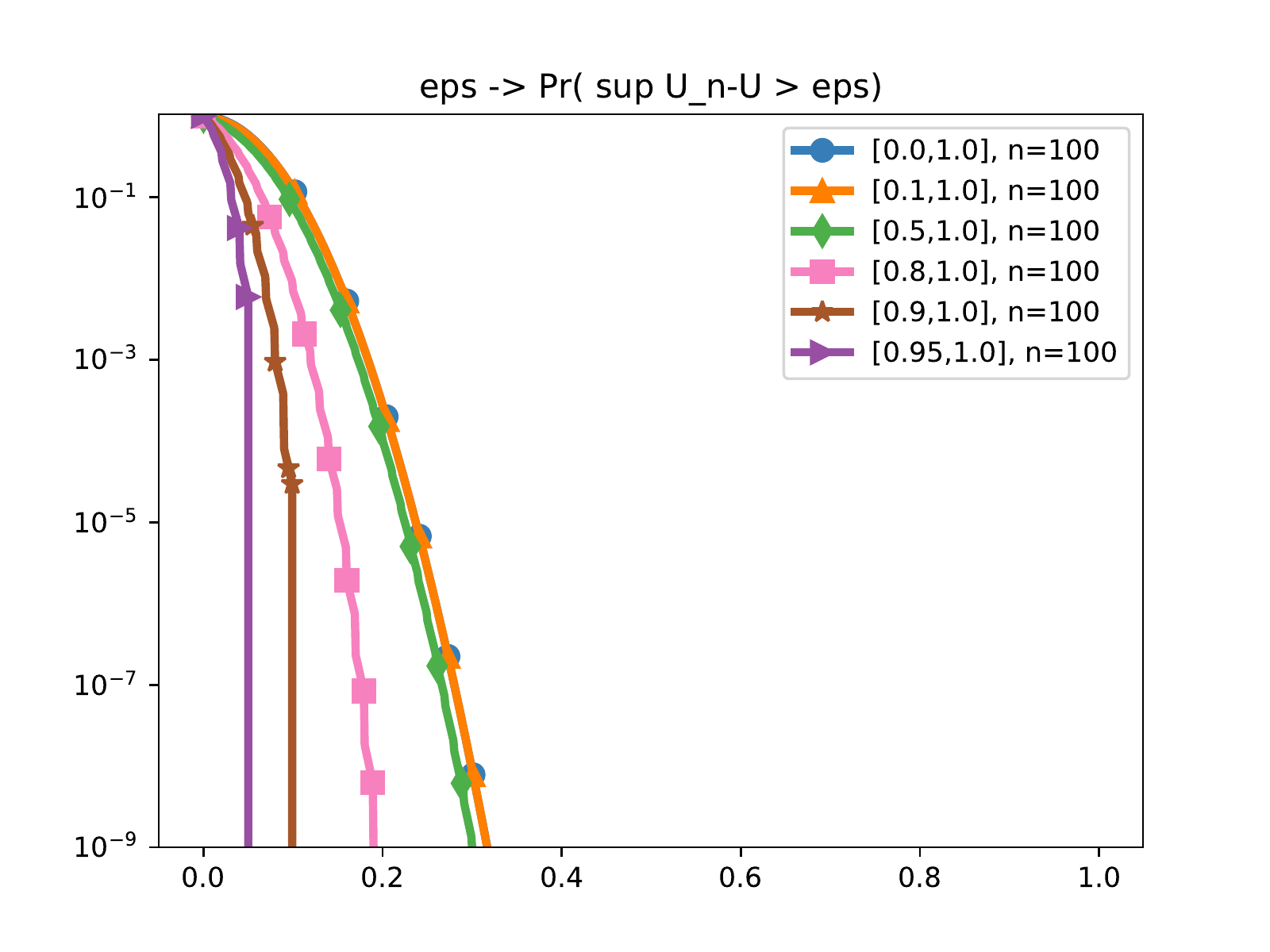}
	\vspace{-0.5cm}
	\caption{MCMC (left) versus Exact (Right) plot of $\epsilon\to\delta_{[\uu,\ou]}(n,\epsilon)$ for various values of $n$ and 
		interval $[\uu,\ou]$ build from $M=10^4$ replicate.}
	\label{fig:MCMC2}
\end{figure}

\section{Technical details regarding the CVAR}\label{app:cvar}

Proposition~\ref{prop:cvar} is a consequence of the following more general results.

\begin{proposition}[Conditional value at risk]\label{pro:CVARrewriting}
	Any solution $x^\star$ to the following problem
	\beqan
	\lCVAR_{1\!-\!\alpha}(\nu) &=& \inf_{x\in \Real}\bigg\{x + \frac{1}{\alpha}\Esp[\max(X\!-\!x,0)]\bigg\}
	\eeqan
	must satisfy $1-\alpha \in [F(x^\star)-\Pr(X=x^\star), F(x^\star)]$. Further, it holds
	\beqan
	\lCVAR_{1\!-\!\alpha}(\nu) = \frac{1}{\alpha}\bigg(\Esp\Big[X\indic{X> x^\star}\Big] + x^\star\Big(\!F(x^\star)-(1\!-\!\alpha\!)\Big)\bigg)\,.
	\eeqan
\end{proposition}

\begin{myproof}{of Proposition~ \ref{pro:CVARrewriting}}
	Let us introduce the function $H:x \mapsto x + \frac{1}{\alpha}\Esp[\max(X\!-\!x,0)]$.
	This is a convex function. Let $\partial H(x)$ denotes its subdifferential  at point $x$. 
	In particular, for $y \in \partial H(x)$, we must have 
	$\forall x', H(x')\leq H(x) + y(x'-x)$, and $x$ is a minimum of $H$ if $0\in \partial H(x)$.
	Using Minkowski set notations, we first have
	\beqan
	\partial H(x) = \{1\} + \frac{1}{\alpha }\partial\Esp[\max(X\!-\!x,0)],
	\eeqan
	hence we focus on computing $\partial\Esp[(X\!-\!x)\indic{X\!>\!x}]$.
	To this end, we look at the $y$ such that 
	\beqan
	\forall x', &&\Esp[(X\!-\!x')\indic{X\!>\!x'}] \geq \Esp[(X\!-\!x)\indic{X\!>\!x}] + y(x'-x)\\
	\text{i.e.} \forall x', && -(x-x')\Esp[\indic{X\!>\!x}]+
	\Esp[(X-x')(\indic{X\!>\!x'}- \indic{X\!>\!x})] \geq y(x'-x)
	\eeqan
	Remarking that if $x>x'$, then $\indic{X\!>\!x'}- \indic{X\!>\!x} = \indic{X\in(x',x]}$, while if $x'>x$ then $\indic{X\!>\!x'}- \indic{X\!>\!x} = -\indic{X\in(x,x']}$, and reorganizing the terms, this means we must have
	\beqan
	\forall x'>x,&& -\Esp[\indic{X\!>\!x}] - \Esp[\frac{(X\!-\!x')}{x'-x}\indic{X\in (x,x']}] \geq  y\\
	\forall x'<x, &&  -\Esp[\indic{X\!>\!x}] + \Esp[\frac{(X\!-\!x')}{x'-x}\indic{X\in (x',x]}] \leq  y\\
	\eeqan
	Further, note that if $x'>x$, then $\frac{(X\!-\!x')}{x'-x}\indic{X\in (x,x']} \in (-\indic{X\in (x,x']}, 0]$, while if $x>x'$, then 
	$\frac{(X\!-\!x')}{x'-x}\indic{X\in (x',x]} \in [-\indic{X\in (x',x]},0)$.
	Hence, we deduce that such $y$ must satisfy 
	\beqan 
	\inf_{x':x'>x}-\Pr(X\!>\!x)+\Pr(X\in(x,x'])\geq y \geq
	\sup_{x':x'<x} -\Pr(X\!>\!x)-\Pr(X\in(x',x])\,.\eeqan
	Hence,  $-\Pr(X\!>\!x)\geq y \geq  -\Pr(X\!\geq\!x)$,
	$\partial\Esp[(X\!-\!x)\indic{X\!>\!x}] \subset  [-\Pr(X\geq x), -\Pr(X> x)] = [F(x)-1- \Pr(X=x), F(x)-1]$, from which we deduce that 
	\beqan
	\partial H(x) \subset \bigg[ \frac{1}{\alpha }(F(x)-\Pr(X=x)-(1-\alpha)),  \frac{1}{\alpha }(F(x)-(1-\alpha))\bigg]
	\,.\eeqan
	
	This means that a minimum $x^\star$ of $H$ should at least satisfy that
	$1-\alpha \in [F(x^\star)-\Pr(X=x^\star), F(x^\star)]$.
	Finally, the value of the optimization is given by
	\beqan
	x^\star + \frac{1}{\alpha}\Esp[X\indic{X>x^\star}] - \frac{x^\star}{\alpha}(1\!-\!F(x^\star))
	&=&\frac{1}{\alpha}\Esp\Big[X\indic{X>x^\star}\Big] + \frac{x^\star}{\alpha}\Big(\!F(x^\star)-(1\!-\!\alpha\!)\Big)\,.
	\eeqan

\end{myproof}

\begin{proposition}[Expected shorfall]\label{pro:ESrewriting}
	Any solution $x^\star$ to the following problem
	\beqan
	\rCVAR_{\alpha}(\nu) &=& \sup_{x\in \Real}\bigg\{ \frac{1}{\alpha}\Esp[\min(X\!-\!x,0)]+x\bigg\}
	\eeqan
	must satisfy $\alpha \in [F(x^\star)-\Pr(X=x^\star),F(x^\star)]$. Further, it holds
	\beqan
	\rCVAR_{\alpha}(\nu) =  \frac{1}{\alpha}\bigg(\Esp\Big[X\indic{X<x^\star}\Big] + x^\star\Big(\! \alpha -F(x^\star)+\Pr(X=x^\star)\Big)\bigg)\,.
	\eeqan
\end{proposition}

\begin{myproof}{of Proposition~\ref{pro:ESrewriting}}
	Let us introduce the function $H:x \mapsto \frac{1}{\alpha}\Esp[\min(X\!-\!x,0)]+x$.
	This is a concave function. Let $\partial H(x)$ denotes its subdifferential  at point $x$. 
	In particular, for $y \in \partial H(x)$, we must have 
	$\forall x', H(x')\geq H(x) + y(x'-x)$, and $x$ is a minimum of $H$ if $0\in \partial H(x)$.	
	Using Minkowski set notations, we first have
	\beqan
	\partial H(x) =  \frac{1}{\alpha }\partial\Esp[\min(X\!-\!x,0)]+ \{1\},
	\eeqan
	hence we focus on computing $\partial\Esp[(X\!-\!x)\indic{X\!<\!x}]$.
	To this end, we look at the $y$ such that 
	\beqan
	\forall x', &&\Esp[(X\!-\!x')\indic{X\!<\!x'}] \leq \Esp[(X\!-\!x)\indic{X\!<\!x}] + y(x'-x)\\
	\text{i.e.} \forall x', && -(x-x')\Esp[\indic{X\!<\!x}]+
	\Esp[(X-x')(\indic{X\!<\!x'}- \indic{X\!<\!x})] \leq y(x'-x)
	\eeqan
	Remarking that if $x>x'$, then $\indic{X\!<\!x'}- \indic{X\!<\!x} = -\indic{X\in[x',x)}$, while if $x'>x$ then $\indic{X\!<\!x'}- \indic{X\!<\!x} = \indic{X\in[x,x')}$, and reorganizing the terms, this means we must have
	\beqan
	\forall x'>x,&& -\Esp[\indic{X\!<\!x}] + \Esp[\frac{(X\!-\!x')}{x'-x}\indic{X\in [x,x')}] \leq  y\\
	\forall x'<x, &&  -\Esp[\indic{X\!<\!x}] - \Esp[\frac{(X\!-\!x')}{x'-x}\indic{X\in [x',x)}] \geq  y\\
	\eeqan
	Further, note that if $x'>x$, then $\frac{(X\!-\!x')}{x'-x}\indic{X\in [x,x')} \in (-\indic{X\in [x,x')}, 0]$, while if $x>x'$, then 
	$\frac{(X\!-\!x')}{x'-x}\indic{X\in [x',x)} \in [-\indic{X\in [x',x)},0)$.
	Hence, we deduce that such $y$ must satisfy 
	\beqan 
	\sup_{x':x'>x}-\Pr(X\!<\!x)-\Pr(X\in[x,x'))\leq y \leq
	\inf_{x':x'<x} -\Pr(X\!<\!x)+\Pr(X\in[x',x))\,.\eeqan
	Hence,  $-\Pr(X\!\leq\!x)\leq y \leq  -\Pr(X\!<\!x)$,
	$\partial\Esp[(X\!-\!x)\indic{X\!<\!x}] \subset  [-\Pr(X\leq x), -\Pr(X< x)] = [-F(x), -F(x)+\Pr(X=x)]$, from which we deduce that 
	\beqan
	\partial H(x) \subset \bigg[ \frac{1}{\alpha }(-F(x)+ \alpha),  \frac{1}{\alpha }(-F(x) + \Pr(X=x)+\alpha)\bigg]
	\,.\eeqan
	
	This means that a minimum $x^\star$ of $H$ should at least satisfy that
	$\alpha \in [F(x^\star)-\Pr(X=x^\star), F(x^\star)]$.
	Finally, the value of the optimization is given by
	\beqan
	\lefteqn{
		x^\star + \frac{1}{\alpha}\Esp[X\indic{X<x^\star}] - \frac{x^\star}{\alpha}(\!F(x^\star)-\Pr(X=x^\star))}\\
	&=&\frac{1}{\alpha}\Esp\Big[X\indic{X<x^\star}\Big] + \frac{x^\star}{\alpha}\Big(\! \alpha -F(x^\star)+\Pr(X=x^\star)\Big)\,.
	\eeqan

\end{myproof}

\begin{proposition}[Integrated and Optimization forms]\label{pro:integrated}
	Let $X$ be a real-valued random variable with distribution $\nu$ and CDF $F$. Let $a,b\in\overline{\Real}$ be such that $\Pr_\nu(a\leq X\leq b)=1$.
	Let $\alpha\in[0,1]$ and $x^\star$ be any solution to the optimization problem
		$\lCVAR_{1\!-\!\alpha}(\nu)$. 
	Let $(\chi_i)_{i\in\mathbb{Z}}$ denotes the discontinuity points of $F$ (empty when $F$ is continuous),
	and let $\kappa=1\!-\!\alpha$. Then, if $a\geq 0$,  the following rewriting holds
	\beqan
	\lCVAR_{\kappa}(\nu) &=&  \frac{1}{\alpha}\int_{a}^b \bigg[\alpha-\max\big(F(x)\!-\!\kappa, F(x^\star)\!-\!\kappa\big)\bigg]dx+ \frac{a}{\alpha}\\
	&& + \frac{1}{\alpha}\sum_{i\in \mathbb{Z}}\Pr(X\!=\!\chi_i)\indic{\chi_i\!>\!x^\star}
	+\frac{x^\star}{\alpha}(F(x^\star)\!-\!\kappa)\,.
	\eeqan	
	In particular if $X$ is continuous, $a\geq 0$ and $b<\infty$, then
	\beqan
	\lCVAR_{\kappa}(\nu)&=&  (b\!-\!a) - \frac{1}{\alpha}\bigg(\int_{a}^b \max(F(x)\!-\!\kappa, 0)dx-a\bigg)\,.
	\eeqan
\end{proposition}

\begin{myproof}{of Proposition~\ref{pro:integrated}}
	Indeed, we first have that 
	\beqan
	\lCVAR_{\kappa}(\nu) =	\frac{1}{\alpha}\bigg(\Esp\Big[X\indic{X> x^\star}\Big] + x^\star\Big(\!F(x^\star)-(\kappa\!)\Big)\bigg)\,.
	\eeqan
	Now, if $a\geq 0$, then $Y=X\indic{X> x^\star}$ is a non-negative random variable, hence we can use the following rewriting $\Esp[Y] = \int_0^b\Pr(Y\geq y)dy = a+\int_a^b\Pr(Y\geq y)dy$.
	Hence,
	\beqan
	\Esp[Y] &=& a+\int_a^b \Pr\bigg( X\indic{X> x^\star} \geq x\bigg) dx\\
	&=& a+\int_a^b \Pr\bigg( X\geq x\bigg)\indic{x> x^\star} +\Pr\bigg( X> x^\star\bigg)\indic{x \leq  x^\star}  dx\\	
	&=& a+\int_a^b (1-F(x))\indic{x> x^\star} +(1-F(x^\star))\indic{x \leq  x^\star}  dx + \int_a^b\Pr(X=x)\indic{x> x^\star} \\
	&=& a+\int_a^b 1-F(x)\indic{x> x^\star}-F(x^\star)\indic{x \leq  x^\star}  dx + \sum_{i\in\mathbb{Z}}\Pr(X=x_i)\indic{x_i> x^\star} \\
	&=& a+\int_a^b 1-\max(F(x), F(x^\star)) dx + \sum_{i\in\mathbb{Z}}\Pr(X=x_i)\indic{x_i> x^\star}\,,
	\eeqan
	where the last line is by monotony of $F$. We conclude remarking that $1-\max(F(x), F(x^\star)) = \alpha-\max\big(F(x)\!-\!(\kappa), F(x^\star)\!-\!(\kappa)\big)$.
\end{myproof}

\begin{proposition}[Integrated and Optimization forms]\label{pro:integrated2}
	Let $X$ be a real-valued random variable with distribution $\nu$ and CDF $F$. Let $a,b\in\overline{\Real}$ be such that $\Pr_\nu(a\leq X\leq b)=1$.
	Let $\alpha\in[0,1]$ and $x^\star$ be any solution to the optimization problem
	$\rCVAR_{\alpha}(\nu)$. 
	Let $(\chi_i)_{i\in\mathbb{Z}}$ denotes the discontinuity points of $F$ (empty when $F$ is continuous). Then, if $a\geq 0$,  the following rewriting holds
	\beqan
	\rCVAR_{\alpha}(\nu)  &=&  \frac{1}{\alpha} \int_a^b (F(x^\star)-F(x))_{+}dx + \frac{a}{\alpha}\\
	&&+ \frac{1}{\alpha}\sum_{i\in\mathbb{Z}}(p_i-\Pr(X=x^\star))\indic{\chi_i <  x^\star} + \frac{x^\star}{\alpha}\Big(\! \alpha -F(x^\star)+\Pr(X=x^\star)\Big)\,.
	\eeqan	
	In particular if $X$ is continuous, $a\geq 0$ and $b<\infty$, then
	\beqan
	\rCVAR_{\alpha}(\nu) &=&   \frac{1}{\alpha} \int_a^b (\alpha-F(x))_{+}dx + \frac{a}{\alpha}\,.
	\eeqan
\end{proposition}

\begin{myproof}{of Proposition~\ref{pro:integrated2}}
	Indeed, we first have that 
	\beqan
	\rCVAR_{\alpha}(\nu) =  \frac{1}{\alpha}\bigg(\Esp\Big[X\indic{X<x^\star}\Big] + x^\star\Big(\! \alpha -F(x^\star)+\Pr(X=x^\star)\Big)\bigg)
	\eeqan
	Now, if $a\geq 0$, then $Y=X\indic{X< x^\star}$ is a non-negative random variable, hence we can use the following rewriting $\Esp[Y] = \int_0^b\Pr(Y\geq y)dy = a+\int_a^b\Pr(Y\geq y)dy$.
	Hence,
	\beqan
	\Esp[Y] &=& a+\int_a^b \Pr\bigg( X\indic{X< x^\star} \geq x\bigg) dx\\
	&=& a+0\times \indic{x\geq x^\star} +\int_a^b \Pr\bigg( x^\star> X\geq x\bigg)\indic{x <  x^\star}  dx\\	
	&=& a + \int_a^b (F(x^\star)-\Pr(X=x^\star)-F(x)+\Pr(X=x))\indic{x <  x^\star}  dx\\
	&=&a + \int_a^b (F(x^\star)-F(x))_{+}dx
	+\sum_{i\in\mathbb{Z}}(p_i-\Pr(X=x^\star))\indic{\chi_i <  x^\star}\,.
	\eeqan
	Hence, we deduce that
	\beqan
	\rCVAR_{\alpha}(\nu) &=&  \frac{1}{\alpha} \int_a^b (F(x^\star)-F(x))_{+}dx + \frac{a}{\alpha}\\
	&&+ \frac{1}{\alpha}\sum_{i\in\mathbb{Z}}(p_i-\Pr(X=x^\star))\indic{\chi_i <  x^\star} + \frac{x^\star}{\alpha}\Big(\! \alpha -F(x^\star)+\Pr(X=x^\star)\Big)\,.
	\eeqan
\end{myproof}

\section{Other result}\label{app:other}

We provide below  for the interested reader  some examples of functions $g$ satisfying  $\sum\limits_{t=1}^\infty \frac{1}{g(t)}\leq 1$.	
\begin{lemma}[Controlled sums]\label{lem:g}
	The following functions $g$ satisfy  $\sum\limits_{t=1}^\infty \frac{1}{g(t)}\leq 1$.	
	\begin{itemize}
		\item $g(t)=3t^{3/2}$, $g(t)=t(t+1)$, $g(t)=\frac{(t+1)\ln^2(t+1)}{\ln(2)}$,
		$g(t) = \frac{(t+2)\ln(t+2)(\ln\ln(t+2))^2}{\ln\ln(3)}$,
		\item For each $m\in\Nat$, $g_m(t)\! =\! C_m(\overline{\ln}^{\bigcirc\!\!\!\!\!{m}}(t))^2\prod_{i=0}^{m-1}\overline{\ln}^{\bigcirc\!\!\!\!{i}}(t)$,
		where $f^{\bigcirc\!\!\!\!\!{m}}$ denotes the $m$-fold composition of function $f$, $\overline{\ln}(t)\!=\! \max\{\ln(t),1\}$, and we introduced the constants 
		$C_1 = 2+\ln(2)+1/e$, $C_2= 2.03+\ln(e^e-1)$ as well as $C_m=2+\ln\Big( \exp^{\bigcirc\!\!\!\!\!{m}}(1) \Big)$ for $m\geq 3$.
	\end{itemize}   

\end{lemma}

\begin{myproof}{of Lemma~\ref{lem:g}}
	Note that
	$\overline{\ln}^{\bigcirc\!\!\!\!\!{m}}(t)\!=\!\ln^{\bigcirc\!\!\!\!\!{m}}(t)$ for $t\!\geq\! \exp^{\bigcirc\!\!\!\!\!{m}}(1)$ and $1$ else.	Using that $g(t)\! =\! C_m(\overline{\ln}^{\bigcirc\!\!\!\!\!{m}}(t))^2\prod_{i=0}^{m-1}\overline{\ln}^{\bigcirc\!\!\!\!{i}}(t)$,  and  that $t\!\mapsto\!-\!\frac{1}{\ln^{\bigcirc\!\!\!\!\!{m}}(t)}$ has derivative 
	$t\!\mapsto\!\frac{1}{(\ln^{\bigcirc\!\!\!\!\!{m}}(t))^2\prod_{i=0}^{m-1}\ln^{\bigcirc\!\!\!\!{i}}(t) }$, it comes
	\beqan
	\sum_{t=1}^\infty \frac{C_m}{g(t)} &=& \sum_{t=1}^{\lceil \exp^{\bigcirc\!\!\!\!\!{m}}(1) \rceil-1}\frac{1}{t} + \sum_{t=\lceil \exp^{\bigcirc\!\!\!\!\!{m}}(1) \rceil}\frac{1}{(\ln^{\bigcirc\!\!\!\!\!{m}}(t))^2\prod_{i=0}^{m-1}\ln^{\bigcirc\!\!\!\!{i}}(t) }\\
	&\leq&1\!+\!\ln\Big(\!\lceil \exp^{\bigcirc\!\!\!\!\!{m}}(1) \rceil\!-\!1\Big) +\\
	&&	\frac{1}{(\ln^{\bigcirc\!\!\!\!\!{m}}(\lceil \exp^{\bigcirc\!\!\!\!\!{m}}(1) \rceil))^2\prod_{i=0}^{m-1}\!\ln^{\bigcirc\!\!\!\!{i}}(\lceil \exp^{\bigcirc\!\!\!\!\!{m}}(1) \rceil) }
	\!+\!
	\frac{1}{\ln^{\bigcirc\!\!\!\!\!{m}}\!\Big(\!\lceil \exp^{\bigcirc\!\!\!\!\!{m}}(1)\! \rceil\Big)}\\
	&=& 2+ \ln\Big(\lceil \exp^{\bigcirc\!\!\!\!\!{m}}(1) \rceil-1\Big)
	+  \frac{1}{\prod_{i=0}^{m-1}\exp^{(m-i)}(1)}\,.
	\eeqan
\end{myproof}

\end{document}